\documentclass[11pt]{article}

\usepackage{amsmath,amsthm,amssymb}
\usepackage[T1]{fontenc}
\usepackage{newtxtext}
\usepackage{newtxmath}
\usepackage{mathrsfs}
\usepackage{bm}
\usepackage{graphicx}
\usepackage{booktabs}
\usepackage{algorithm}
\usepackage{algpseudocode}
\usepackage[margin=1in]{geometry}
\usepackage{natbib}
\bibpunct[, ]{(}{)}{,}{a}{}{,}
\usepackage[colorlinks=true,linkcolor=blue,citecolor=blue,urlcolor=blue]{hyperref}

\theoremstyle{plain}
\newtheorem{theorem}{Theorem}
\newtheorem{lemma}{Lemma}
\newtheorem{proposition}{Proposition}
\newtheorem{corollary}{Corollary}
\newtheorem{assumption}{Assumption}
\theoremstyle{remark}
\newtheorem{remark}{Remark}

\title{Optimal Batching with Abandonment: Threshold Policies Under Risk-Neutral and Risk-Sensitive Criteria}
\author{Shuren Wang\\
Rice University\\
\url{sw79@rice.edu}
\and
Guodong Pang\\
Rice University\\
\url{gdpang@rice.edu}}
\date{}

\begin{document}

\maketitle

\begin{abstract}
We study optimal batching in a discrete-time Markovian system with Bernoulli arrivals and independent item abandonment. Each batch earns a per-item processing reward and incurs holding, abandonment, and fixed stopping costs. We formulate the batching decision as an optimal stopping problem and characterize optimal threshold policies under risk-neutral, risk-seeking, and risk-averse criteria. Under risk neutrality and the entropic risk-seeking criterion, the value function admits a forward-recursion representation, and the resulting reduced objectives are strictly unimodal. In the risk-seeking case, the objective is strictly discretely concave, yielding an explicit marginal-cost condition. For risk aversion, the unregularized entropic continuation branch may diverge, precluding standard fixed-point arguments. We therefore study a regularized robust surrogate formulated as a minimax problem with a relative-entropy penalty and capped likelihood ratio. We establish the existence of an optimal threshold policy and derive an explicit truncated-Gibbs worst-case distribution. We also study repeated batch clearing under both the long-run average and discounted criteria, for which the threshold objectives remain strictly unimodal. 
Numerical experiments examine threshold sensitivity to system parameters, costs, and risk preferences.
\end{abstract}


\bigskip
\section{Introduction}

The coordination of batch processing and inventory accumulation is a fundamental operational challenge. In manufacturing, components arriving from upstream stations are consolidated into batches before transfer to downstream lines \citep{deb1973optimal}; in logistics, customer orders accumulate at distribution centers until a vehicle is dispatched \citep{ccetinkaya2002optimal}; in telecommunications, data packets are buffered until a transmission block is formed \citep{xiong2000control}; and in healthcare, platelet concentrates from multiple donations are pooled into transfusion doses, with production and inventory decisions coordinated to balance uncertain demand against outdating \citep{haijema2007blood}.
At the same time, when items accumulate in a system, they are typically subject to forces that reduce the inventory level over time. Customers waiting in a queue may abandon out of impatience \citep{baccelli1981queues, garnett2002designing}. Perishable goods in a cold chain may spoil \citep{graves1982application, kaspi1983inventory, nahmias2011perishable}.
Data packets in a router buffer may expire when their time-to-live counters reach zero. In each setting, processing items individually sacrifices economies of scale, whereas waiting to form larger batches increases holding costs
and the risk that items leave the system before processing.

We study a discrete-time batch accumulation problem that captures this trade-off. Items arrive stochastically, at most one per period, and each waiting item independently abandons before processing at a constant per-period rate. At the beginning of each period, the decision maker either processes the accumulated batch---paying a fixed stopping cost and earning a per-item reward---or continues accumulating, incurring holding costs on the waiting items and a penalty on each item that abandons. We analyze the resulting optimal stopping problem under three criteria: a discounted risk-neutral criterion, an undiscounted risk-seeking entropic criterion, and a robustly regularized risk-averse criterion.

Existing methods do not directly cover this problem, for two reasons. First, in the classical clearing and batch-service models cited below, items are not lost one by one before clearing or service; here, abandonment destroys that cumulative structure, and the usual monotonicity/submodularity sufficient conditions for threshold optimality  cannot be satisfied. 
Second, the risk-sensitive stopping results closest to our entropic criterion assume bounded running/terminal costs or comparable regularity, whereas here both grow linearly on an unbounded state space. 
Under risk aversion, the unregularized continuation branch lacks uniform control and the standard unique-fixed-point apparatus is unavailable. The literature review below substantiates both points.

Our main result is that threshold policies are optimal under all three criteria: there exists a critical inventory level at which the decision maker processes the batch and continues accumulating below it (Theorems~\ref{thm:RN}, \ref{thm:RS}, and \ref{thm:RA}). Beyond optimality, we characterize the threshold exactly in each regime. Under the risk-neutral criterion, the value of every threshold policy is given by an exact, strictly unimodal scalar objective, so the optimal threshold is identified by a discrete ascent that terminates after finitely many evaluations. Under the risk-seeking criterion, the reduced objective is strictly discretely concave, and the optimal threshold is marked by an explicit marginal-cost condition: the first level at which the marginal risk-sensitive cost of accumulation reaches the marginal processing reward. Under the risk-averse criterion, we replace the unregularized entropic continuation branch by a capped robust dual formulation---a minimax stochastic game in which an adversary perturbs the transition law subject to a relative-entropy penalty and a capped likelihood ratio---and show that it admits a unique fixed point, an optimal threshold policy, and a worst-case transition law of explicit truncated-Gibbs form.

A central contribution of this paper is a unified method for establishing threshold optimality on an unbounded state space. Rather than relying on state-space truncation or backward induction, our approach exploits the skip-free upward dynamics of the inventory process and the linear structure of the one-period cost. We characterize the value of each threshold policy through forward recursions and reduce the original stopping problem to a one-dimensional optimization. A coupling argument establishes that restricting attention to threshold policies entails no loss of optimality. Together, these results provide an exact characterization of the optimal policy while avoiding artificial restrictions on the inventory level.

The three regimes also exhibit a clear structural relationship: risk aversion accelerates processing relative to the risk-neutral benchmark, whereas risk seeking delays it, and both risk-sensitive thresholds converge to an optimal threshold of the undiscounted risk-neutral problem as the risk-sensitivity parameter vanishes, with the risk-averse limit selecting the smallest risk-neutral optimal threshold when multiple optima exist (Proposition~\ref{prop:regime-comparison}). The ordering further yields an exact truncation point for the numerical solution of the regularized risk-averse problem. 

Our numerical study confirms the comparative statics implied by the theory, including the finding that the response of the optimal threshold to the per-item reward can reverse sign depending on the discount factor. The risk-sensitive experiments quantify the risk-seeking premium and the risk-averse discount: risk seeking shifts the optimal threshold upward by deflating the marginal cost of accumulation, whereas regularized risk aversion shifts it downward by closing the stopping gap earlier. Both thresholds converge to the risk-neutral benchmark as the risk-sensitivity parameter approaches 0, consistent with the theoretical ordering and continuity results. The experiments also characterize the worst-case abandonment distortion under the likelihood-ratio cap and illustrate the operational cost of over-batching: as the buffer grows, abandonment increasingly offsets arrivals, cycle lengths grow disproportionately, and repeated-clearing performance deteriorates.

Finally, we return to repeated clearing: in the applications above, clearing is a recurring event rather than a terminal one. Section~\ref{sec:renewal} develops a repeated-clearing formulation in which each clearing resets the buffer to zero and the controller faces the same batching decision again, optimizing either the long-run average reward or the total discounted reward over infinitely many regenerative cycles. The first-passage decomposition extends with a single additional object---a time-accumulation function that prices the expected cycle length---and yields scalar reduced objectives for both repeated-clearing criteria. Threshold policies remain optimal among all stationary clearing policies. Moreover, both repeated-clearing objectives are strictly unimodal over the full threshold range: their increments admit a single sign change, characterized by a strictly decreasing marginal-rate function that prices the net reward of extending the cycle against its marginal length. The optimal threshold for repeated clearing is therefore computable exactly by a discrete ascent that terminates at the first nonpositive marginal improvement, as in the risk-neutral analysis---with no state-space truncation and no a priori bound on the search range. Furthermore, this threshold never exceeds the corresponding optimal-stopping threshold: the latter prices the erosion of value within a cycle but not the opportunity cost of the cycle length itself, and therefore systematically over-batches.

\subsection{Literature Review}

\noindent{\it Batch Clearing Systems.}
Operationally, our model is a batch clearing system without service: the controller decides when to clear the accumulated inventory, but clearing neither occupies a server nor delays future arrivals. \citet{stidham1974stochastic}
and \citet{whitt1981stationary} analyzed clearing processes in stochastic input systems, characterizing the stationary distribution of the accumulated quantity under exogenous clearing schemes. \citet{ross1969optimal, ross1971infinitesimal} studied optimal clearing of Poisson processes and established the optimality of threshold policies under linear holding costs, while \citet{kella2003stochastic} extended the analysis to clearing systems with both Brownian and compound Poisson input. These works share the feature that inventory is not lost item by item before clearing: the accumulated input is removed only at clearing epochs. 
The optimal stopping problems studied here, however, differ from this long-run average-cost formulation; the repeated-clearing extension of Section~\ref{sec:renewal} bridges the two, adopting the long-run average and discounted criteria of the clearing literature while retaining the item-by-item abandonment that distinguishes the present model.

\noindent{\it Batch Service and Dispatching in Queueing Systems.}
A related queueing literature studies batch service, where the controller repeatedly decides when to initiate service of a batch while arrivals continue during the service itself---in effect, a repeated optimal stopping problem in which each service completion renews the stopping decision.
\citet{deb1973optimal} established the optimality of threshold policies for batch-service queues, and \citet{deb1976optimal} extended this analysis to settings with switching costs. \citet{weiss1979computation} developed methods for computing optimal control limits for such queues, while \citet{xia2002optimal} analyzed the optimal control of parallel queues with batch service. \citet{vanleeuwen2017optimal} studied optimal batch dispatching in a tandem of two single-server Markovian queues with batch service at the first station. \citet{browne1992gated} studied gated service in infinite-server queues with Poisson arrivals and uniform service times, in which service begins for a batch whenever all jobs in the preceding batch have been completed and the queue is nonempty.
The closest antecedent is \citet{larranaga2015efficient}, who studied optimal batching in a single-server queue with abandonment under a long-run average-cost criterion. 
The renewal formulation of Section~\ref{sec:renewal} adopts the same long-run average criterion as \citet{larranaga2015efficient} and the classical batch-service models, but differs in that there is no service in our model (no server is occupied),
so the cycle length is governed by the accumulation-and-abandonment dynamics alone.
Related work also considers vehicle or shuttle dispatching between two stations with transportation times \citep{deb1978optimal}.

\noindent{\it Optimal Stopping of Discrete-Time Markov Chains.}
The classical theory of optimal stopping for discrete-time Markov chains with time-additive cost structures \citep{shiryaev2008optimal, peskir2006optimal} delivers the existence of optimal stopping rules through the theory of excessive functions and Snell envelopes, and it yields control-limit (threshold) policies when the underlying dynamic program exhibits submodularity or supermodularity---monotonicity-based sufficient conditions that, as noted above, do not apply directly to our cost structure. Moreover, although Snell-envelope theory exists on general state spaces under integrability conditions, generic backward induction on a countably infinite grid is not directly implementable: without additional structure, computation uses a finite truncation whose level must be justified or controlled \citep{puterman2014markov}. We avoid this issue by exploiting the skip-free dynamics and linear cost structure of the present problem, which establish the crucial value-function decomposition that leads to a strictly unimodal reduced objective and yields the optimal threshold exactly, without truncation or backward induction.

\noindent{\it Risk-Sensitive Optimal Stopping of Markov Chains.}
Risk-sensitive Markov decision problems, in which the expected additive cost is replaced by an exponential (multiplicative) criterion, equivalently an entropic criterion, are well studied \citep{howard1972risk, denardo1979optimal, borkar2002risk, ruszczynski2010risk, pichler2022risk, ren2025structural}. Risk-sensitive optimal stopping in particular has attracted growing attention: \citet{nagai2007stopping} analyzed stopping problems for multiplicative functionals, \citet{bauerle2018risk} treated risk-sensitive stopping of continuous-time Markov chains, and \citet{jelito2021risk} developed a general theory of entropic optimal stopping; a finite-horizon treatment of the associated stopping problems was developed by \citet{palczewski2010finite}. These results, however, are developed under boundedness or regularity conditions on the costs that do not cover our setting, in which both the running cost and the terminal reward grow linearly in the inventory level and are therefore unbounded on an unbounded state space. The obstacle is not merely technical: for risk-sensitive optimal stopping with an unbounded terminal cost, \citet{jelito2022risk} showed that the associated Bellman equation may admit multiple solutions, of which only the minimal and maximal ones carry a probabilistic interpretation.
In the risk-averse regime, our model circumvents this difficulty through the equivalent minimax stochastic game formulation of the entropic criterion, regularized via the Donsker--Varadhan variational formula \citep{donsker1975asymptotic, dupuis2011weak}: the stopping problem is recast as a minimax stochastic game in which an adversary perturbs the reference transition law subject to an entropy penalty or constraint \citep{james1992asymptotic, dai1996connections, petersen2000minimax}. The robust-MDP and distributionally robust literature is related but not equivalent. Classical rectangular robust MDPs typically assume finite state and action spaces or rely on bounded rewards and contraction properties \citep{nilim2005robust,iyengar2005robust,wiesemann2013robust}. Subsequent general-state and approximation results exist, but they do not provide the threshold characterization sought here. In our setting, the likelihood-ratio cap is not primarily a model-uncertainty set; it is a regularization that restores uniform integrability and a uniform Foster--Lyapunov drift estimate for the regularized operator.

\subsection{Organization of the Paper}
Section~\ref{sec:model} formulates the model and presents the main results: Section~\ref{sec:rn} treats the risk-neutral criterion, Section~\ref{sec:risk-sensitive} the risk-seeking criterion, and Section~\ref{sec:rs-averse} the regularized risk-averse criterion. Section~\ref{sec:regime-comparison} compares the optimal thresholds across regimes. Section~\ref{sec:renewal} develops the renewal formulation for repeated clearing under the long-run average and discounted repeated-clearing criteria, and Section~\ref{sec:numerical} presents the numerical study. Sections~\ref{sec:RN}--\ref{sec:RA} prove Theorems~\ref{thm:RN}--\ref{thm:RA}, with supporting lemmas and the proofs of the repeated-clearing results deferred to the Appendix, and Section~\ref{sec:conclusion} concludes.

\section{Model and Results}\label{sec:model}

\subsection{Problem Formulation}\label{sec:formulation}
Consider an infinite-horizon, discrete-time Markovian system $\{X_t\}_{t\ge0}$ taking values in $\mathcal{S}:= \mathbb{Z}_{+}$, with dynamics as follows.

\noindent{\it State dynamics.}
\begin{enumerate}
    \item[(i)] \textit{Arrivals.} An external arrival occurs with probability $q \in (0,1)$. Let $\{D_t\}_{t \ge 0}$ be a sequence of i.i.d.\ Bernoulli random variables with $\mathbb{P}(D_t = 1) = q$, representing arrivals.

    \item[(ii)] \textit{Abandonment.} Each of the items present independently survives the period with probability $1-p$, where $p \in (0,1)$. Conditional on $X_t = x$, the number of departures in period $t$ is $A_t(x) \sim \operatorname{Binomial}(x, p)$. Abandonment occurs independently of arrivals.
\end{enumerate}
Conditional on continuation, the state transition dynamics are
\begin{equation}\label{eq:dynamics}
X_{t+1} = X_t - A_t(X_t) + D_t, \qquad t = 0, 1, 2, \dots.
\end{equation}
$X_t$ is the inventory level at the beginning of period $t$.
We write $\mathcal{P}(x, \cdot)$ for the one-step transition kernel induced by \eqref{eq:dynamics}. Because $0 < p, q < 1$, this kernel assigns positive probability to every state in $\{0, 1, \dots, x+1\}$ and zero probability elsewhere. Throughout, $A_x$ denotes a generic $\operatorname{Binomial}(x,p)$ random variable and $D$ a generic $\operatorname{Bernoulli}(q)$ random variable, independent of $A_x$.

\noindent{\it Action space and admissible policies.} At the beginning of each period $t$, the decision maker chooses an action $a_t \in \mathcal{A} = \{\textsc{stop}, \textsc{continue}\}$. Choosing \textsc{stop} terminates the process and triggers a batch operation; choosing \textsc{continue} allows the system to evolve for another period. A stationary deterministic Markov policy is a mapping $\pi: \mathcal{S} \to \mathcal{A}$. Of particular interest is the class of stationary threshold policies indexed by $\theta \in \mathbb{Z}_{+}$:
\begin{equation}\label{eq:threshold-policy}
\pi_\theta(x) = \begin{cases} \textsc{continue}, & x < \theta, \\ \textsc{stop}, & x \ge \theta. \end{cases}
\end{equation}
Under $\pi_\theta$, the associated stopping time is the first-entrance time into $\{\theta, \theta+1, \dots\}$, namely $\tau_\theta := \inf\{t \ge 0 : X_t \ge \theta\}$. Restricting attention to this class entails no loss of optimality under any criterion analyzed below: Propositions~\ref{prop:threshold-wlog} and~\ref{prop:rs-threshold-wlog} establish this for the risk-neutral and risk-seeking problems, respectively, and Theorem~\ref{thm:RA}(ii) does so for the regularized risk-averse problem.

\noindent{\it Economic primitives.} Let $\delta \in (0,1]$ denote the per-period discount factor. The cost and reward structure consists of a fixed stopping cost $\alpha > 0$ incurred upon stopping; a reward coefficient $\beta > 0$ earned per item processed in the batch; a per-period holding cost $c > 0$ per item in the buffer; and an abandonment penalty $\gamma > 0$ per item departing before processing. Holding costs are charged on beginning-of-period inventory, whereas abandonments are realized at the end of the period. The abandonment penalty is therefore discounted by one additional factor of $\delta$ relative to the holding cost in \eqref{eq:objective}. We write $\mathcal{C} := c + \delta \gamma p > 0$ for the effective expected one-stage cost rate and $S(x) := \beta x - \alpha$ for the stopping (terminal) reward.

\noindent{\it The optimization problems.} For each $x \in \mathcal{S}$, we write $\mathbb{P}_x$ for the law of the chain \eqref{eq:dynamics} started at $X_0 = x$ and $\mathbb{E}_x$ for the associated expectation. For a stopping time $\tau$, define the total discounted reward along the trajectory by
\begin{equation}\label{eq:objective}
\mathcal{R}_{\delta}(\tau) := \delta^{\tau} (\beta X_{\tau} - \alpha) - \sum_{t=0}^{\tau-1} \delta^t \bigl(c X_t + \delta \gamma A_t(X_t)\bigr).
\end{equation}
The risk-neutral and entropic risk-sensitive criteria from the initial state $x \in \mathcal{S}$ are, respectively,
\begin{equation}\label{eq:obj-rn}
J_{rn}(x, \tau) := \mathbb{E}_x\!\left[\mathcal{R}_{\delta}(\tau) \right],
\end{equation}
\begin{equation}\label{eq:obj-rs}
J_r(x, \tau) := \frac{1}{\eta} \log \mathbb{E}_x\!\left[e^{\eta \mathcal{R}_{\delta}(\tau)} \right], \qquad \eta \in \mathbb{R} \setminus \{0\}.
\end{equation}
The sign of $\eta$ determines the risk posture: $\eta > 0$ is risk seeking, whereas $\eta < 0$ is risk averse. The corresponding optimization problems are
\begin{equation}\label{eq:opt-problem-rn}
V_{rn}(x) := \sup_{\tau \in \mathfrak{T}} J_{rn}(x, \tau), \quad x \in \mathcal{S},
\end{equation}
\begin{equation}\label{eq:opt-problem-rs}
V_r(x) := \sup_{\tau \in \mathfrak{T}} J_r(x, \tau), \qquad x \in \mathcal{S},
\end{equation}
where $\mathfrak{T}$ denotes the set of almost surely finite stopping times adapted to the natural filtration of $\{X_t\}$.

\subsubsection{The Skip-Free Property}

The chain $\{X_t\}_{t\ge 0}$ is skip-free to the right, so every upward first
passage traverses all intermediate levels.

\begin{lemma}\label{lem:skip-free}
The DTMC $\{X_t\}$ governed by \eqref{eq:dynamics} is skip-free to the right: for
every $x \in \mathcal{S}$ and every $t \ge 0$,
\[
\mathbb{P}(X_{t+1} \le x + 1 \mid X_t = x) = 1,
\]
\[
\mathbb{P}(X_{t+1} = x+1 \mid X_t = x) = (1-p)^x q > 0.
\]
Moreover, every state $y \in \{0, 1, \dots, x+1\}$ is reachable from $x$ in one
step.
\end{lemma}

Two consequences are used repeatedly. First, when started from any state $x \le \theta$, the sample path visits every integer in $\{x, x+1, \dots, \sup_t X_t\}$; hence, under $\mathbb{P}_x$ with $x \le \theta$, the first-entrance time $\tau_\theta$ coincides with the first-passage time $T_\theta := \inf\{t \ge 0 : X_t = \theta\}$, and on $\{T_\theta < \infty\}$ one has $X_{T_\theta} = \theta$ and $T_\theta \ge \theta - x$ almost surely. 
Second, $T_\theta$ is almost surely finite: from every state $y < \theta$ the one-step probability of moving up one level is at least $(1-p)^{\theta-1} q > 0$, so sufficiently many consecutive upward moves occur eventually with probability one.

\subsection{Risk-Neutral Optimal Batching}\label{sec:rn}

\subsubsection{Bellman Optimality Equation}\label{sec:rn-bellman}

For $\rho > 1$, let $\mathcal{B}_\rho$ denote the Banach space of functions $u: \mathcal{S} \to \mathbb{R}$ with finite weighted supremum norm $\|u\|_\rho := \sup_{x \in \mathcal{S}} \rho^{-x} |u(x)|$. Define the Bellman operator $\mathcal{T}_{rn}$ by
\begin{equation}\label{eq:bellman-operator}
\mathcal{T}_{rn} (u)(x) = \max\!\Bigl\{ \,\underbrace{\beta x - \alpha}_{\textsc{stop}},\; \underbrace{- \mathcal{C}\, x + \delta \, \mathbb{E}[u(\hat{X}_x)]}_{\textsc{continue}} \Bigr\},
\end{equation}
where $\hat{X}_x := x - A_x + D$ is the next state given that the current state is $x$. The one-stage cost in the continuation branch is the expected running cost in \eqref{eq:objective}: since $\mathbb{E}[A_x] = px$, conditioning gives $\mathbb{E}_x[\sum_{t<\tau} \delta^t (c X_t + \delta \gamma A_t(X_t))]= \mathcal{C}\, \mathbb{E}_x[\sum_{t<\tau} \delta^t X_t]$ for every stopping time $\tau$ (the tower-property computation is carried out in Section~\ref{app:RN-main}), and the case $\tau \equiv 1$ yields the continuation branch of \eqref{eq:bellman-operator}. The \textit{Bellman equation} associated with problem \eqref{eq:opt-problem-rn} is the fixed-point equation
$V_{rn} = \mathcal{T}_{rn}V_{rn}$, that is,
\begin{equation}\label{eq:bellman-equation}
V_{rn}(x) = \max\!\Bigl\{ \,{\beta x - \alpha},\; {- \mathcal{C}\, x + \delta \, \mathbb{E}_{A_x, D}\!\bigl[V_{rn}(\hat{X}_x)\bigr]}\Bigr\}.
\end{equation}

\begin{proposition}[Well-posedness of the Bellman equation]\label{prop:rn-bellman}
\begin{enumerate}
    \item[(i)] Let $\delta < 1$ and $\rho \in \bigl(1,\, 1 + \frac{1-\delta}{\delta q}\bigr)$. Then $\mathcal{T}_{rn}$ maps $\mathcal{B}_\rho$ into itself and is a contraction with modulus $\delta(1 - q + q\rho) < 1$. The value function $V_{rn}$ of \eqref{eq:opt-problem-rn} belongs to $\mathcal{B}_\rho$ and is the unique solution of the Bellman equation~\eqref{eq:bellman-equation} in $\mathcal{B}_\rho$; value iteration $u_{n+1} = \mathcal{T}_{rn} u_n$ converges to $V_{rn}$ from any $u_0 \in \mathcal{B}_\rho$.
    \item[(ii)] Let $\delta = 1$. Then $V_{rn}$ is the minimal fixed point of $\mathcal{T}_{rn}$ among functions $u: \mathcal{S} \to \mathbb{R}$ dominating the stopping reward $x \mapsto \beta x - \alpha$, and the iterates defined by $u_0(x) = \beta x - \alpha$ and $u_{n+1} = \mathcal{T}_{rn} u_n$ increase pointwise to $V_{rn}$.
\end{enumerate}
\end{proposition}

Direct solution of \eqref{eq:bellman-equation} is hindered by a structural coupling: the continuation value at state $x$ depends on the values at $\{0, 1, \dots, x+1\}$, which in turn depend on the unknown optimal stopping boundary. We therefore first reduce the search to the threshold family \eqref{eq:threshold-policy} (Section~\ref{sec:rn-increments}) and then evaluate any fixed threshold by a forward recursion (Section~\ref{sec:decoupling}). Proofs for this subsection and Section~\ref{sec:rn-increments} are deferred to Section~\ref{app:RN-main}.

\subsubsection{Increment Bound and Threshold Structure}\label{sec:rn-increments}

The reduction to threshold policies rests on a uniform bound on the increments of the value function.

\begin{lemma}[Increments of the value function]\label{lem:increments}
For all $x \in \mathcal{S}$, $V_{rn}(x+1) - V_{rn}(x) \le \beta$.
\end{lemma}

\begin{proposition}[Threshold structure of optimal policies]\label{prop:threshold-wlog}
The optimal stopping region $\mathcal{S}^* := \{x \in \mathcal{S} : V_{rn}(x) = \beta x - \alpha\}$ is up-closed: $x \in \mathcal{S}^*$ implies $x+1 \in \mathcal{S}^*$. Consequently, whenever the supremum in \eqref{eq:opt-problem-rn} is attained by the first-entrance rule into $\mathcal{S}^*$, that rule is the stationary threshold policy $\pi_{\theta^*}$ with $\theta^* = \min \mathcal{S}^*$.
\end{proposition}

Proposition~\ref{prop:threshold-wlog} reduces problem \eqref{eq:opt-problem-rn} to a one-dimensional search over the integer threshold $\theta$. The remainder of this section makes the search efficient: Section~\ref{sec:decoupling} derives a reduced objective that is computable by a forward recursion, and Section~\ref{sec:rn-threshold} shows that this objective is strictly unimodal with an attained maximizer, so the search terminates at the first local maximizer.

\subsubsection{Value-Function Decomposition}\label{sec:decoupling}

Fix a threshold $\theta \in \mathbb{Z}_{+}$. Because the one-stage cost is linear in the state and, by Lemma~\ref{lem:skip-free}, the first passage to $\theta$ traverses every intermediate level, the value function under $\pi_\theta$ decomposes additively along the trajectory. The decomposition is anchored at state $0$ and captured by two auxiliary sequences. Write
\[
\omega_x := 1-\delta(1-p)^x(1-q), \qquad x \ge 0.
\]
Define $f: \mathcal{S} \to \mathbb{R}$, the \textit{cost-accumulation function}, by $f(0) = 0$ and the forward recurrence
\begin{equation}\label{eq:f-recurrence}
f(x+1)
= \frac{\mathcal{C}x+\omega_x f(x)-\delta F(x)}
       {\delta(1-p)^xq},
\quad x \ge 0,
\end{equation}
where
\begin{equation}\label{eq:F-operator}
F(x) = \mathbb{E}_{A_x, D}\!\bigl[\mathbf{1}_{A_x \ge 1} \cdot f(\hat{X}_x)\bigr]
\end{equation}
aggregates the contribution of periods with strictly positive attrition. Define $g: \mathcal{S} \to \mathbb{R}$, the \textit{discount-accumulation function}, by $g(0) = 1$ and
\begin{equation}\label{eq:g-recurrence}
g(x+1) = \frac{\omega_x g(x)-\delta G(x)}
       {\delta(1-p)^xq}, \quad x \ge 0,
\end{equation}
where
\begin{equation}\label{eq:G-operator}
G(x) = \mathbb{E}_{A_x, D}\!\bigl[\mathbf{1}_{A_x \ge 1} \cdot g(\hat{X}_x)\bigr].
\end{equation}
Both are proper forward recursions: on $\{A_x \ge 1\}$ one has $\hat{X}_x \le x$, so $F(x)$ and $G(x)$ depend only on values of $f$ and $g$ at states $0, \dots, x$.

\begin{remark}[Fixed-point form of the recurrences]\label{rem:f-rearrange}
The recurrences \eqref{eq:f-recurrence} and~\eqref{eq:g-recurrence} are equivalent
to
\begin{equation}\label{eq:f-fixed}
f(x) = -\mathcal{C}x+\delta\,\mathbb{E}[f(\hat{X}_x)],
\end{equation}
\begin{equation}\label{eq:g-fixed}
g(x) = \delta\,\mathbb{E}[g(\hat{X}_x)].
\end{equation}
Indeed, conditioning on $A_x$ and solving for $f(x+1)$ and $g(x+1)$ gives the
original recurrences. Binomial thinning also implies that, for every
$u:\mathbb{Z}_{+}\to\mathbb{R}$ and $k\ge1$,
\begin{equation}\label{eq:binomial-difference}
\Delta^k\,\mathbb{E}[u(Y_x)]
=(1-p)^k\,\mathbb{E}[\Delta^k u(Y_x)],
\end{equation}
where $\Delta u(x) := u(x+1) - u(x)$ is the forward-difference operator,
$\Delta^k$ its $k$-th iterate, and
$Y_x := x - A_x \sim \operatorname{Binomial}(x, 1-p)$.
\end{remark}

\begin{proposition}\label{prop:value-decomposition-and-objective}
Fix $\theta \in \mathbb{Z}_{+}$, let $V_{rn,\pi_\theta}(x) := J_{rn}(x, \tau_\theta)$ denote the value function under $\pi_\theta$ started from $x$, and write $v_{rn}(\theta) := V_{rn,\pi_\theta}(0)$ for the \textit{reduced objective}.
\begin{enumerate}
    \item[(i)] (Affine decomposition) For every initial state $x \in \mathcal{S}$ and every threshold $\theta \ge x$,
    \begin{equation}\label{eq:rn-affine-in-v}
    V_{rn,\pi_\theta}(x) = g(x)\, v_{rn}(\theta) + f(x),
    \end{equation}
    while over thresholds $\theta \le x$ the policy stops immediately and $V_{rn,\pi_\theta}(x) = \beta x - \alpha$ is independent of $\theta$. For fixed $x$, the right-hand side of \eqref{eq:rn-affine-in-v} is a strictly increasing affine transformation of the reduced objective, so every maximizer $\theta^* \ge x$ of $v_{rn}$ also maximizes $\theta \mapsto V_{rn,\pi_\theta}(x)$ over all $\theta \in \mathbb{Z}_{+}$.
    \item[(ii)] (Discounted reduced objective) For $\delta \in (0, 1)$, $g(x) \ge 1$ and
    \begin{equation}\label{eq:reduced-objective}
    v_{rn}(\theta) \;=\; \frac{\beta \theta - \alpha - f(\theta)}{g(\theta)}.
    \end{equation}
    \item[(iii)] (Undiscounted reduced objective) For $\delta = 1$, $g(x) \equiv 1$ and
    \begin{equation}\label{eq:undiscounted-reduced-objective}
    v_{rn}(\theta) \;=\; \beta \theta - \alpha - f(\theta).
    \end{equation}
\end{enumerate}
\end{proposition}

Maximizing the reduced objective anchored at the reference state $0$ therefore yields, in a single one-dimensional search, a threshold policy that is optimal among threshold policies from every initial state simultaneously.

\begin{remark}[Operational interpretations of $f$ and $g$]\label{rem:fg-meaning}
Let
\[
C_\theta := \mathbb{E}_0\!\left[\sum_{t=0}^{T_\theta-1} \delta^t \bigl(c X_t + \delta\gamma A_t(X_t)\bigr)\right]
\]
denote the expected present value of the holding and abandonment costs accrued during the first passage from the reference state $0$ to $\theta$. Directly from
\eqref{eq:objective}, using $X_{T_\theta} = \theta$,
\[
V_{rn,\pi_\theta}(0) = (\beta\theta - \alpha)\, \mathbb{E}_0\!\left[\delta^{T_\theta}\right] - C_\theta.
\]
Because $f$ and $g$ do not depend on $\alpha$ or $\beta$, matching the coefficient of $\beta\theta - \alpha$ and the constant term with those in \eqref{eq:reduced-objective}
identifies
\[
\mathbb{E}_0\!\left[\delta^{T_\theta}\right] = \frac{1}{g(\theta)}, \qquad C_\theta = \frac{f(\theta)}{g(\theta)}.
\]
Thus $g(\theta)$ is the reciprocal of the expected discount factor over the first passage, and $f(\theta)$ is the expected present value of the first-passage costs expressed in units of $\mathbb{E}_0[\delta^{T_\theta}]$. Only when $\delta = 1$ do $g(\theta)$ and $f(\theta)$ coincide with $\mathbb{E}_0[\delta^{T_\theta}]$ and $C_\theta$ themselves.
\end{remark}

By Proposition~\ref{prop:threshold-wlog} and Proposition~\ref{prop:value-decomposition-and-objective}(i), the search for an optimal policy from any initial state reduces to maximizing the reduced objective
$v_{rn}$, whose evaluation requires only the forward iterates of \eqref{eq:f-recurrence}--\eqref{eq:g-recurrence}; the remainder of this section analyzes $v_{rn}$ alone.

\subsubsection{The Optimal Threshold}\label{sec:rn-threshold}

To rule out the degenerate regime in which no batch is ever profitable, we impose a standing assumption.

\begin{assumption}[Profitability]\label{ass:profitability}
There exists at least one $\theta \in \mathbb{Z}_{+}$ such that
\begin{equation}\label{eq:profitability}
N_{rn}(\theta) \;:=\; \beta\theta - \alpha - f(\theta) \;>\; 0.
\end{equation}
\end{assumption}

Write $\Theta_+ := \{\theta \in \mathbb{Z}_{+} : N_{rn}(\theta) > 0\}$ for the positive domain of the numerator and $\overline\theta := \max\Theta_+$ for its largest element whenever $\Theta_+$ is finite; Assumption~\ref{ass:profitability} is imposed throughout the remainder of this subsection.

\begin{theorem}[Optimal risk-neutral policy]\label{thm:RN}
Suppose that Assumption~\ref{ass:profitability} holds.
\begin{enumerate}
     \item[(i)] (Strict unimodality) On its positive domain $\Theta_+$, the reduced objective is strictly unimodal: there exists $\theta^u \in \Theta_+$ such that $\Delta v_{rn}(\theta) > 0$ for $\theta < \theta^u$, $\Delta v_{rn}(\theta^u) \le 0$, and $\Delta v_{rn}(\theta) < 0$ for $\theta^u < \theta \le \overline\theta$.
    \item[(ii)] (Tail behavior) For $\delta \in (0, 1)$, define
    \begin{equation}\label{eq:c-infty}
        C_\infty := \frac{\delta q\, \mathcal{C}}{(1-\delta)\bigl(1-\delta(1-p)\bigr)}.
    \end{equation}
    Then
    \begin{equation}\label{eq:rn-tail}
        \lim_{\theta \to \infty} v_{rn}(\theta) =
        \begin{cases}
            -C_\infty, & \text{if } \delta < 1, \\[6pt]
            -\infty,   & \text{if } \delta = 1.
        \end{cases}
    \end{equation}
    \item[(iii)] (Optimal threshold) The index
    \begin{equation}\label{eq:optimal-theta}
    \theta^*_{rn} = \min \bigl\{\theta \in \mathbb{Z}_{+} : \Delta v_{rn}(\theta) \le 0 \bigr\}
    \end{equation}
    is well defined and finite, coincides with $\theta^u$, and is the smallest optimal threshold: at most two thresholds are optimal, and they differ by one. The threshold policy $\pi_{\theta^*_{rn}}$ is optimal among all admissible policies from every initial state: its value function $V_{rn,\pi_{\theta^*_{rn}}}$ solves the Bellman equation~\eqref{eq:bellman-equation} and coincides with $V_{rn}$, and the optimal stopping region is $\mathcal{S}^* = \{\theta^*_{rn}, \theta^*_{rn}+1, \dots\}$. Moreover,
    \begin{equation}\label{eq:rn-attainment}
    V_{rn}(0) \;=\; v_{rn}(\theta^*_{rn}) \;=\; \max_{\theta \in \mathbb{Z}_{+}} v_{rn}(\theta).
    \end{equation}
\end{enumerate}
\end{theorem}

The proof is given in Section~\ref{sec:RN}. Theorem~\ref{thm:RN} justifies a discrete ascent---evaluate $v_{rn}(0), v_{rn}(1), \dots$ and stop at the first index with $\Delta v_{rn}(\theta) \le 0$---that terminates after finitely many evaluations at the smallest optimal threshold.

\begin{remark}[The degenerate case]\label{rem:degenerate}
If \eqref{eq:profitability} fails, then no batch recoups the fixed cost of processing it: $N_{rn}(\theta) \le 0$ for every $\theta$, the positive domain $\Theta_+$ is empty, and the unimodality statement of Theorem~\ref{thm:RN}(i) is vacuous. In this regime, operating the system is never profitable, and the interesting question concerns whether, rather than when, to clear the buffer. We exclude it throughout: the comparison and regularization results of Sections~\ref{sec:rs-averse} and~\ref{sec:regime-comparison} that involve $\theta^*_{rn}$---in particular the truncation certificate $X_{\max} = \theta^*_{rn}$ of Remark~\ref{rem:ra-numerics}---are all stated under Assumption~\ref{ass:profitability}.
\end{remark}

\subsection{Risk-Seeking Optimal Batching}\label{sec:risk-sensitive}

We now turn to the entropic criterion \eqref{eq:obj-rs}.

\begin{assumption}\label{ass:undiscounted}
The risk-sensitive analysis is restricted to the undiscounted case $\delta = 1$.
\end{assumption}

\begin{remark}\label{rem:why-delta-1}
Under discounting, stationary multiplicative separability is lost: for $\delta<1$, the certainty equivalent still satisfies a Bellman recursion, but the risk parameter is discounted together with the reward, $\eta \mapsto \eta\delta$, so the value function is defined on the augmented state $(x,\lambda)$ and the optimal stopping region is a boundary in the $(x,\lambda)$-plane rather than a fixed threshold \citep{chung1987discounted}. The decomposition of Proposition~\ref{prop:risk-V0} requires the undiscounted setting, in which the same $\eta$ appears at every step; the discounted case can be solved numerically (Section~\ref{app:discounted-rs}).
\end{remark}

\subsubsection{Bellman Equation}\label{sec:rs-criterion}
Throughout the risk-seeking analysis (Sections~\ref{sec:rs-criterion}--\ref{sec:rs-seeking}), we write $V_{rs} := V_r$ for the optimal entropic value~\eqref{eq:opt-problem-rs} with $\eta > 0$. Let
\begin{equation}\label{eq:rs-continuation}
Q_r(x)
:=-cx+\frac{1}{\eta}\log
\mathbb{E}\!\left[e^{\eta(V_r(\hat{X}_x)-\gamma A_x)}\right]
\end{equation}
denote the continuation branch. By multiplicative separability, the associated
value function satisfies the Bellman equation
\begin{equation}\label{eq:bellman-risk-primal}
V_r(x)=\max\!\left\{S(x),\,Q_r(x)\right\}.
\end{equation}
Well-posedness of \eqref{eq:bellman-risk-primal} for $\eta > 0$ rests on a pathwise envelope for the undiscounted reward, which exists because holding costs accumulate superlinearly along any excursion while the terminal reward grows only linearly.

\begin{proposition}[Well-posedness of the risk-seeking Bellman equation]\label{prop:rs-dpp}
Let $\eta > 0$ and define
\[
\bar U(x) := \frac{(\beta + c/2)^2}{2c} - \alpha + \frac{c}{2}\,x(x-1), \qquad x \in \mathbb{Z}_{+}.
\]
\begin{enumerate}
    \item[(i)] (Pathwise envelope) For every \(x\) and every almost surely finite stopping time \(\tau\),
    \[
    \mathcal{R}_1(\tau) \le \bar U(x)
    \qquad \mathbb{P}_x\text{-almost surely}.
    \]    
    Consequently, \(\beta x-\alpha\le V_{rs}(x)\le\bar U(x)\) for all \(x\). Moreover, along the chain under any policy, for every \(t\ge0\), \(\mathbb{P}_x\)-almost surely, 
    \[
    V_{rs}(X_t)
    -\sum_{s=0}^{t-1}\bigl(cX_s+\gamma A_s(X_s)\bigr)
    \le \bar U(x).
    \]
    \item[(ii)] (Dynamic programming) $V_{rs}$ satisfies the Bellman equation~\eqref{eq:bellman-risk-primal}, and all expectations therein are finite.
    \item[(iii)] (Minimality and value iteration) Let $\mathcal{T}_{rs}$ denote the Bellman operator defined by the right-hand side of \eqref{eq:bellman-risk-primal}. The iterates $u_0(x) = \beta x - \alpha$ and $u_{n+1} = \mathcal{T}_{rs} u_n$ increase pointwise to $V_{rs}$. Moreover, $V_{rs}$ is the minimal solution of \eqref{eq:bellman-risk-primal} dominating the stopping reward: every solution $\widehat V$ with $\widehat V(x) \ge \beta x - \alpha$ for all $x$ satisfies $\widehat V \ge V_{rs}$.
\end{enumerate}
\end{proposition}

Proofs of the propositions and lemmas in Sections~\ref{sec:rs-criterion},~\ref{sec:rs-decoupling}, and~\ref{sec:rs-seeking} are deferred to Section~\ref{app:RS}.

\subsubsection{Risk-Seeking Value-Function Decomposition}\label{sec:rs-decoupling}

We extend the decomposition method of Section~\ref{sec:decoupling} to the risk-seeking case ($\eta > 0$). Write
\[
\omega_{rs}(x)
:=e^{\eta c x}-(1-p)^x(1-q), \qquad x\ge0.
\]
Define $f_{rs}:\mathbb{Z}_{+}\to\mathbb{R}$, the \textit{risk-sensitive cost-accumulation function}, by $f_{rs}(0)=0$ and for $x \ge 0$
\begin{equation}\label{eq:fr-recurrence}
f_{rs}(x+1)
=\frac{1}{\eta}\log\!\left(
\frac{e^{\eta f_{rs}(x)}\omega_{rs}(x)-F_{rs}(x)}{(1-p)^xq}
\right),
\end{equation}
where the \textit{risk-adjusted attrition operator} $F_{rs}: \mathbb{Z}_{+} \to \mathbb{R}$ is
\begin{equation}\label{eq:Fr-operator}
F_{rs}(x) = \mathbb{E}_{A_x, D}\!\left[ \mathbf{1}_{A_x \ge 1} \cdot e^{\eta (f_{rs}(\hat{X}_x) - \gamma A_x)}\right].
\end{equation}
The recurrence \eqref{eq:fr-recurrence} is the risk-sensitive analogue of \eqref{eq:f-recurrence}: the factor $e^{\eta c x}$ captures the certain holding cost, while $F_{rs}(x)$ aggregates the uncertain attrition-adjusted continuation. It is a proper forward recursion by the same argument as in Section~\ref{sec:decoupling}, and Proposition~\ref{prop:rs-forward}(i) shows that the argument of the logarithm is strictly positive at every state, so the recursion is well posed.

\begin{proposition}[Risk-seeking value decomposition]\label{prop:risk-V0}
Fix $\theta \in \mathbb{Z}_{+}$, let $V_{rs,\pi_\theta}(x) := J_r(x, \tau_\theta)$ denote the risk-seeking value function under $\pi_\theta$ started from $x$, and write $v_{rs}(\theta) := V_{rs,\pi_\theta}(0)$ for the \textit{risk-seeking reduced objective}.
\begin{enumerate}
    \item[(i)] (Additive decomposition) For every initial state $x \in \mathcal{S}$ and every threshold $\theta > x$,
    \begin{equation}\label{eq:risk-V-decomp}
    V_{rs,\pi_\theta}(x) = f_{rs}(x) + v_{rs}(\theta),
    \end{equation}
    while over thresholds $\theta \le x$, the policy stops immediately and $V_{rs,\pi_\theta}(x) = \beta x - \alpha$ is independent of $\theta$. For fixed $x$, the right-hand side of \eqref{eq:risk-V-decomp} differs from the reduced objective by the $\theta$-independent constant $f_{rs}(x)$, so every maximizer $\theta^* \ge x$ of $v_{rs}$ also maximizes $\theta \mapsto V_{rs,\pi_\theta}(x)$ over all $\theta \in \mathbb{Z}_{+}$.
    \item[(ii)] (Reduced objective) The function $f_{rs}$ does not depend on $\theta$, and
    \begin{equation}\label{eq:risk-V0-closed-form}
    v_{rs}(\theta) = \beta \theta - \alpha - f_{rs}(\theta).
    \end{equation}
\end{enumerate}
\end{proposition}

\begin{remark}[Certainty-equivalent interpretation]\label{rem:rs-ce}
Write $C_{x \to y}$ for the cost accumulated by the always-continue policy from state $x$ until the first passage to $y > x$. By \eqref{eq:fr-ce}, $f_{rs}(x)$ is the entropic certainty equivalent of the first-passage cost $C_{0 \to x}$, and the marginal cost $\Delta f_{rs}(x)$ is that of a single one-level excursion $C_{x \to x+1}$ (Proposition~\ref{prop:rs-ce}). For $\eta > 0$, risk seeking deflates marginal costs relative to the undiscounted risk-neutral benchmark,
$\Delta f_{rs} \le \Delta f$, and consequently pushes the optimal threshold upward
(Proposition~\ref{prop:rs-ce}(ii)--(iii)).
\end{remark}

\subsubsection{Optimal Threshold Under the Risk-Seeking Criterion}\label{sec:rs-seeking}

The increment bound of Lemma~\ref{lem:increments} extends to the entropic criterion, reducing the search for an optimal policy to the threshold family; the reduced objective \eqref{eq:risk-V0-closed-form} then yields the optimal threshold exactly. The increment $\Delta f_{rs}(x)$ is the marginal risk-sensitive cost of increasing the state from $x$ to $x+1$.

\begin{lemma}[Increment bound]\label{lem:rs-increments}
For all $x \ge 0$, $V_{rs}(x+1) - V_{rs}(x) \le \beta$.
\end{lemma}

\begin{proposition}[Threshold policies suffice]\label{prop:rs-threshold-wlog}
The stopping region $\{x \ge 0 : V_{rs}(x) = \beta x - \alpha\}$ is nonempty and up-closed; writing $\theta^*_{rs}$ for its smallest element, the threshold policy $\pi_{\theta^*_{rs}}$ is optimal among all admissible policies. In particular, $V_{rs}(0) = \max_{\theta \ge 0} v_{rs}(\theta)$.
\end{proposition}

\begin{theorem}[Optimal risk-seeking policy]\label{thm:RS}
Let $\eta > 0$.
\begin{enumerate}
    \item[(i)] The reduced objective is strictly discretely concave: $\Delta^2 v_{rs}(\theta) \le -c$ for all $\theta \ge 0$.
    \item[(ii)] The smallest optimal threshold is
    \begin{equation}\label{eq:optimal-theta-rs}
    \theta^*_{rs} = \min \bigl\{ \theta \in \mathbb{Z}_{+} : \Delta f_{rs}(\theta) \ge \beta \bigr\},
    \end{equation}
    and it satisfies $1 \le \theta^*_{rs} \le \lceil \beta / c \rceil$.

    \item[(iii)] The optimal threshold is unique unless $\Delta f_{rs}(\theta^*_{rs}) = \beta$, in which case exactly the two adjacent thresholds $\theta^*_{rs}$ and $\theta^*_{rs} + 1$ are optimal. Moreover, $V_{rs}(0) = v_{rs}(\theta^*_{rs}) = \max_{\theta \ge 0} v_{rs}(\theta)$, and $\pi_{\theta^*_{rs}}$ is optimal among all admissible policies.
\end{enumerate}
Equivalently, $\theta^*_{rs}$ is the smallest state at which the marginal risk-sensitive cost of increasing the state reaches the marginal processing reward.
\end{theorem}

\subsection{Risk-Averse Optimal Threshold: Robust Dual Formulation}\label{sec:rs-averse}

For $\eta<0$, the unregularized entropic value remains finite, but the standard bounded-cost/contractive fixed-point apparatus does not apply on the unbounded space, and the continuation branch has no uniform lower control; indeed, Lemma~\ref{lem:ra-continuation-diverges} shows that it can diverge downward under strong risk aversion. We therefore introduce a capped KL-robust dual as a tractable regularization: it upper-bounds the entropic problem and restores tractability, admits a unique fixed point in the natural linear-growth class, preserves threshold optimality, and yields an explicit truncated-Gibbs worst-case law. Throughout this section, Assumptions~\ref{ass:undiscounted} and~\ref{ass:profitability} are in force; in particular, $\delta=1$. 

\begin{lemma}[Downward divergence of the unregularized risk-averse continuation branch]
\label{lem:ra-continuation-diverges}
Let $\eta<0$. Define
\[
\Lambda_\eta
:=e^{\eta\beta}
\Bigl(1-p+p e^{-\eta(\beta+\gamma)}\Bigr).
\]
For each $K\in\mathbb R$, define the reward class
\[
\mathscr L_K
:=\left\{
u:\mathbb Z_+\to\mathbb R:
\beta x-\alpha\le u(x)\le\beta x+K
\quad\text{for every }x\in\mathbb Z_+
\right\},
\]
and write the unregularized risk-averse continuation branch as
\[
\widetilde Q_{ra}(u)(x)
:=-cx+\frac1\eta\log
\mathbb E\!\left[
e^{\eta\{u(\hat X_x)-\gamma A_x\}}
\right],
\qquad
\hat X_x=x-A_x+D .
\]
If $\Lambda_\eta>1$, then, uniformly over $u\in\mathscr L_K$,
\[
\widetilde Q_{ra}(u)(x)
\le
-\left(c-\frac{\log\Lambda_\eta}{\eta}\right)x
+K+\frac{\log(1-q+qe^{\eta\beta})}{\eta}
\longrightarrow-\infty
\qquad(x\to\infty).
\]
Moreover, $\Lambda_\eta>1$ for all sufficiently large $|\eta|$.
\end{lemma}

\noindent{\it Dual representation via the Donsker--Varadhan variational principle \citep{dupuis2011weak}.} For any bounded random variable $Z$, probability measure $\mathbb{P}$, and $\eta < 0$,
\begin{equation}\label{eq:dv}
\frac{1}{\eta} \log \mathbb{E}_\mathbb{P}\!\left[e^{\eta Z}\right] = \inf_{\mathbb{Q} \ll \mathbb{P}} \left\{ \mathbb{E}_\mathbb{Q}[Z] - \frac{1}{\eta} D_{\mathrm{KL}}(\mathbb{Q} \,\|\, \mathbb{P}) \right\},
\end{equation}
where $D_{\mathrm{KL}}(\mathbb{Q} \,\|\, \mathbb{P}) = \mathbb{E}_\mathbb{Q}[\log(d\mathbb{Q}/d\mathbb{P})]$ is the Kullback--Leibler divergence (relative entropy) and the infimum is attained by the Gibbs measure $d\mathbb{Q}^*/d\mathbb{P} \propto e^{\eta Z}$.

Applying \eqref{eq:dv} to the continuation branch gives the uncapped dual continuation branch
\begin{equation}
\begin{aligned}
\widetilde{Q}_{ra}(x)
=-cx+\inf_{\mathbb{Q}_x\in\mathcal{P}(\cdot\mid x)}
\Bigl\{
\mathbb{E}_{\mathbb{Q}_x}[V_r(\hat{X}_x)-\gamma A_x]-\frac{1}{\eta}D_{\mathrm{KL}}(\mathbb{Q}_x\,\|\,\mathbb{P}_x)
\Bigr\}.
\end{aligned}
\end{equation}
Hence
\begin{equation}\label{eq:bellman-dual}
 V_r(x)=\max\!\left\{S(x),\,\widetilde{Q}_{ra}(x)\right\},
\end{equation}
where $\mathbb{P}_x$ is the nominal one-step law of the pair $(A_x, D)$, or equivalently of $\hat{X}_x$---throughout the risk-averse analysis the symbol $\mathbb{P}_x$ denotes this one-step kernel, a slight overload of the path-law notation of Section~\ref{sec:formulation}---and $\mathcal{P}(\cdot\mid x)$ is the set of probability measures absolutely continuous with respect to $\mathbb{P}_x$. Equation~\eqref{eq:bellman-dual} is the Bellman equation of a stochastic game: the decision maker chooses a stopping time, and an adversary (Nature) chooses at each state a perturbed transition kernel $\mathbb{Q}_x$ that distorts the abandonment and arrival distributions so as to minimize the decision maker's payoff, with each deviation from the nominal model penalized by the positive term $-(1/\eta) D_{\mathrm{KL}}(\mathbb{Q}_x \,\|\, \mathbb{P}_x)$.

\begin{remark}\label{rem:regimes}
The dual representation exposes the structural difference between the risk regimes. For $\eta > 0$, the infimum in \eqref{eq:dv} becomes a supremum, and the adversary cooperates with the decision maker, seeking scenarios with high reward. For $\eta < 0$, the adversary is antagonistic, seeking scenarios with high cost, but is penalized for deviating excessively from the nominal model.
\end{remark}

\noindent{\it Regularization via uncertainty sets.} 
The uncapped dual formulation \eqref{eq:bellman-dual} is finite at each state, but it does not provide the uniform control needed on the unbounded state space: with the adversary constrained only by the KL penalty, distorted one-step payoffs need not satisfy a uniform integrability or drift estimate across all $x$, and the associated fixed-point problem is outside the standard contraction framework. We therefore cap the likelihood ratio of the adversary's distortion. 
For $\kappa > 1$, define
\begin{equation}\label{eq:uncertainty-set}
\mathcal{Q}_\kappa(x) = \left\{ \mathbb{Q}_x \ll \mathbb{P}_x : 0 \le \frac{d\mathbb{Q}_x}{d\mathbb{P}_x} \le \kappa \right\},
\end{equation}
where $\mathbb{P}_x = \mathbb{P}_x^{(A)} \otimes \mathbb{P}_x^{(D)}$ is the nominal joint distribution of the abandonment count $A_x \sim \operatorname{Binomial}(x,p)$ and the arrival $D \sim \operatorname{Bernoulli}(q)$. The cap serves two purposes: it keeps every distorted one-step cost uniformly integrable, and it yields the uniform Foster--Lyapunov drift estimate behind the well-posedness argument (Lemma~\ref{lem:ra-drift}). In addition, the KL penalty forces the worst-case measure to have full support (Lemma~\ref{lem:ra-fullsupport})---the proof exploits the boundary singularity of $m \log m$ at $m = 0$, not the cap---and this full-support property drives the uniqueness and verification arguments (Lemma~\ref{lem:ra-uniqueness} and Proposition~\ref{prop:ra-verification}).

\begin{remark}\label{rem:ra-Qx}
The nominal kernel is the product law $\mathbb{P}_x = \mathbb{P}_x^{(A)} \otimes \mathbb{P}_x^{(D)}$ on the finite support $\mathcal{X}_x = \{0, 1, \dots, x\} \times \{0, 1\}$. Because $0 < p, q < 1$, every atom carries positive nominal probability, so the absolute-continuity requirement $\mathbb{Q}_x \ll \mathbb{P}_x$ imposes no restriction: $\mathbb{Q}_x$ may be any probability vector on the $2(x+1)$ atoms. Two consequences are worth noting. First, $\mathbb{Q}_x$ need not be a product measure: the adversary may couple abandonment and arrivals, for instance by concentrating mass on atoms where heavy abandonment coincides with an arrival, inflating the abandonment penalty and the subsequent holding load simultaneously. Second, parametric perturbations $\operatorname{Binomial}(x, p') \otimes \operatorname{Bernoulli}(q')$ with $(p', q') \neq (p, q)$ form only a small subfamily of $\mathcal{P}(\cdot \mid x)$: the worst-case kernel \eqref{eq:ra-gibbs} reweights each atom $(a, d)$ individually according to the exponentiated continuation value $V_{ra}^{\kappa}(x - a + d) - \gamma a$, and therefore generally lies outside any such parametric family.
\end{remark}

Define the regularized continuation branch $Q_{ra}$ by capping the likelihood ratio in the uncapped branch $\widetilde{Q}_{ra}$:
\begin{equation}\label{eq:ra-continuation}
\begin{aligned}
Q_{ra}(x)
:=-cx+\inf_{\mathbb{Q}_x\in\mathcal{Q}_\kappa(x)}
\Bigl\{
\mathbb{E}_{\mathbb{Q}_x}\!\bigl[
V_{ra}^{\kappa}(\hat{X}_x)-\gamma A_x
\bigr]-\frac{1}{\eta}D_{\mathrm{KL}}(\mathbb{Q}_x\,\|\,\mathbb{P}_x)
\Bigr\}.
\end{aligned}
\end{equation}
The regularized Bellman equation is then
\begin{equation}\label{eq:bellman-robust}
V_{ra}^{\kappa}(x)
=\max\!\left\{S(x),\,Q_{ra}(x)\right\}.
\end{equation}
We refer to the fixed point of \eqref{eq:bellman-robust} as the \textit{regularized value function} $V_{ra}^{\kappa}$.

\begin{theorem}[Optimal risk-averse policy under regularization]\label{thm:RA}
Let $\eta < 0$ and $\kappa > 1$.
\begin{enumerate}
    \item[(i)] (Well-posedness) The regularized Bellman equation~\eqref{eq:bellman-robust} admits a unique fixed point $V_{ra}^{\kappa}$ in the class of functions $u : \mathcal{S} \to \mathbb{R}$ satisfying $u \ge S$ and $u(x) \le \beta x + K$ for some constant $K$; the value iterates $V_0 = S$, $V_{n+1} = \mathcal{T}_{ra} V_n$ increase pointwise to $V_{ra}^{\kappa}$.
    \item[(ii)] (Threshold optimality) With $\theta^*_{ra} = \min\{ x : V_{ra}^{\kappa}(x) = S(x) \}$, the threshold policy $\pi_{\theta^*_{ra}}$ is optimal among all admissible stopping rules, and $V_{ra}^{\kappa}$ is the value of the robust stopping game. The optimal threshold is unique unless $Q_{ra}(\theta^*_{ra}) = S(\theta^*_{ra})$, in which case exactly the two adjacent policies $\pi_{\theta^*_{ra}}$ and $\pi_{\theta^*_{ra}+1}$ are optimal.
    \item[(iii)] (Worst-case characterization) At every state $x$, the infimum in \eqref{eq:bellman-robust} is attained by a unique $\mathbb{Q}^*_x \in \mathcal{Q}_\kappa(x)$, with truncated exponential-tilting likelihood ratio
    \begin{equation}\label{eq:ra-gibbs}
    \frac{d\mathbb{Q}^*_x}{d\mathbb{P}_x}(a,d)
    =\min\!\left\{
    \kappa,\,
    \lambda_x e^{\eta[V_{ra}^{\kappa}(\hat{X}_x)-\gamma a]}
    \right\},
    \end{equation}
    where $\hat{X}_x=x-a+d$ on the atom $(a,d)$, and $\lambda_x>0$ is chosen so that $\mathbb{Q}^*_x$ is a probability measure. In particular, $\mathbb{Q}^*_x$ has full support.
\end{enumerate}
\end{theorem}
The proof is given in Section~\ref{sec:RA}.

\noindent{\it The role and choice of $\kappa$.} The likelihood-ratio cap $\kappa$ calibrates the amount of model uncertainty against which the decision maker hedges. At $\kappa = 1$ the uncertainty set \eqref{eq:uncertainty-set} collapses to the nominal kernel and the regularized criterion reduces to the risk-neutral one; as $\kappa \to \infty$ the likelihood-ratio constraint becomes inactive state-by-state and the capped dual formally approaches the uncapped entropic dual; we use this only as a formal consistency limit, since the uncapped representation is justified here only under the one-step boundedness and exponential-integrability conditions noted above.
The next proposition records that the regularized criterion interpolates monotonically---in value and in policy---between these two extremes.

\begin{proposition}[Monotonicity and consistency in the robustness bound]\label{prop:kappa-limit}
Fix $\eta < 0$ and suppose Assumption~\ref{ass:profitability} holds.
\begin{enumerate}
    \item[(i)] (Monotonicity in $\kappa$) The regularized value function $V_{ra}^{\kappa}(x)$ is nonincreasing in $\kappa$ for every $x \in \mathbb{Z}_{+}$, and the smallest optimal threshold $\theta^*_{ra}(\kappa)$ is nonincreasing in $\kappa$: a larger uncertainty set strengthens the adversary, depresses the continuation value, and shifts stopping earlier.
    \item[(ii)] (Nominal limit) As $\kappa \downarrow 1$, $V_{ra}^{\kappa}(x) \uparrow V_{rn}(x)$ for every $x \in \mathbb{Z}_{+}$, and the smallest optimal thresholds coincide for all $\kappa$ sufficiently close to $1$: there exists $\bar\kappa > 1$ such that $\theta^*_{ra}(\kappa) = \theta^*_{rn}$ for all $\kappa \in (1, \bar\kappa)$.
\end{enumerate}
\end{proposition}

In applications, $\kappa$ can be calibrated to the statistical plausibility of alternative models---for instance, as the maximal likelihood ratio between the nominal kernel and the kernels estimated at the boundary of a confidence region for the nominal parameters $(p, q)$---so that the adversary distorts the model only within statistically credible limits; $\theta^*_{ra}(\kappa)$ then traces out the risk-averse discount as a function of the uncertainty budget (Section~\ref{sec:num-robust-pmf}, Figure~\ref{fig:robust_pmf}).

\begin{remark}[Numerical solution]\label{rem:ra-numerics}
The regularized robust Bellman equation~\eqref{eq:bellman-robust} is solved numerically by monotone value iteration on a truncated state space, and the ordering of Proposition~\ref{prop:regime-comparison}(i) provides the truncation point used in the implementation $X_{\max} := \theta^*_{rn}$. Implementation details are deferred to Section~\ref{app:algorithms}.
\end{remark}

\subsection{Threshold Ordering and Convergence Across Risk Regimes}\label{sec:regime-comparison}

\begin{proposition}[Comparison across risk regimes]\label{prop:regime-comparison}
Let $\delta = 1$, and suppose Assumption~\ref{ass:profitability} holds. For $\eta > 0$, write $V_{rs,\eta}$ and $\theta^*_{rs,\eta}$ for the value function and the smallest optimal threshold of the risk-seeking problem; for $\eta < 0$, write $V_{ra,\eta}^{\kappa}$ for the risk-averse regularized value function and $\theta^*_{ra,\eta}$ for the optimal threshold of the risk-averse problem, with the robustness parameter $\kappa > 1$ fixed.
\begin{enumerate}
    \item[(i)] (Threshold ordering) The risk-neutral threshold lies between the risk-averse and risk-seeking thresholds:
    \[
    \theta^*_{ra,\eta'}
    \le \theta^*_{rn}
    \le \theta^*_{rs,\eta},
    \qquad \eta' < 0 < \eta.
    \]

    \item[(ii)] (Risk-seeking limit) For every $x \in \mathbb{Z}_+$,
    \[
    V_{rs,\eta}(x)\downarrow V_{rn}(x)
    \qquad\text{as }\eta\downarrow0.
    \]
    Moreover, the risk-seeking threshold stabilizes: $\theta^*_{rs,\eta}$ is constant for all sufficiently small $\eta > 0$, and its limit is an optimal risk-neutral threshold---namely $\theta^*_{rn}$ when the undiscounted risk-neutral optimal threshold is unique, and $\theta^*_{rn} + 1$, the second optimal risk-neutral threshold, in the tie case $\Delta f(\theta^*_{rn}) = \beta$.

    \item[(iii)] (Risk-averse limit) For every $x \in \mathbb{Z}_+$ and fixed $\kappa > 1$,
    \[
    V_{ra,\eta}^{\kappa}(x)\uparrow V_{rn}(x)
    \qquad\text{as }\eta\uparrow0.
    \]
    Moreover, there exists $\underline{\eta} < 0$ such that
    \[
    \theta^*_{ra,\eta} = \theta^*_{rn}
    \qquad\text{for all }\eta \in (\underline{\eta}, 0):
    \]
    the risk-averse limit always selects the smallest optimal risk-neutral threshold, with no uniqueness assumption required.
\end{enumerate}
\end{proposition}

\section{Extension: Repeated Clearing}\label{sec:renewal}

In the applications motivating this study---manufacturing transfer batches, order consolidation, packet aggregation, and platelet pooling---clearing is often a recurring event rather than a terminal one: the buffer is emptied, accumulation resumes, and the controller faces the same decision again. Some existing clearing and batch-service literature accordingly adopts long-run average criteria, in contrast to the optimal stopping criteria of Sections~\ref{sec:rn}--\ref{sec:rs-averse}.
This section shows that the first-passage construction of Section~\ref{sec:decoupling} also yields reduced objectives for repeated clearing. The repeated-clearing formulation prices not only the cost accrued within a cycle but also the length of the cycle itself.

\noindent{\it Timing convention and regenerative cycles.}
At the beginning of a period in state $x$, the controller may continue or clear. Continuation at discount factor $\delta$ incurs the realized holding and abandonment costs $c x+\delta\gamma A_x$ and produces the transition $\hat X_x=x-A_x+D$. Clearing is permitted only when $x\ge1$: it generates the terminal reward $\beta x-\alpha$, resets the buffer immediately to $0$, and no holding or abandonment cost is charged at the clearing epoch. The next cycle then begins from state $0$. Under a stationary clearing policy, the strong Markov property at each clearing epoch makes the cycles independent; under a threshold policy $\pi_\theta$, they are i.i.d. with common length $T_\theta$.

For the threshold policy $\pi_\theta$, define the discounted cycle rewards by
\[
R_\theta(\delta) := \delta^{T_\theta}(\beta\theta-\alpha) -\sum_{t=0}^{T_\theta-1}\delta^t \bigl(cX_t+\delta\gamma A_t(X_t)\bigr).
\]
We consider the long-run average reward generated by repeated cycles and the total discounted reward over infinitely many cycles. Let $(T_{\theta,k},R_{\theta,k}(\delta))_{k\ge1}$ be i.i.d.\ copies of
$(T_\theta,R_\theta(\delta))$, and set $S_0:=0$ and
$S_k:=\sum_{j=1}^{k}T_{\theta,j}$. Thus $S_{k-1}$ is the clearing epoch at
which cycle $k$ begins.

For $\delta\in(0,1)$, the total discounted reward over infinitely many cycles
is
\begin{equation}\label{eq:discounted-cycle-value}
V_\theta(\delta)
:=\mathbb{E}\!\left[
\sum_{k=1}^{\infty}\delta^{S_{k-1}}R_{\theta,k}(\delta)
\right]
=\frac{\mathbb{E}[R_\theta(\delta)]}
{1-\mathbb{E}[\delta^{T_\theta}]},
\end{equation}
where independence between the past cycle lengths and the reward in cycle
$k$ gives
\[
\mathbb{E}\!\left[\delta^{S_{k-1}}R_{\theta,k}(\delta)\right]
=\bigl(\mathbb{E}[\delta^{T_\theta}]\bigr)^{k-1}
 \mathbb{E}[R_\theta(\delta)].
\]
The geometric series converges because $T_\theta\ge1$ and $\delta<1$.

For the long-run average criterion, take $\delta=1$ within each cycle and
define
\begin{equation}\label{eq:average-cycle-value}
\overline R_\theta
:=\lim_{K\to\infty}
\frac{\sum_{k=1}^{K}R_{\theta,k}(1)}
{\sum_{k=1}^{K}T_{\theta,k}}
=\frac{\mathbb{E}[R_\theta(1)]}{\mathbb{E}[T_\theta]}
\qquad\text{a.s.},
\end{equation}
provided $\mathbb{E}[T_\theta]<\infty$ and
$\mathbb{E}|R_\theta(1)|<\infty$. The second equality is the
renewal--reward theorem. Hence both criteria reduce to one-cycle
calculations: $\overline R_\theta>0$ precisely when
$\mathbb{E}[R_\theta(1)]>0$, while $V_\theta(\delta)>0$ precisely when
$\mathbb{E}[R_\theta(\delta)]>0$.
To distinguish the one-cycle primitives by discount factor, write $f(\cdot;\delta)$ and $g(\cdot;\delta)$ for the functions in Section~\ref{sec:decoupling}, and define
\[
N_\delta(\theta) :=
\beta\theta-\alpha-f(\theta;\delta).
\]
We say that \emph{Assumption~\ref{ass:profitability} holds at discount factor $\delta$} if $N_\delta(\theta)>0$ for some $\theta\ge1$; since $N_\delta(0)=-\alpha<0$, quantifying over $\theta\ge1$ is equivalent to quantifying over $\theta\in\mathbb Z_+$ used in Section~\ref{sec:rn-threshold}.

\begin{lemma}[Time-accumulation function]\label{lem:t-recursion}
Define $t:\mathbb Z_+\to\mathbb R_+$ by $t(0)=0$ and the forward recurrence
\begin{equation}\label{eq:t-recurrence}
t(x+1) = \frac{\omega_x\,t(x)+1-F_t(x)}{(1-p)^xq}, \qquad
F_t(x) = \mathbb E_{A_x,D}\!\left[ \mathbf 1_{\{A_x\ge1\}}\,t(\hat X_x) \right],
\end{equation}
where $\omega_x=1-(1-p)^x(1-q)$. Equivalently, the sequence defined by
\eqref{eq:t-recurrence} satisfies
\begin{equation}\label{eq:t-fixed}
t(x)=-1+\mathbb E\!\left[t(\hat X_x)\right], \qquad x\ge0.
\end{equation}
For every $\theta\ge1$,
\begin{equation}\label{eq:t-hitting-time}
t(\theta)=\mathbb E_0[T_\theta].
\end{equation}
\end{lemma}

\begin{proposition}[Repeated-clearing reduced objectives]\label{prop:renewal-objectives}
Fix $\theta\ge1$.
\begin{enumerate}
\item[(i)] (Average reward) The long-run average reward of $\pi_\theta$ is
\begin{equation}\label{eq:avg-objective}
\bar v(\theta) = \frac{N_1(\theta)}{t(\theta)}.
\end{equation}
\item[(ii)] (Discounted repeated clearing) For $\delta\in(0,1)$, the total discounted reward of $\pi_\theta$ over infinitely many cycles, started from $0$, is
\begin{equation}\label{eq:disc-renewal-objective}
v^{\mathrm{ren}}_\delta(\theta) = \frac{v_{rn}(\theta;\delta)}{1-\mathbb E_0[\delta^{T_\theta}]} = \frac{N_\delta(\theta)}{g(\theta;\delta)-1}.
\end{equation}
\item[(iii)] (Abelian consistency) For every fixed $\theta\ge1$,
\begin{equation}\label{eq:abelian-consistency}
(1-\delta)\,v^{\mathrm{ren}}_\delta(\theta) \longrightarrow \bar v(\theta)
\qquad \text{as } \delta \uparrow 1.
\end{equation}
\end{enumerate}
Parts~(i)--(iii) are policy-evaluation identities and require no profitability assumption; Assumption~\ref{ass:profitability} enters only through the optimization results below (Theorem~\ref{thm:renewal-unimodal} and Proposition~\ref{prop:renewal-threshold-wlog}).
\end{proposition}

The reduced objectives above evaluate threshold policies. The next proposition records the dynamic-programming step showing that restricting attention to threshold policies entails no loss of optimality for either repeated-clearing criterion.

\begin{proposition}[Threshold policies suffice]\label{prop:renewal-threshold-wlog}
Under Assumption~\ref{ass:profitability} at $\delta=1$, there exists a threshold policy $\pi_\theta$ with $\theta\ge1$ that is optimal among all stationary clearing policies for the average-reward regenerative problem. The same conclusion holds for the discounted problem at every $\delta\in(0,1)$ at which Assumption~\ref{ass:profitability} holds.
\end{proposition}

\begin{theorem}[Strict unimodality of the repeated-clearing objectives]\label{thm:renewal-unimodal}
Suppose that Assumption~\ref{ass:profitability} holds at the relevant discount
factor.
\begin{enumerate}
\item[(i)] (Average reward) The long-run average objective $\bar v$ is strictly unimodal over all thresholds $\theta\ge1$: there exists $\theta^u_1\ge1$ such that $\Delta\bar  v(\theta)>0$ for $1\le\theta<\theta^u_1$, $\Delta\bar v(\theta^u_1)\le0$, and $\Delta\bar v(\theta)<0$ for $\theta>\theta^u_1$.
\item[(ii)] (Discounted repeated clearing) For every $\delta\in(0,1)$, the objective $v^{\mathrm{ren}}_\delta$ is strictly unimodal over all thresholds $\theta\ge1$ in the same sense, with crossing index $\theta^u_\delta$.
\item[(iii)] (Optimal thresholds) The indices
\begin{equation}\label{eq:renewal-argmax}
\theta^*_{\mathrm{avg}}=\min\bigl\{\theta\ge1:\Delta\bar v(\theta)\le0\bigr\}, \qquad \theta^*_{\mathrm{ren}}(\delta)=\min\bigl\{\theta\ge1:\Delta v^{\mathrm{ren}}_\delta(\theta)\le0\bigr\}
\end{equation}
are well defined and finite, coincide with $\theta^u_1$ and $\theta^u_\delta$, and are the smallest optimal thresholds for repeated clearing; at most two thresholds are optimal for each criterion, and they differ by one. The corresponding threshold policies are optimal among all stationary clearing policies by Proposition~\ref{prop:renewal-threshold-wlog}.
\end{enumerate}
\end{theorem}

The proof is given in Section~\ref{app:renewal}. Theorem~\ref{thm:renewal-unimodal} justifies, for both repeated-clearing criteria, the same discrete-ascent procedure that Theorem~\ref{thm:RN} justifies for the optimal stopping objective: the first local maximizer encountered by a forward pass is the global one. 

\begin{remark}[Marginal-versus-average interpretation]\label{rem:marginal-average}
The crossing condition compares the average rate with the marginal rate of a one-level cycle extension: $\Delta\bar v(\theta)\le0$ exactly when $\bar v(\theta)\ge\xi_1(\theta)$, and 
$\Delta v^{\mathrm{ren}}_\delta(\theta) \le 0$ exactly when $v^{\mathrm{ren}}_\delta(\theta)\ge\xi_\delta(\theta)$ (see \eqref{eq:renewal-sign-xi} in the proof). The marginal rate $\xi_1(\theta)$ prices the net reward of extending the cycle from $\theta$ to $\theta+1$ against the marginal expected cycle length, and its limit $-\mathcal C_1q/p$ is the negative of the stationary holding-and-abandonment cost rate of the always-continue chain, whose stationary mean inventory level is $q/p$. The optimal batch size is the first level at which the marginal rate of a larger batch falls below the average rate generated by the current batch---a marginal-versus-average rule in the spirit of classical clearing theory.
\end{remark}

The repeated-clearing criteria also sharpen the interpretation of the optimal stopping threshold. Let $\theta^*_{rn}(1)$ and $\theta^*_{rn}(\delta)$ denote the smallest optimal stopping thresholds at $\delta=1$ and at $\delta\in(0,1)$, respectively.

\begin{proposition}\label{prop:renewal-ordering}
Suppose that Assumption~\ref{ass:profitability} holds. The smallest optimal thresholds for repeated clearing satisfy
\begin{equation}\label{eq:renewal-ordering}
\theta^*_{\mathrm{avg}}\le\theta^*_{rn}(1), \qquad \theta^*_{\mathrm{ren}}(\delta)\le\theta^*_{rn}(\delta).
\end{equation}
\end{proposition}

\begin{remark}
The optimal stopping objective prices erosion within a single cycle. The repeated-clearing objectives additionally price the opportunity cost of the cycle length through $t(\theta)$ in \eqref{eq:avg-objective} and through $g(\theta;\delta)-1$ in \eqref{eq:disc-renewal-objective}. Proposition~\ref{prop:renewal-ordering} formalizes the resulting over-batching bias: the optimal stopping threshold is an upper bound on the corresponding repeated-clearing threshold. 
\end{remark}

\begin{remark}[Entropic criteria]
The additive cycle-reward reduction is specific to additive criteria. Under an entropic criterion, cycle rewards combine multiplicatively, so i.i.d. cycles do not yield the scalar formulas \eqref{eq:avg-objective} and \eqref{eq:disc-renewal-objective}. A risk-sensitive theory for repeated clearing would require an ergodic multiplicative-functional or recursive-entropic formulation \citep{borkar2002risk}; we do not pursue that extension in this paper.
\end{remark}

\section{Numerical Study} \label{sec:numerical}

\subsection{Experimental Design and Baseline Calibration}\label{sec:num-setup}

\noindent{\it Baseline parameters.}
The baseline is chosen so that the profitability condition \eqref{eq:profitability} holds with a nontrivial optimal threshold, the buffer displays appreciable stochastic variability, and the three risk regimes admit interior, mutually comparable solutions. Unless otherwise stated, every experiment below fixes the baseline parameters of Table~\ref{tab:baseline}; the table also reports the grid ranges used in the comparison and sensitivity experiments.

\begin{table}[htbp]
    \centering
    \caption{Baseline calibration and grid ranges used in the sensitivity analysis.}
    \label{tab:baseline}
    \begin{tabular}{lcc}
        \toprule
        Parameter & Baseline & Grid range \\
        \midrule
        Stopping cost $\alpha$ & 100.0 & [10.0, 1000.0] \\
        Reward coefficient $\beta$ & 10.0 & [1.0, 100.0] \\
        Holding cost $c$ & 0.01 & [0.001, 0.1] \\
        Abandonment penalty $\gamma$ & 0.5 & [0.1, 10.0] \\
        Arrival probability $q$ & 0.8 & [0.5, 1] \\
        Abandonment probability $p$ & 0.01 & [0.001, 0.1] \\
        Discount factor $\delta$ & 0.99 & [0.9, 1] \\
        \bottomrule
    \end{tabular}
\end{table}

\noindent{\it Solution procedures.}
\begin{enumerate}
    \item[(a)] \textit{Risk-neutral forward recursion}: forward iteration of \eqref{eq:f-recurrence}--\eqref{eq:g-recurrence}, followed by the discrete ascent on $v_{rn}$ justified by Theorem~\ref{thm:RN}.
    \item[(b)] \textit{Risk-seeking forward recursion}: forward iteration of \eqref{eq:fr-recurrence}, terminated at the first index satisfying $\Delta f_{rs}(\theta) \ge \beta$, per Theorem~\ref{thm:RS}.
    \item[(c)] \textit{Regularized risk-averse value iteration}: monotone value iteration of the robust operator in \eqref{eq:bellman-robust}, initialized at the stopping reward and truncated at $X_{\max} = \theta^*_{rn}$ per Proposition~\ref{prop:regime-comparison}(i), with the inner Donsker--Varadhan dual solved by bisection on $\log \lambda_x$ (Section~\ref{app:algorithms}).
    \item[(d)] \textit{Repeated-clearing forward recursion}: forward iteration of \eqref{eq:f-recurrence}--\eqref{eq:g-recurrence} at the relevant discount factor and, for the average-reward criterion, of the time-accumulation recursion \eqref{eq:t-recurrence}, followed by a discrete ascent on the repeated-clearing objectives \eqref{eq:avg-objective} and \eqref{eq:disc-renewal-objective} terminated at the first nonpositive marginal improvement, per Theorem~\ref{thm:renewal-unimodal}.
\end{enumerate}
All computations are performed in double precision; convergence tolerances and truncation rules are reported with each experiment.

\subsection{Sensitivity to the Economic Primitives}\label{sec:num-sensitivity}
\begin{figure}[htbp]
    \centering
    \includegraphics[scale=0.5]{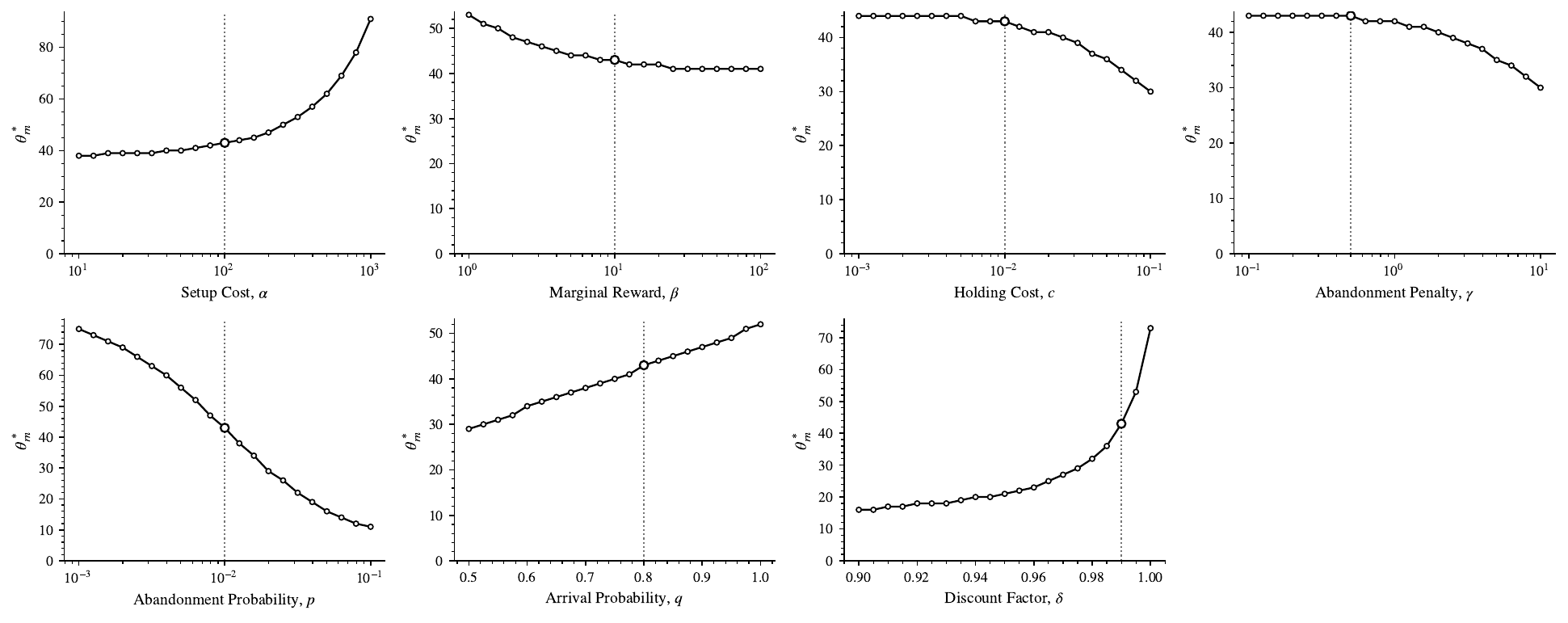}
    \caption{Sensitivity of $\theta^*_{rn}$ to the economic primitives. }
    \label{fig:sensitivity}
\end{figure}

Starting from the baseline, we vary one primitive at a time over the grid ranges reported in Table~\ref{tab:baseline}, holding all remaining parameters fixed.
Figure~\ref{fig:sensitivity} reports comparative statics for the optimal risk-neutral threshold $\theta^*_{rn}$ with respect to the seven primitives.
The optimal threshold increases with the stopping cost $\alpha$---from $\theta^*_{rn}=43$ at the baseline to $91$ at $\alpha=1000$---because larger fixed costs require larger batches for amortization (Lemma~\ref{lem:sensitivity-costs}(ii)). It decreases with the holding cost $c$ and the abandonment penalty $\gamma$, both of which raise the effective cost of accumulation (Lemma~\ref{lem:sensitivity-costs}(i)), and with the abandonment probability $p$, which additionally steepens the downward drift. It increases with the arrival probability $q$, which raises the upward drift of the inventory process.

The response to $\delta$ is strong: $\theta^*_{rn}$ rises from $16$ at $\delta=0.9$ to $73$ at $\delta=1$, as greater patience lowers the opportunity cost of extending accumulation. The sign of the response in the $\beta$ panel is not determined by cost monotonicity alone: at the baseline $\delta=0.99$, the threshold decreases from $53$ at $\beta=1$ to $41$ at $\beta=100$.
The sign of this response is governed by a reward--discount phase transition, discussed in Section~\ref{sec:num-phase}.

\subsection{Operational Implications}\label{sec:num-operational}
\begin{figure}[htbp]
    \centering
    \includegraphics[scale=0.65]{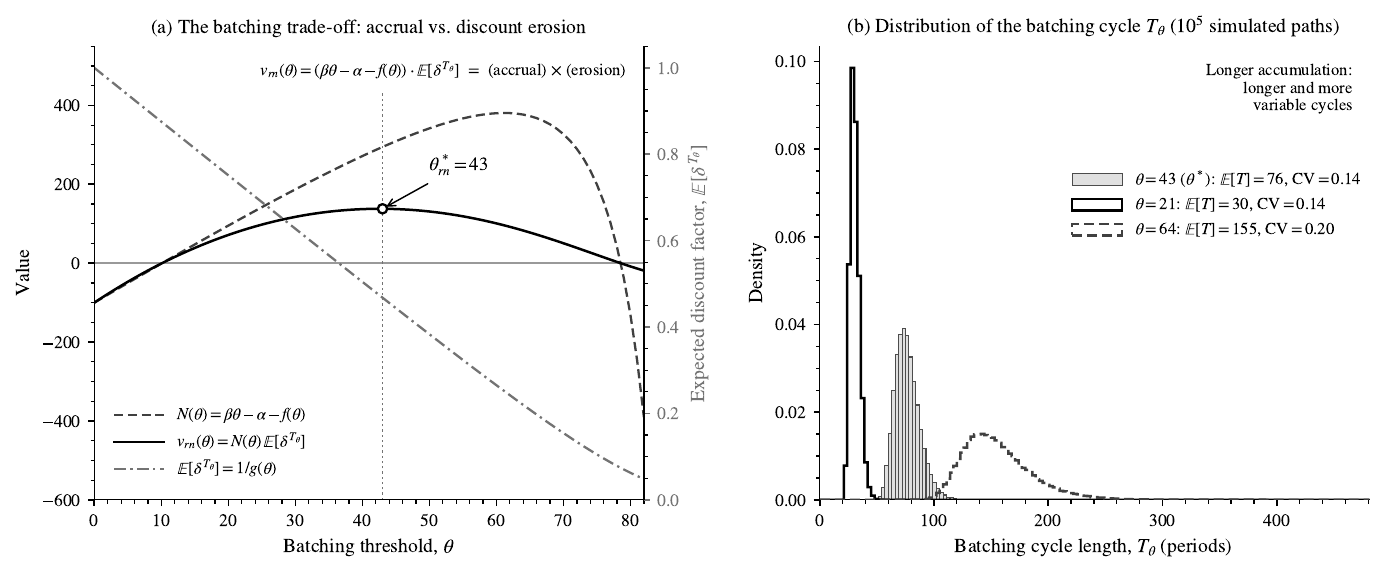}
    \caption{Operational implications at the discounted baseline. (a) Factorization of $v_{rn}(\theta)$ into accrued nominal value $N_{rn}(\theta)$ (solid) and expected discount factor $\mathbb{E}[\delta^{T_\theta}]=1/g(\theta)$ (dashed); the product is maximized at $\theta^*_{rn}=43$. (b) Simulated batching-cycle distribution for thresholds $0.5\theta^*$, $\theta^*$, and $1.5\theta^*$ based on $10^5$ independent paths.}
    \label{fig:operational}
\end{figure}

We illustrate the factorization of the discounted reduced objective from
Remark~\ref{rem:fg-meaning} and discuss its operational content. Recall that
\[
v_{rn}(\theta) = N_{rn}(\theta)\,\mathbb{E}\!\left[\delta^{T_\theta}\right] = \frac{\beta\theta-\alpha-f(\theta)}{g(\theta)},
\]
where
\(N_{rn}(\theta)=\beta\theta-\alpha-f(\theta)\) and
\(\mathbb{E}[\delta^{T_\theta}]=1/g(\theta)\).
In Figure~\ref{fig:operational}(a), the accrued nominal value $N_{rn}(\theta)$ rises almost linearly at first---the reward coefficient $\beta$ dominates the accumulating holding and abandonment costs---peaks at approximately $381$ near $\theta=61$, and then declines, crossing zero between $\theta=78$ and $\theta=79$: beyond that point, accrued costs outweigh the batch reward even before discounting. Over the same range, the expected discount factor $\mathbb{E}[\delta^{T_\theta}]$ falls from $1$ to below $0.05$. The optimum $\theta^*_{rn}=43$ therefore lies well inside the region where nominal accrual is still increasing: discount erosion binds before cost accrual does. At the optimum, the policy accrues nominal value $N_{rn}(43)\approx295$; waiting erodes more than half of it, with $\mathbb{E}[\delta^{T_{43}}]\approx0.47$, leaving $v_{rn}(43) \approx 138$.

Figure~\ref{fig:operational}(b) reports the simulated distribution of the batching cycle $T_\theta$ under $\pi_\theta$ for $\theta\in\{21,43,64\}$, using $10^5$ independent paths per threshold. The mean cycle length grows faster than the threshold: halving $\theta^*_{rn}$ from $43$ to $21$ shortens the mean cycle from $76$ to $30$ periods, whereas raising it to $64$ lengthens the mean to $155$. The mechanism is the abandonment hazard: as inventory accumulates, abandonment drains the net drift $q-p\,\mathbb{E}[X_t]$, so the chain climbs more slowly the higher it must go. Over-batching is therefore disproportionately expensive in capital lock-up, a convexity that the scalar objective alone does not display.

\subsection{The Risk-Seeking Premium and the Risk-Averse Discount}\label{sec:num-premium}
\begin{figure}[htbp]
    \centering
    \includegraphics[scale=0.6]{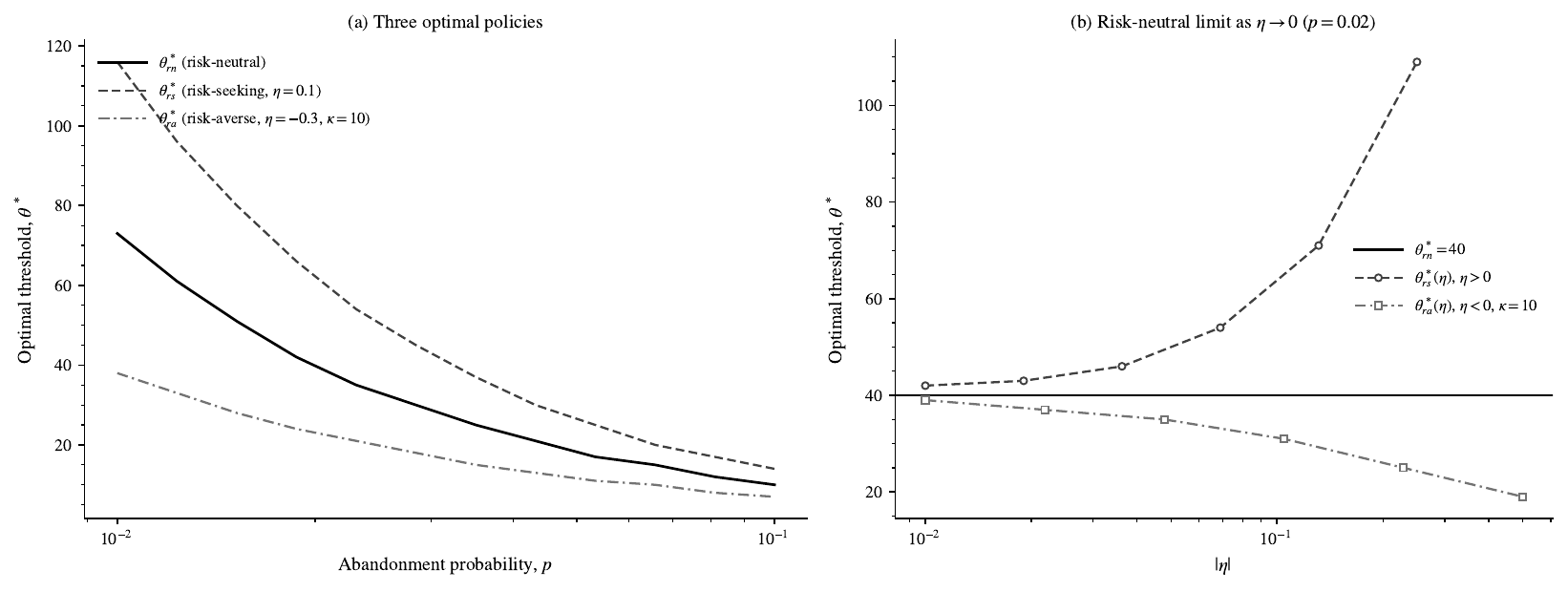}
    \caption{Threshold comparison under the undiscounted calibration. (a) Optimal thresholds under risk-neutral ($\theta^*_{rn}$), risk-seeking ($\theta^*_{rs}$, $\eta>0$), and KL-robust risk-averse ($\theta^*_{ra}$, $\eta<0$, $\kappa=10$) criteria versus abandonment probability $p$; shaded regions mark the risk-seeking premium and risk-averse discount. (b) Convergence of both risk-adjusted thresholds to $\theta^*_{rn}$ as $\eta\to 0$ at fixed $p=0.02$.}
    \label{fig:convergence}
\end{figure}

This subsection compares the three criteria under the undiscounted calibration $\delta=1$ of Assumption~\ref{ass:undiscounted}. We define the \textit{risk-seeking premium} as $\theta^*_{rs}-\theta^*_{rn}$ and the \textit{risk-averse discount} as $\theta^*_{rn}-\theta^*_{ra}$. Figure~\ref{fig:convergence} summarizes the comparison: panel~(a) traces the three thresholds as functions of $p$, and panel~(b) tracks the risk-adjusted thresholds as $\eta\to 0$; the ordering is that of Proposition~\ref{prop:regime-comparison}(i). Throughout, $\theta^*_{ra}$ is computed under the regularized criterion of Section~\ref{sec:rs-averse} at a fixed cap $\kappa$, reported in each figure.

\smallskip

\begin{figure}[htbp]
    \centering
    \includegraphics[scale=0.6]{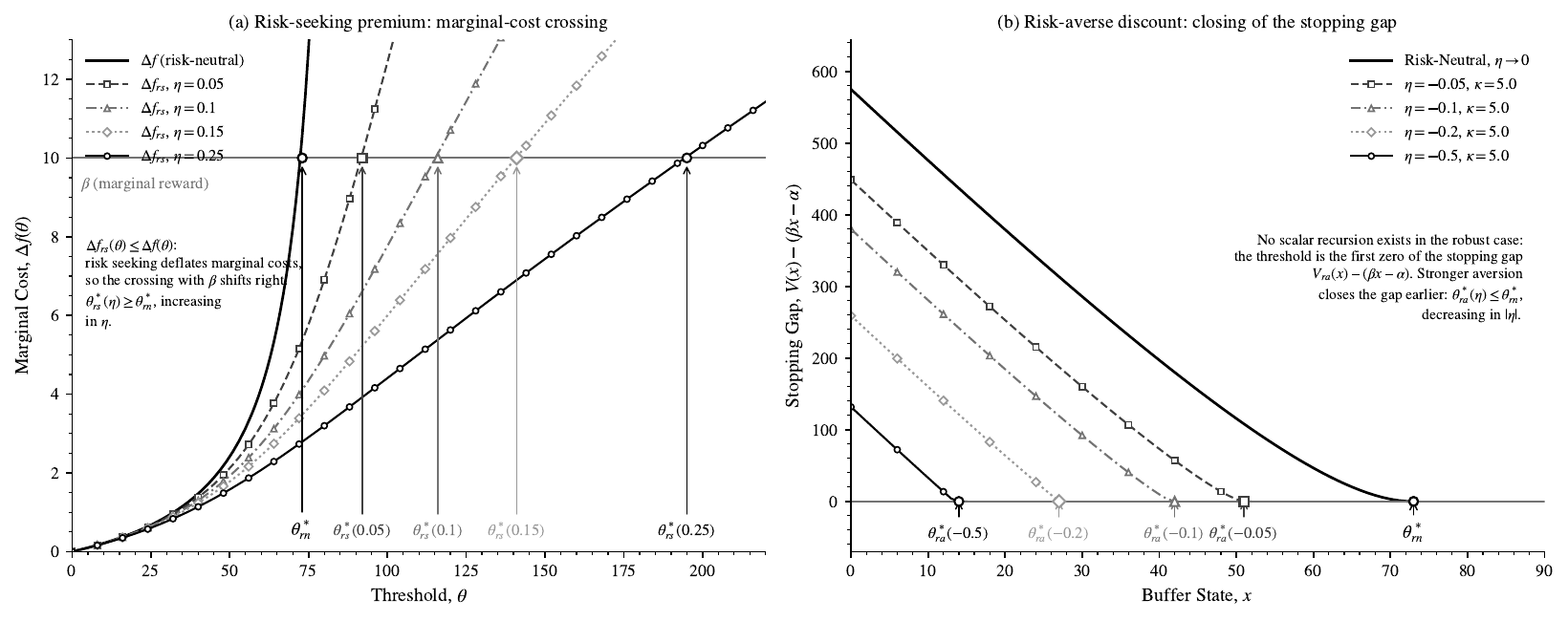}
    \caption{Mechanisms behind the risk-seeking premium and risk-averse discount ($\delta=1$, $\kappa=5$). (a) Risk seeking deflates marginal continuation cost, $\Delta f_{rs}(\theta)\le\Delta f(\theta)$, shifting the crossing with $\beta$ to the right. (b) Risk aversion shrinks the stopping gap $V^\kappa_{ra}(x)-S(x)$, shifting its first zero to the left.}
    \label{fig:marginal_cost_rs}
\end{figure}

\noindent{\it Risk-seeking premium ($\eta>0$).}
The entropic tilt deflates the marginal accumulation cost pointwise, $\Delta f_{rs}(\theta)\le\Delta f(\theta)$ (Proposition~\ref{prop:rs-ce}(ii)), shifting the first crossing of the marginal-cost curve with the reward coefficient $\beta$ to the right: at the undiscounted calibration, the crossing moves from $\theta^*_{rn}=73$ to $\theta^*_{rs}=92,116,141,195$ as $\eta$ ranges over
$\{0.05,0.10,0.15,0.25\}$ (Figure~\ref{fig:marginal_cost_rs}(a)).

\smallskip
\noindent{\it Risk-averse discount ($\eta<0$).}
Because the robust inner optimization destroys the scalar forward recursion available in the risk-seeking case, the threshold is characterized through the stopping gap $V^\kappa_{ra}(x)-S(x)$ rather than through a marginal-cost crossing; Figure~\ref{fig:marginal_cost_rs}(b) shows that the gap is positive on the continuation region and closes at the optimal threshold. As $|\eta|$ grows, the gap closes earlier: at $\kappa=5$, the first zero moves from $\theta^*_{rn}=73$ to $\theta^*_{ra}=51,42,27,14$ for $\eta\in\{-0.05,-0.10,-0.20,-0.50\}$.

\smallskip

\noindent{\it Continuity at $\eta=0$.}
Figure~\ref{fig:convergence}(b) confirms the threshold continuity of Proposition~\ref{prop:regime-comparison}(ii)--(iii): the risk-averse threshold approaches $\theta^*_{rn}$ from below as $\eta\uparrow0$, the risk-seeking threshold approaches it from above as $\eta\downarrow 0$, and the two coincide with $\theta^*_{rn}$ for sufficiently small $|\eta|$ on the reported grid.

\subsection{Worst-Case Distortion of the Abandonment Distribution}\label{sec:num-robust-pmf}
\begin{figure}[htbp]
    \centering
    \includegraphics[scale=0.6]{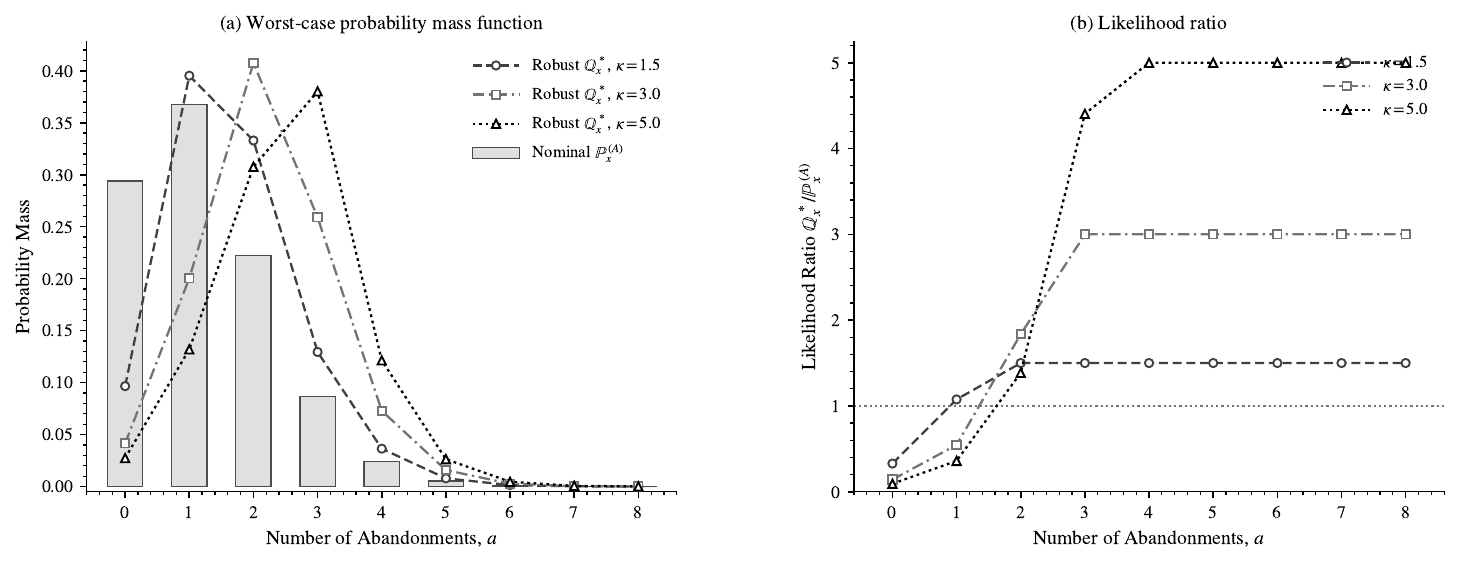}
    \caption{Adversarial distortion of the abandonment distribution at $x=30$ with $p=0.04$ and $\eta=-0.3$. Left: nominal PMF $\mathbb{P}^{(A)}_x$ versus the abandonment marginals $\mathbb{Q}^{*,(A)}_x$ of the robust joint laws for $\kappa\in\{1.5,3.0,5.0\}$. Right: likelihood ratio $d\mathbb{Q}^{*,(A)}_x/d\mathbb{P}^{(A)}_x$ of the abandonment marginal, obtained by averaging the joint tilt \eqref{eq:ra-gibbs} over the arrival atom, with flat segments where the cap $\kappa$ binds.}
    \label{fig:robust_pmf}
\end{figure}

Figure~\ref{fig:robust_pmf} visualizes the adversary's optimal distortion of the abandonment distribution at the representative state $x=30$, with $p=0.04$ and $\eta=-0.3$. The left panel compares the nominal abandonment PMF $\mathbb{P}^{(A)}_x$ with the abandonment marginals $\mathbb{Q}^{*,(A)}_x$ of the robust joint laws $\mathbb{Q}^*_x$ for likelihood-ratio caps $\kappa\in\{1.5,3.0,5.0\}$. As $\kappa$ increases, the adversary shifts more mass toward the right tail, increasing the weight on severe abandonment events and lowering the continuation value.

The right panel plots the likelihood ratio of the abandonment marginal, $d\mathbb{Q}^{*,(A)}_x/d\mathbb{P}^{(A)}_x$, against the abandonment count $a$; it is obtained by averaging the joint truncated tilt \eqref{eq:ra-gibbs} over the arrival atom $d$. The profile has two phases: an exponentially tilted region and a capped region. High-abandonment atoms are scaled up until the likelihood ratio reaches $\kappa$, after which the constraint binds and the ratio is flat. The active constraints at $\kappa=1.5,3.0,5.0$ show how the cap limits the pessimism that the adversary can introduce.

\subsection{Repeated-Clearing Calibration}\label{sec:num-renewal}
Table~\ref{tab:renewal} contrasts the optimal stopping and repeated-clearing optima at the baseline calibration. The repeated-clearing thresholds are computed by the discrete ascent of Theorem~\ref{thm:renewal-unimodal}. For comparability, the discounted repeated-clearing objective is reported in average-equivalent form, $(1-\delta)v^{\mathrm{ren}}_\delta(\theta)$, so that the optimal values in both columns are measured per period.

As an internal consistency check, the analytical mean cycle lengths agree with the simulated means in Figure~\ref{fig:operational}(b): $t(43)=76.1$ and $t(64)=155.0$, compared with simulated means of approximately $76$ and $155$, respectively. Using the optimal stopping threshold rather than the corresponding repeated-clearing optimum sacrifices $47.6\%$ of the attainable average reward and $18.3\%$ of the average-equivalent discounted repeated-clearing value at $\delta=0.99$.

\begin{table}[htbp]
\centering
\caption{Optimal stopping versus repeated-clearing optima at the baseline calibration. The discounted objective is scaled by $1-\delta$ to obtain an average-equivalent per-period value.}
\label{tab:renewal}
\begin{tabular}{lcc}
\toprule
 & Average reward ($\delta=1$) & Discounted ($\delta=0.99$) \\
\midrule
Optimal stopping threshold            & $73$   & $43$   \\
Optimal repeated-clearing threshold   & $35$   & $28$   \\
Optimal value (per period)            & $4.11$ & $3.17$ \\
Loss at the optimal stopping threshold & $47.6\%$ & $18.3\%$ \\
\bottomrule
\end{tabular}
\end{table}

\section{Proof of Theorem~\ref{thm:RN}}\label{sec:RN}
This section proves Theorem~\ref{thm:RN}. The fixed-point form of the recurrences (Remark~\ref{rem:f-rearrange}) yields the monotonicity and curvature properties of $f$ and $g$ (Lemmas~\ref{lem:monotone-increasing}--\ref{lem:strict-g-log-convex});
Lemma~\ref{lem:f-express-g} relates the two sequences explicitly, and Lemma~\ref{lem:verification} provides the Bellman verification argument. Throughout this section, $\mathcal{C}:=c+\delta\gamma p>0$ and $Y_x:=x-A_x\sim\operatorname{Binomial}(x,1-p)$ (Remark~\ref{rem:f-rearrange}). Proofs of the supporting results are deferred to Section~\ref{app:RN}.

\begin{lemma}\label{lem:monotone-increasing}
The functions $f$ and $g$ are nondecreasing on $\mathbb{Z}_{+}$.
\end{lemma}

\begin{lemma}[Curvature of the cost-accumulation function]\label{lem:strict-f-convex}
The cost-accumulation function $f$ is strictly convex on $\mathbb{Z}_{+}$: $\Delta^2 f(x) > 0$ for all $x$.
\end{lemma}

\begin{lemma}[Curvature of the discount-accumulation function]\label{lem:strict-g-log-convex}
The discount-accumulation function $g$ is convex on $\mathbb{Z}_{+}$: $\Delta^2 g(x) \ge 0$ for all $x$. Moreover, for $\delta \in (0,1)$, $g$ is strictly log-convex: $g(x+1)^2 < g(x)\, g(x+2)$ for all $x \in \mathbb{Z}_{+}$.
\end{lemma}

\begin{lemma}\label{lem:f-express-g}
Let $\delta \in (0,1)$. Then $f$ admits the explicit representation
\[
f(x) = K g(x) - \bar c x - K, \quad
\bar c := \frac{\mathcal{C}}{1 - \delta(1-p)},
\]
\[
K := \frac{\delta q\, \mathcal{C}}{(1-\delta)\big(1 - \delta(1-p)\big)}.
\]
In particular, $K = C_\infty$, the never-stop cost of Theorem~\ref{thm:RN}(ii).
\end{lemma}

\begin{lemma}[Verification]\label{lem:verification}
Suppose that $\theta^*:=\theta^*_{rn}$ is the smallest maximizer of $v_{rn}$, and let $v^*:=v_{rn}(\theta^*)>0$. Set $U(x):=g(x)v^*+f(x)$ and define
\[
\widehat V(x):=
\begin{cases}
U(x), & x\le\theta^*,\\
S(x), & x\ge\theta^*.
\end{cases}
\]
Then $\widehat V$ solves the Bellman equation~\eqref{eq:bellman-equation}, and $\widehat V\ge S$ on $\mathcal{S}$, with equality exactly for $x\ge\theta^*$.
\end{lemma}

\begin{proof}[Proof of Theorem~\ref{thm:RN}]
We first prove part~(ii), since part~(i) uses its tail behavior.

(ii) Suppose $\delta<1$. By Remark~\ref{rem:fg-meaning},
\[
v_{rn}(\theta)
=(\beta\theta-\alpha)\mathbb{E}_0[\delta^{T_\theta}]
-\mathcal{C}\,
\mathbb{E}_0\!\left[\sum_{t=0}^{T_\theta-1}\delta^tX_t\right].
\]
Since $T_\theta\ge\theta$, the first term vanishes as $\theta\to\infty$. The second expectation increases to
\[
C_\infty
=\mathcal{C}\sum_{t\ge0}\delta^t\mathbb{E}_0[X_t].
\]
Taking expectations in \eqref{eq:dynamics} gives
\[
\mathbb{E}_0[X_t]=\frac{q}{p}\bigl(1-(1-p)^t\bigr),
\]
and hence
\begin{equation}
\begin{aligned}
C_\infty
&=\frac{\mathcal{C}q}{p}
\left(
\frac{1}{1-\delta}
-\frac{1}{1-\delta(1-p)}
\right)=\frac{\delta q\,\mathcal{C}}
{(1-\delta)(1-\delta(1-p))}.
\end{aligned}
\end{equation}
Thus $v_{rn}(\theta)\to-C_\infty$.

Now let $\delta=1$. Then $v_{rn}=\beta\theta-\alpha-f(\theta)$, so it suffices to show that $\Delta f(\theta)\to\infty$. Since $\Delta f$ is strictly increasing, $L:=\lim_{\theta\to\infty}\Delta f(\theta)$ exists. If $L<\infty$, then $f(x)/x\to L$. Dividing \eqref{eq:f-fixed} by $x$ and letting $x\to\infty$ gives
\[
L=-\mathcal{C}+(1-p)L,
\]
because $\hat{X}_x/x\to1-p$ almost surely and $f$ has linear growth. Hence $pL=-\mathcal{C}<0$, contradicting $L\ge0$. Therefore $L=\infty$ and $v_{rn}(\theta)\to-\infty$.

(i) Write
\[
v_{rn}(\theta)=\frac{N_{rn}(\theta)}{g(\theta)}.
\]
By Lemma~\ref{lem:strict-f-convex}, $N_{rn}$ is strictly concave, so $\Theta_+$ is an interval of integers. Moreover, $N_{rn}(0)=-\alpha<0$, so $\underline\theta:=\min\Theta_+\ge1$. Part~(ii) implies that $\Theta_+$ is finite: for $\delta<1$, $v_{rn}(\theta)\to-C_\infty<0$ and $g(\theta)\ge\delta^{-\theta}$, so $N_{rn}=v_{rn}g\to-\infty$; for $\delta=1$, $N_{rn}=v_{rn}\to-\infty$. Write $\Theta_+=\{\underline\theta,\dots,\overline\theta\}$.

For $\theta,\theta+1\in\Theta_+$,
\begin{equation}\label{eq:rn-sign-factorization}
\operatorname{sgn}\Delta v_{rn}(\theta)
=\operatorname{sgn}\bigl(R(\theta)-\Upsilon(\theta)\bigr),
\end{equation}
where
\[
R(\theta):=\frac{N_{rn}(\theta+1)}{N_{rn}(\theta)},
\qquad
\Upsilon(\theta):=\frac{g(\theta+1)}{g(\theta)}.
\]
Strict concavity and positivity of $N_{rn}$ imply
\[
N_{rn}(\theta)^2>N_{rn}(\theta-1)N_{rn}(\theta+1),
\]
so $R$ is strictly decreasing. By Lemma~\ref{lem:strict-g-log-convex}, $\Upsilon$ is strictly increasing when $\delta<1$ and is identically $1$ when $\delta=1$. Thus $R-\Upsilon$ is strictly decreasing, so $\Delta v_{rn}$ changes sign at most once on $\Theta_+$, from positive to negative.

A crossing occurs: $\Delta v_{rn}(\underline\theta-1)>0$ and $\Delta v_{rn}(\overline\theta)<0$. Finally, $\Delta v_{rn}(\theta)>0$ for every $\theta<\underline\theta$: the increments of $N_{rn}$ are positive there, so if $N_{rn}(\theta+1)>0$, then $v_{rn}$ changes from nonpositive to positive; otherwise both values are nonpositive, while $|N_{rn}|$ decreases and $g$ increases. This proves the strict unimodality statement, with a two-point tie exactly when $\Delta v_{rn}(\theta^u)=0$.

(iii) Since $\Delta v_{rn}(\overline\theta)<0$, the set in \eqref{eq:optimal-theta} is nonempty. Part~(i) then identifies $\theta^*_{rn}$ with $\theta^u$, the smallest maximizer of $v_{rn}$, and the only other possible maximizer is $\theta^*_{rn}+1$.

Let $\theta^*:=\theta^*_{rn}$ and $\widehat V:=V_{rn,\pi_{\theta^*}}$. By Proposition~\ref{prop:value-decomposition-and-objective}, $\widehat V$ is precisely the function in Lemma~\ref{lem:verification}. Hence $\widehat V$ solves \eqref{eq:bellman-equation}, dominates $S$, and equals $S$ exactly on $\{x\ge\theta^*\}$. For $\delta<1$, $\widehat V\in\mathcal{B}_\rho$, so uniqueness in Proposition~\ref{prop:rn-bellman}(i) gives $\widehat V=V_{rn}$. For $\delta=1$, minimality in Proposition~\ref{prop:rn-bellman}(ii) gives $\widehat V\ge V_{rn}$, while the policy-value interpretation gives $\widehat V\le V_{rn}$; hence equality holds in both cases. The remaining claims follow from $g(0)=1$ and $f(0)=0$.
\end{proof}

\section{Proof of Theorem~\ref{thm:RS}}\label{sec:RS}

This section proves Theorem~\ref{thm:RS}. Proposition~\ref{prop:rs-forward} shows that the recursion \eqref{eq:fr-recurrence} is well posed and gives its exponentiated linear form; Proposition~\ref{prop:rs-ce} provides the certainty-equivalent representation of $f_{rs}$; and Proposition~\ref{prop:fr-strict-convex} establishes $\Delta^2f_{rs}\ge c$. These results reduce the optimization of $v_{rs}$ to the first-order condition \eqref{eq:optimal-theta-rs} and yield $\theta^*_{rs}\le\lceil\beta/c\rceil$.

Throughout this section, $\eta>0$ is fixed. Write $T_\theta:=\inf\{t\ge0:X_t=\theta\}$, and let $C_{x\to y}$ denote the first-passage cost from $x$ to $y$. The always-continue chain is irreducible and has drift $q-px<0$ for $x>q/p$; hence it is positive recurrent, so $T_y$ and $C_{x\to y}$ are almost surely finite. Proofs of the supporting results are deferred to Section~\ref{app:RS}.

\begin{proposition}[Well-posedness and exponentiated form of the recursion]\label{prop:rs-forward}
Let $H(x) := e^{\eta f_{rs}(x)}$, where $f_{rs}$ is defined by $f_{rs}(0) = 0$ and~\eqref{eq:fr-recurrence}.
\begin{enumerate}
    \item[(i)] The recurrence \eqref{eq:fr-recurrence} is well posed: for every $x\ge0$, write
    \[
    \underline{\omega}_{rs}(x)
    :=e^{\eta cx}-1+(1-p)^xq.
    \]
    Then
    \begin{equation*}
    H(x)\omega_{rs}(x)-F_{rs}(x)
    \ge H(x)\underline{\omega}_{rs}(x)>0.
    \end{equation*}
    Moreover, $f_{rs}(1) = 0$, and $H(x+1) \ge H(x) \ge 1$ for all $x$; equivalently, $f_{rs}$ is nonnegative and nondecreasing with $\Delta f_{rs}(0) = 0$.
    \item[(ii)] $H$ satisfies the exponentiated linear recursion
    \begin{equation*}
    H(x) = e^{-\eta c x}\, \mathbb{E}_{A_x, D}\!\left[ e^{-\eta \gamma A_x}\, H(\hat{X}_x) \right], \qquad x \ge 0.
    \end{equation*}
    \item[(iii)] (Multiplicative decomposition) For every $\theta \in \mathbb{Z}_{+}$ and every $x \in \{0, 1, \dots, \theta\}$,
    \begin{equation*}
    e^{\eta V_{rs,\pi_\theta}(x)} = H(x)\, e^{\eta V_{rs,\pi_\theta}(0)}.
    \end{equation*}
\end{enumerate}
\end{proposition}

\begin{proposition}[Certainty-equivalent representation]\label{prop:rs-ce}
For every $\theta \ge 0$,
\begin{equation}\label{eq:dfr-ce}
\Delta f_{rs}(\theta) = -\frac{1}{\eta} \log \mathbb{E}\!\left[ e^{-\eta C_{\theta \to \theta + 1}} \right].
\end{equation}
Consequently:
\begin{enumerate}
    \item[(i)] $f_{rs}(1) = 0$ and $0 \le f_{rs}(x) \le f(x)$ for all $x$, where $f$ is the risk-neutral cost-accumulation function at $\delta = 1$;
    \item[(ii)] $\Delta f_{rs}(\theta) \le \Delta f(\theta)$ for all $\theta$;
    \item[(iii)] Under Assumption~\ref{ass:profitability}, $\theta^*_{rs} \ge \theta^*_{rn}$: the risk-seeking threshold is no smaller than the undiscounted risk-neutral threshold of Theorem~\ref{thm:RN}.
\end{enumerate}
\end{proposition}

\begin{proposition}[Strict convexity of $f_{rs}$]\label{prop:fr-strict-convex}
For all $\theta \ge 0$, $\Delta^2 f_{rs}(\theta) \ge c$. In particular, $f_{rs}$ is strictly convex, $\Delta f_{rs}(\theta) \ge c \theta$, and $f_{rs}(\theta) \ge c \, \theta (\theta - 1) / 2$.
\end{proposition}

\begin{proof}[Proof of Theorem~\ref{thm:RS}]
By Proposition~\ref{prop:rs-threshold-wlog}, some threshold policy is optimal among all admissible policies, so it suffices to maximize $v_{rs}(\theta) = \beta \theta - \alpha - f_{rs}(\theta)$ over $\theta \in \mathbb{Z}_{+}$.

(i) Since $\Delta v_{rs}(\theta) = \beta - \Delta f_{rs}(\theta)$, Proposition~\ref{prop:fr-strict-convex} gives $\Delta^2 v_{rs}(\theta) = -\Delta^2 f_{rs}(\theta) \le -c$ for every $\theta \ge 0$: $v_{rs}$ is strictly discretely concave and $\Delta v_{rs}$ is strictly decreasing. In particular, $v_{rs}(\theta) \le \beta \theta - \alpha - c\,\theta(\theta-1)/2$, so $v_{rs}(\theta) \to -\infty$ at least quadratically as $\theta \to \infty$.

(ii) The first-order condition $\Delta v_{rs}(\theta) \le 0$ is equivalent to $\Delta f_{rs}(\theta) \ge \beta$. By Proposition~\ref{prop:fr-strict-convex}, $\Delta f_{rs}$ is strictly increasing with $\Delta f_{rs}(\theta) \ge c\,\theta$, so $\Delta f_{rs}(\theta) \ge \beta$ for every $\theta \ge \lceil \beta / c \rceil$: the set in \eqref{eq:optimal-theta-rs} is nonempty and its minimum satisfies $\theta^*_{rs} \le \lceil \beta / c \rceil$, while $\Delta f_{rs}(0) = 0 < \beta$ (Proposition~\ref{prop:rs-forward}(i)) gives $\theta^*_{rs} \ge 1$. Because $\Delta v_{rs}$ is strictly decreasing by part~(i), the smallest maximizer of $v_{rs}$ is the first index at which the forward difference is nonpositive:
\begin{equation*}
\begin{aligned}
\theta^*_{rs}
&=\min\bigl\{\theta:\Delta v_{rs}(\theta)\le0\bigr\}=\min\bigl\{\theta:\Delta f_{rs}(\theta)\ge\beta\bigr\}.
\end{aligned}
\end{equation*}

(iii) If $\Delta f_{rs}(\theta^*_{rs}) = \beta$, then $\Delta v_{rs}(\theta^*_{rs}) = 0$, so $v_{rs}(\theta^*_{rs}) = v_{rs}(\theta^*_{rs} + 1)$ and both thresholds are optimal, while strict concavity precludes any third optimizer; otherwise the maximizer is unique. The remaining claims restate Proposition~\ref{prop:rs-threshold-wlog}: $V_{rs}(0) = v_{rs}(\theta^*_{rs}) = \max_{\theta \ge 0} v_{rs}(\theta)$, and $\pi_{\theta^*_{rs}}$ is optimal among all admissible policies.
\end{proof}

\section{Proof of Theorem~\ref{thm:RA}}\label{sec:RA}

This section proves Theorem~\ref{thm:RA}: Lemma~\ref{lem:ra-drift} establishes a uniform Foster--Lyapunov drift condition over the adversary's measures,
Lemmas~\ref{lem:ra-concavity}--\ref{lem:ra-fullsupport} and Proposition~\ref{prop:ra-wellposed} give attainment, full support, and monotone value iteration, Lemma~\ref{lem:ra-increments} and Corollary~\ref{cor:ra-threshold-wlog} yield the threshold structure,
Lemma~\ref{lem:ra-uniqueness} gives uniqueness, and Proposition~\ref{prop:ra-verification} identifies the fixed point with the robust game value. A finite-dimensional KKT argument then gives Theorem~\ref{thm:RA}(iii), from which Proposition~\ref{prop:kappa-limit} follows. Proofs of the supporting lemmas and of Proposition~\ref{prop:ra-verification} are deferred to Section~\ref{app:RA}.

Throughout, $\eta<0$ is fixed and $w(x)=x+L$, with $L>0$ chosen below. Write
$
\mathcal{K}_\eta(\mathbb{Q}\,\|\,\mathbb{P}) \;:=\; -\frac{1}{\eta}\, D_{\mathrm{KL}}(\mathbb{Q}\,\|\,\mathbb{P}) \;\ge\; 0
$
for the KL penalty at risk parameter $\eta$ (recall $-\eta^{-1} > 0$ for $\eta < 0$). For $u\ge S$ and $\mathbb{Q}_x\in\mathcal{Q}_\kappa(x)$, define the one-step objective
\begin{equation}\label{eq:ra-fu}
f_u(\mathbb{Q}_x)
:=\mathbb{E}_{\mathbb{Q}_x}\!\left[
u(\hat{X}_x)-\gamma A_x
\right]
+\mathcal{K}_\eta(\mathbb{Q}_x\,\|\,\mathbb{P}_x).
\end{equation}
The robust Bellman operator is
\begin{equation}\label{eq:ra-operator}
(\mathcal{T}_{ra}u)(x)
=\max\!\left\{
S(x),\,
-cx+\inf_{\mathbb{Q}_x\in\mathcal{Q}_\kappa(x)}
f_u(\mathbb{Q}_x)
\right\}.
\end{equation}
At $\delta=1$, the corresponding risk-neutral operator is
\[
(\mathcal{T}_{rn}u)(x)
=\max\!\left\{
S(x),\,
-\mathcal{C}x+\mathbb{E}[u(\hat{X}_x)]
\right\}.
\]
Since $\mathbb{P}_x\in\mathcal{Q}_\kappa(x)$ and $\mathcal{K}_\eta(\mathbb{P}_x\,\|\,\mathbb{P}_x)=0$,
\begin{equation}\label{eq:ra-le-rn}
\mathcal{T}_{ra}u\le\mathcal{T}_{rn}u,
\qquad u\ge S.
\end{equation}

\begin{lemma}[Uniform drift under the adversary]\label{lem:ra-drift}
There is a constant $C_\kappa > 0$ such that, for every $x \in \mathbb{Z}_{+}$ and every $\mathbb{Q}_x \in \mathcal{Q}_\kappa(x)$,
\begin{equation*}
\mathbb{E}_{\mathbb{Q}_x}[A_x] \ge p x - C_\kappa \sqrt{x}.
\end{equation*}
Consequently, for every $\xi \in (1 - p, 1)$ there exist $M \in \mathbb{Z}_{+}$ and $b > 0$ such that, uniformly over $\mathbb{Q}_x \in \mathcal{Q}_\kappa(x)$,
\begin{equation}\label{eq:ra-drift}
\mathbb{E}_{\mathbb{Q}_x}\!\left[ w(X_1) \right] \le \xi w(x) + b \, \mathbf{1}_{\{x \le M\}}, \qquad x \in \mathbb{Z}_{+},
\end{equation}
and iterating \eqref{eq:ra-drift} gives $\mathbb{E}_{\mathbb{Q}}\!\left[ w(X_m) \right] \le \xi^m w(x) + b/(1-\xi)$ under any history-dependent adversarial strategy.
\end{lemma}

\begin{lemma}[Attainment of the one-step robust infimum]\label{lem:ra-concavity}
Let $\mathbb{P}$ be a probability measure with full support on a finite set $\mathcal{X}$, and let $\kappa > 1$. The one-step robust operator
\[
\mathcal{R}(Z)=\inf_{\mathbb{Q}\in\mathcal{Q}_\kappa}
\Bigl\{\mathbb{E}_{\mathbb{Q}}[Z]
-\frac{1}{\eta}D_{\mathrm{KL}}(\mathbb{Q}\,\|\,\mathbb{P})\Bigr\}, \quad Z \in \mathbb{R}^{\mathcal{X}},
\]
is nondecreasing and concave, and the infimum is attained by some $\mathbb{Q}^* \in \mathcal{Q}_\kappa$ for every $Z$.
\end{lemma}

\begin{lemma}[Full support of the worst-case measure]\label{lem:ra-fullsupport}
Let $Z$ be a bounded random variable on the finite support of $\mathbb{P}_x$. Every $\mathbb{Q}^* \in \mathcal{Q}_\kappa(x)$ attaining $\inf_{\mathbb{Q} \in \mathcal{Q}_\kappa(x)} \{ \mathbb{E}_{\mathbb{Q}}[Z] - \eta^{-1} D_{\mathrm{KL}}(\mathbb{Q} \| \mathbb{P}_x) \}$ satisfies $\mathbb{Q}^*(a,d) > 0$ at every atom $(a,d)$ with $\mathbb{P}_x(a,d) > 0$. In particular, the Markov kernel $\mathbb{Q}^*_x(y) := \mathbb{Q}^*(\{(a,d) : x - a + d = y\})$ charges every state $y \in \{0, 1, \dots, x+1\}$.
\end{lemma}

\begin{proposition}[Well-posedness: existence and value iteration]\label{prop:ra-wellposed}
Let $\eta < 0$ and $\kappa > 1$, and suppose Assumption~\ref{ass:profitability} holds.
\begin{enumerate}
    \item[(i)] The value iterates $V_0 = S$, $V_{n+1} = \mathcal{T}_{ra} V_n$ form a nondecreasing sequence bounded above by the risk-neutral value function $V_{rn}$ of \eqref{eq:opt-problem-rn} at $\delta = 1$, and converge pointwise to the regularized value function $V_{ra}^{\kappa}$, a fixed point of $\mathcal{T}_{ra}$ with $S \le V \le V_{rn}$.

    \item[(ii)] $V$ is the minimal fixed point of $\mathcal{T}_{ra}$ dominating $S$: every fixed point $\tilde V$ of $\mathcal{T}_{ra}$ with $\tilde V \ge S$ satisfies $\tilde V \ge V$.
    \item[(iii)] For every $u_0$ with $S \le u_0 \le V$, the iterates $\mathcal{T}_{ra}^n u_0$ converge pointwise to $V$.
\end{enumerate}
\end{proposition}

\begin{proof}[Proof of Proposition~\ref{prop:ra-wellposed}]
(i) $\mathcal{T}_{ra}$ is order preserving, and $\mathcal{T}_{ra} S \ge S$ since the outer maximum includes $S$; hence $V_n = \mathcal{T}_{ra}^n S$ is nondecreasing. Moreover, $V_0 = S \le V_{rn}$ (stopping immediately is admissible in the risk-neutral problem), and if $V_n \le V_{rn}$, then by \eqref{eq:ra-le-rn} and monotonicity of both operators, $V_{n+1} = \mathcal{T}_{ra} V_n \le \mathcal{T}_{rn} V_n \le \mathcal{T}_{rn} V_{rn} = V_{rn}$. Hence $V_n \uparrow V$ pointwise for some $V$ with $S \le V \le V_{rn}$. Note that $V_{rn}$ is finite-valued with $V_{rn}(x) \le \beta x + K$ for some $K$: at $\delta = 1$ the risk-neutral optimal stopping region is $\{\theta^*_{rn}, \theta^*_{rn} + 1, \dots\}$ with $\theta^*_{rn} < \infty$ by Theorem~\ref{thm:RN}, so $V_{rn}$ equals $S$ above $\theta^*_{rn}$ and is finite on the finitely many states below it.

It remains to show $V = \mathcal{T}_{ra} V$.
Fix $x$ and recall $f_u$ from \eqref{eq:ra-fu}, so that $(\mathcal{T}_{ra} u)(x) = \max\{ S(x),\, -cx + \inf_{\mathbb{Q}_x \in \mathcal{Q}_\kappa(x)} f_u(\mathbb{Q}_x) \}$ for $u \ge S$. Since $V_n \uparrow V$, monotone convergence gives $f_{V_n}(\mathbb{Q}_x) \uparrow f_V(\mathbb{Q}_x)$ for every $\mathbb{Q}_x$ (the abandonment term is bounded by $\gamma x$). Consequently, $\inf_{\mathbb{Q}_x} f_{V_n}(\mathbb{Q}_x) \uparrow \inf_{\mathbb{Q}_x} f_V(\mathbb{Q}_x)$: the inequality $\le$ is immediate, and for the reverse, let $\mathbb{Q}_n$ attain $\inf f_{V_n}$ (Lemma~\ref{lem:ra-concavity}); then
\[
\begin{aligned}
\inf f_V-\inf f_{V_n}\le f_V(\mathbb{Q}_n)-f_{V_n}(\mathbb{Q}_n)
&=\mathbb{E}_{\mathbb{Q}_n}\!\left[
(V-V_n)(\hat{X}_x)
\right]\\
&\le\kappa\,\mathbb{E}_{\mathbb{P}_x}\!\left[
(V-V_n)(\hat{X}_x)
\right]
\to0.
\end{aligned}
\]
Here dominated convergence applies because the differences $(V - V_n)(\hat{X}_x)$ decrease pointwise to $0$ and are bounded by the $\mathbb{P}_x$-integrable envelope $(V - S)(\hat{X}_x)$, since $V_n \ge S$ and
\[
\mathbb{E}_{\mathbb{P}_x}\!\left[ (V - S)(\hat{X}_x) \right]
\le \beta (x+1) + K + \alpha < \infty;
\]
the one-step support $\{0, \dots, x+1\}$ of $\mathbb{P}_x$ is finite, so the envelope is integrable and therefore uniformly integrable. Taking the maximum with $S(x)$ gives $\mathcal{T}_{ra} V_n \uparrow \mathcal{T}_{ra} V$; but $\mathcal{T}_{ra} V_n = V_{n+1} \uparrow V$, so $\mathcal{T}_{ra} V = V$.

(ii) If $\tilde V = \mathcal{T}_{ra} \tilde V \ge S$, then by monotonicity $\tilde V = \mathcal{T}_{ra}^n \tilde V \ge \mathcal{T}_{ra}^n S = V_n$ for every $n$, and letting $n \to \infty$ gives $\tilde V \ge V$.

(iii) If $S \le u_0 \le V$, then $\mathcal{T}_{ra}^n S \le \mathcal{T}_{ra}^n u_0 \le \mathcal{T}_{ra}^n V = V$; the left side increases to $V$ by (i), so $\mathcal{T}_{ra}^n u_0 \to V$.
\end{proof}

\begin{lemma}[Increment bound]\label{lem:ra-increments}
Let $V = V_{ra}^{\kappa}$ be the fixed point of Proposition~\ref{prop:ra-wellposed}. Then $\Delta V(x) := V(x+1) - V(x) \le \beta$ for all $x$, and the continuation branch satisfies the strict decrement bound $\Delta Q_{ra}(x) \le \beta - c - \gamma p$ for all $x$.
\end{lemma}

\begin{corollary}[Structure of the stopping region]\label{cor:ra-threshold-wlog}
Let $V = V_{ra}^{\kappa}$ and $\theta^* = \min\{x : V(x) = S(x)\}$.
\begin{enumerate}
    \item[(i)] (Nonempty) The stopping region $\{x : V(x) = S(x)\}$ is nonempty, so $\theta^* < \infty$.
    \item[(ii)] (Up-closed) $V(x) = S(x)$ for all $x \ge \theta^*$.
    \item[(iii)] (Strict inequalities) Continuation is strictly optimal below $\theta^*$: $Q_{ra}(x) > S(x)$ for $x < \theta^*$; and stopping is strictly optimal above $\theta^*$: $S(x) > Q_{ra}(x)$ for $x > \theta^*$.
\end{enumerate}
\end{corollary}

\begin{lemma}[Uniqueness of the fixed point]\label{lem:ra-uniqueness}
Let $V$ be the fixed point of Proposition~\ref{prop:ra-wellposed}, and let $\tilde V$ be any fixed point of $\mathcal{T}_{ra}$ with $\tilde V \ge S$ and $\tilde V(x) \le \beta x + K$ for some constant $K$. Then $\tilde V = V$.
\end{lemma}

\begin{proposition}[Game value and threshold optimality]\label{prop:ra-verification}
Let $V$ be the fixed point of Proposition~\ref{prop:ra-wellposed} and $\theta^* = \min\{x : V(x) = S(x)\}$. Define an \textit{adversarial strategy} $\mathbb{Q} = (\mathbb{Q}_t)_{t \ge 0}$ as a sequence of one-step kernels $\mathbb{Q}_t \in \mathcal{Q}_\kappa(X_t)$ adapted to $\mathcal{F}_t = \sigma(X_0, \dots, X_t)$, and define the per-step payoff
\begin{equation*}
r(X_t, \mathbb{Q}_t) := -c X_t - \gamma\, \mathbb{E}_{\mathbb{Q}_t}[A] - \frac{1}{\eta} D_{\mathrm{KL}}(\mathbb{Q}_t \| \mathbb{P}_{X_t}),
\end{equation*}
the truncated reward $R_N(\tau) := \sum_{t = 0}^{\tau \wedge N - 1} r(X_t, \mathbb{Q}_t) + S(X_{\tau \wedge N})$, and the robust policy values
\[
\begin{aligned}
J(x,\tau,\mathbb{Q})
&:=\liminf_{N\to\infty}
\mathbb{E}^{\mathbb{Q}}_x[R_N(\tau)]
\in[-\infty,\infty],\quad
\underline{J}(x,\tau):=\inf_{\mathbb{Q}}J(x,\tau,\mathbb{Q})
\end{aligned}
\]
for $\tau \in \mathfrak{T}$ (almost surely finite stopping times). Then the following hold.
\begin{enumerate}
    \item[(i)] $V(x) \ge \underline{J}(x, \tau)$ for every $\tau \in \mathfrak{T}$.
    \item[(ii)] $V(x) = \sup_{\tau \in \mathfrak{T}} \underline{J}(x, \tau)$.
    \item[(iii)] $\underline{J}(x, \pi_{\theta^*}) = V(x)$, where $\pi_{\theta^*}$ is the threshold policy with stopping time $\tau_{\theta^*} = \inf\{ t : X_t \ge \theta^* \}$. If the tie case $Q_{ra}(\theta^*) = S(\theta^*)$ holds, then also $\underline{J}(x, \pi_{\theta^* + 1}) = V(x)$, and these are the only optimal threshold policies.
\end{enumerate}
\end{proposition}

\begin{proof}[Proof of Theorem~\ref{thm:RA}]
(i) Existence of a fixed point $V$ with $S \le V \le V_{rn}$, its minimality above $S$, and the convergence of the value iterates from $S$ follow from Proposition~\ref{prop:ra-wellposed}; uniqueness in the stated class follows from Lemma~\ref{lem:ra-uniqueness}.

(ii) The threshold structure (nonempty, up-closed stopping region, strict inequalities) follows from Corollary~\ref{cor:ra-threshold-wlog}; the game-value identity and the optimality of $\pi_{\theta^*_{ra}}$ (and of $\pi_{\theta^*_{ra}+1}$ in the tie case) follow from Proposition~\ref{prop:ra-verification}.

(iii) Fix $x$, enumerate the finitely many atoms $i=(a,d)$ of $\mathbb{P}_x$, and write $p_i=\mathbb{P}_x(i)$, $\varphi_i=V(x-a+d)-\gamma a$, and $k=-1/\eta>0$. In terms of the likelihood ratio $m=d\mathbb{Q}/d\mathbb{P}_x$, the one-step minimization is
\[
\begin{aligned}
\min_{\{m_i\}} \quad
& \sum_i p_i \left( m_i \varphi_i + k\, m_i \log m_i \right) \\
\text{s.t.} \quad
& \sum_i p_i m_i = 1, \\
& 0 \le m_i \le \kappa, \qquad \forall i.
\end{aligned}
\]
The objective is strictly convex, so the minimizer $m^*$ is unique; Lemma~\ref{lem:ra-fullsupport} gives $m_i^*>0$ at every atom. With $\mu_x$ the multiplier for the normalization constraint and $\nu_i\ge0$ the multiplier for $m_i\le\kappa$, the KKT conditions are
\[
\varphi_i+k(\log m_i^*+1)+\mu_x+\nu_i=0,
\quad
\nu_i(m_i^*-\kappa)=0.
\]
Set $\lambda_x:=e^{\eta\mu_x-1}>0$. If $\nu_i=0$, stationarity gives $m_i^*=\lambda_x e^{\eta \varphi_i}$; if $\lambda_x e^{\eta \varphi_i}>\kappa$, the upper constraint is active, $m_i^*=\kappa$, and $\nu_i=k\log\bigl(\lambda_x e^{\eta \varphi_i}/\kappa\bigr)>0$. Thus, $m_i^*=\min\bigl\{\kappa,\lambda_x e^{\eta \varphi_i}\bigr\}$, where $\lambda_x$ is determined by the normalization $\sum_i p_i\min\bigl\{\kappa,\lambda_x e^{\eta \varphi_i}\bigr\}=1$. A unique such $\lambda_x$ exists because the left-hand side is continuous, tends to $0$ as $\lambda_x\downarrow0$ and to $\kappa>1$ as $\lambda_x\uparrow\infty$, and is strictly increasing before all upper constraints become active. Hence, $m_i^*>0$ at every atom, so $\mathbb{Q}_x^*$ has full support.
\end{proof}

\section{Concluding Remarks}\label{sec:conclusion}
We have developed a methodology for optimal batch clearing with item-level abandonment, combining threshold verification with forward first-passage recursions. The approach explicitly characterizes optimal policies under risk-neutral, risk-seeking, and regularized risk-averse criteria, as well as under repeated-clearing objectives, and identifies the worst-case transition law in truncated-Gibbs form. More broadly, the analysis shows that attrition couples batch size to cycle duration, making delay, risk posture, and repeated operation integral to optimal batching decisions. 

Our results carry several implications for the management of batch clearing systems. First, the optimal clearing rule is simple to implement: process the batch as soon as the inventory reaches a directly computable threshold. The threshold moves in the intuitive directions---it increases with the stopping cost and the arrival probability and decreases with the holding cost, the abandonment penalty, and the abandonment probability---so attrition-prone environments call for smaller, more frequent batches.
Second, at the optimum it is time, not cost, that binds: the nominal value of the batch is still increasing when clearing is triggered, and what forces the decision is the erosion of that value while waiting. Because abandonment can slow accumulation progressively as the buffer grows, over-batching is disproportionately expensive---cycle lengths grow superlinearly in the threshold---so when in doubt, managers should favor smaller batches.
Third, risk posture enters the decision as an explicit, priced adjustment: risk seeking deflates the perceived marginal cost of accumulation and justifies a larger batch---the risk-seeking premium---whereas risk aversion acts as insurance against abandonment-heavy distortions of the model, with the likelihood-ratio cap bounding how pessimistic the hedge may become and the resulting risk-averse discount quantifying the price of that protection. 

Several directions for future research follow. On the modeling side, positive service times, capacity constraints, partial batches, compound arrivals, multiple item classes under joint capacity constraints, and nonstationary arrival dynamics are natural extensions, with applications to perishable-goods logistics and cloud task aggregation. On the methodological side, the discounted entropic criterion requires an augmented-state analysis; a discounted recursive entropic criterion is different again and is treated separately in Section~\ref{app:discounted-rs}. 
A risk-sensitive renewal theory is a further open direction: under an entropic criterion, cycle rewards combine multiplicatively, so the additive renewal--reward reduction does not apply and an ergodic multiplicative-functional or recursive-entropic formulation would be required.
The robust dual formulation also suggests a distributionally robust version with parameters learned online and the likelihood-ratio cap calibrated from data. More generally, the two-step strategy used here---a uniform increments bound to isolate the threshold family, followed by a first-passage decomposition---extends beyond the present cost structure; identifying further classes of controlled skip-free processes in which it yields exact reduced objectives may be of independent interest.

\section{Code and Data Disclosure}\label{sec:Code and Data Disclosure}
The code to support the numerical experiments in this paper can be found at \url{https://github.com/Ahkfor/Optimal-Batching-with-Abandonment.git}.

\bibliographystyle{plainnat}
\bibliography{references}

@article{ross1971infinitesimal,
  title={Infinitesimal look-ahead stopping rules},
  author={Ross, Sheldon M},
  journal={The Annals of Mathematical Statistics},
  volume={42},
  number={1},
  pages={297--303},
  year={1971},
  publisher={JSTOR}
}

@article{ross1969optimal,
  title={Optimal dispatching of a Poisson process},
  author={Ross, Sheldon M},
  journal={Journal of Applied Probability},
  volume={6},
  number={3},
  pages={692--699},
  year={1969},
  publisher={Cambridge University Press}
}

@article{whitt1981stationary,
  title={The stationary distribution of a stochastic clearing process},
  author={Whitt, Ward},
  journal={Operations Research},
  volume={29},
  number={2},
  pages={294--308},
  year={1981},
  publisher={INFORMS}
}

@article{stidham1974stochastic,
  title={Stochastic clearing systems},
  author={Stidham Jr, Shaler},
  journal={Stochastic Processes and their Applications},
  volume={2},
  number={1},
  pages={85--113},
  year={1974}
}

@article{kella2003stochastic,
  title={A stochastic clearing model with a Brownian and a compound Poisson component},
  author={Kella, Offer and Perry, David and Stadje, Wolfgang},
  journal={Probability in the Engineering and Informational Sciences},
  volume={17},
  number={1},
  pages={1--22},
  year={2003},
  publisher={Cambridge University Press}
}

@article{deb1973optimal,
  title={Optimal control of batch service queues},
  author={Deb, Rajat K and Serfozo, Richard F},
  journal={Advances in Applied Probability},
  volume={5},
  number={2},
  pages={340--361},
  year={1973},
  publisher={Cambridge University Press}
}

@article{deb1976optimal,
  title={Optimal control of batch service queues with switching costs},
  author={Deb, Rajat K},
  journal={Advances in Applied Probability},
  volume={8},
  number={1},
  pages={177--194},
  year={1976},
  publisher={Cambridge University Press}
}

@article{howard1972risk,
  title={Risk-sensitive Markov decision processes},
  author={Howard, Ronald A and Matheson, James E},
  journal={Management Science},
  volume={18},
  number={7},
  pages={356--369},
  year={1972}
}

@article{denardo1979optimal,
  title={Optimal stopping, exponential utility, and linear programming},
  author={Denardo, Eric V and Rothblum, Uriel G},
  journal={Mathematical Programming},
  volume={16},
  number={1},
  pages={228--244},
  year={1979},
  publisher={Springer}
}

@article{borkar2002risk,
  title={Risk-sensitive optimal control for Markov decision processes with monotone cost},
  author={Borkar, Vivek S and Meyn, Sean P},
  journal={Mathematics of Operations Research},
  volume={27},
  number={1},
  pages={192--209},
  year={2002},
  publisher={INFORMS}
}

@article{ruszczynski2010risk,
  title={Risk-averse dynamic programming for Markov decision processes},
  author={Ruszczy{\'n}ski, Andrzej},
  journal={Mathematical Programming},
  volume={125},
  number={2},
  pages={235--261},
  year={2010}
}

@article{pichler2022risk,
  title={Risk-averse stochastic programming: Time consistency and optimal stopping},
  author={Pichler, Alois and Liu, Rui Peng and Shapiro, Alexander},
  journal={Operations Research},
  volume={70},
  number={4},
  pages={2439--2455},
  year={2022},
  publisher={INFORMS}
}

@article{ren2025structural,
  title={On Structural Properties of Risk-Averse Optimal Stopping Problems},
  author={Ren, Xingyu and Fu, Michael C and Marcus, Steven I},
  journal={arXiv preprint arXiv:2511.01022},
  year={2025}
}

@book{shiryaev2008optimal,
  title={Optimal stopping rules},
  author={Shiryaev, Albert N},
  year={2008},
  publisher={Springer}
}

@book{peskir2006optimal,
  title={Optimal stopping and free-boundary problems},
  author={Peskir, Goran and Shiryaev, Albert},
  year={2006},
  publisher={Springer}
}

@book{nahmias2011perishable,
  title={Perishable inventory systems},
  author={Nahmias, Steven},
  year={2011},
  publisher={Springer Science \& Business Media}
}

@book{dupuis2011weak,
  title={A weak convergence approach to the theory of large deviations},
  author={Dupuis, Paul and Ellis, Richard S},
  year={2011},
  publisher={John Wiley \& Sons}
}

@article{ccetinkaya2002optimal,
  title={Optimal outbound dispatch policies: Modeling inventory and cargo capacity},
  author={{\c{C}}etinkaya, S{\i}la and Lee, Chung-Yee},
  journal={Naval Research Logistics (NRL)},
  volume={49},
  number={6},
  pages={531--556},
  year={2002},
  publisher={Wiley Online Library}
}

@article{xiong2000control,
  title={Control architecture in optical burst-switched {WDM} networks},
  author={Xiong, Yijun and Vandenhoute, Marc and Cankaya, Hakki C},
  journal={IEEE Journal on Selected Areas in Communications},
  volume={18},
  number={10},
  pages={1838--1851},
  year={2000}
}

@article{haijema2007blood,
  title={Blood platelet production: Optimization by dynamic programming and simulation},
  author={Haijema, Ren{\'e} and van der Wal, Jan and van Dijk, Nico M},
  journal={Computers \& Operations Research},
  volume={34},
  number={3},
  pages={760--779},
  year={2007},
  publisher={Elsevier}
}

@article{garnett2002designing,
  title={Designing a call center with impatient customers},
  author={Garnett, Ofer and Mandelbaum, Avishai and Reiman, Martin},
  journal={Manufacturing \& Service Operations Management},
  volume={4},
  number={3},
  pages={208--227},
  year={2002},
  publisher={INFORMS}
}

@incollection{baccelli1981queues,
  title={On queues with impatient customers},
  author={Baccelli, Fran{\c{c}}ois and H{\'e}buterne, G{\'e}rard},
  booktitle={Performance '81},
  editor={Kylstra, F. J.},
  pages={159--179},
  year={1981},
  publisher={North-Holland}
}

@article{vanleeuwen2017optimal,
  title={Optimal dispatching in a tandem queue},
  author={van Leeuwen, Daphne and N{\'u}{\~n}ez Queija, Rudesindo},
  journal={Queueing Systems},
  volume={87},
  number={3--4},
  pages={269--291},
  year={2017},
  publisher={Springer}
}

@inproceedings{larranaga2015efficient,
  title={Efficient content delivery in the presence of impatient jobs},
  author={Larra{\~n}aga, Maialen and Boxma, Onno J and N{\'u}nez-Queija, Rudesindo and Squillante, Mark S},
  booktitle={2015 27th International Teletraffic Congress},
  pages={73--81},
  year={2015},
  organization={IEEE}
}

@book{puterman2014markov,
  title={Markov decision processes: discrete stochastic dynamic programming},
  author={Puterman, Martin L},
  year={2014},
  publisher={John Wiley \& Sons}
}

@article{jelito2022risk,
  title={Risk-sensitive optimal stopping with unbounded terminal cost function},
  author={Jelito, Damian and Stettner, {\L}ukasz},
  journal={Electronic Journal of Probability},
  volume={27},
  pages={1--30},
  year={2022},
  publisher={The Institute of Mathematical Statistics and the Bernoulli Society}
}

@article{donsker1975asymptotic,
  title={Asymptotic evaluation of certain Markov process expectations for large time, I},
  author={Donsker, Monroe D and Varadhan, S.R. Srinivasa},
  journal={Communications on Pure and Applied Mathematics},
  volume={28},
  number={1},
  pages={1--47},
  year={1975},
  publisher={Wiley Online Library}
}

@article{nilim2005robust,
  title={Robust control of Markov decision processes with uncertain transition matrices},
  author={Nilim, Arnab and El Ghaoui, Laurent},
  journal={Operations Research},
  volume={53},
  number={5},
  pages={780--798},
  year={2005},
  publisher={INFORMS}
}

@article{iyengar2005robust,
  title={Robust dynamic programming},
  author={Iyengar, Garud N},
  journal={Mathematics of Operations Research},
  volume={30},
  number={2},
  pages={257--280},
  year={2005},
  publisher={INFORMS}
}

@article{petersen2000minimax,
  title={Minimax optimal control of stochastic uncertain systems with relative entropy constraints},
  author={Petersen, Ian R and James, Matthew R and Dupuis, Paul},
  journal={IEEE Transactions on Automatic Control},
  volume={45},
  number={3},
  pages={398--412},
  year={2000},
  publisher={IEEE}
}

@article{dai1996connections,
  title={Connections between stochastic control and dynamic games},
  author={Dai Pra, Paolo and Meneghini, Lorenzo and Runggaldier, Wolfgang J},
  journal={Mathematics of Control, Signals and Systems},
  volume={9},
  number={4},
  pages={303--326},
  year={1996},
  publisher={Springer}
}

@article{wiesemann2013robust,
  title={Robust Markov decision processes},
  author={Wiesemann, Wolfram and Kuhn, Daniel and Rustem, Ber{\c{c}}},
  journal={Mathematics of Operations Research},
  volume={38},
  number={1},
  pages={153--183},
  year={2013},
  publisher={INFORMS}
}

@article{james1992asymptotic,
  title={Asymptotic analysis of nonlinear stochastic risk-sensitive control and differential games},
  author={James, Matthew R},
  journal={Mathematics of Control, Signals and Systems},
  volume={5},
  number={4},
  pages={401--417},
  year={1992},
  publisher={Springer}
}

@article{chung1987discounted,
  title={Discounted {MDP}'s: distribution functions and exponential utility maximization},
  author={Chung, Kun-Jen and Sobel, Matthew J},
  journal={SIAM Journal on Control and Optimization},
  volume={25},
  number={1},
  pages={49--62},
  year={1987}
}

@article{bouakiz1992inventory,
  title={Inventory control with an exponential utility criterion},
  author={Bouakiz, Mokrane and Sobel, Matthew J},
  journal={Operations Research},
  volume={40},
  number={3},
  pages={603--608},
  year={1992},
  publisher={INFORMS}
}

@article{di1999risk,
  title={Risk-sensitive control of discrete-time Markov processes with infinite horizon},
  author={Di Masi, Giovanni B and Stettner, Lukasz},
  journal={SIAM Journal on Control and Optimization},
  volume={38},
  number={1},
  pages={61--78},
  year={1999},
  publisher={SIAM}
}

@article{bauerle2014more,
  title={More risk-sensitive Markov decision processes},
  author={B{\"a}uerle, Nicole and Rieder, Ulrich},
  journal={Mathematics of Operations Research},
  volume={39},
  number={1},
  pages={105--120},
  year={2014},
  publisher={INFORMS}
}

@article{weiss1979computation,
  title={The computation of optimal control limits for a queue with batch services},
  author={Weiss, Howard J},
  journal={Management Science},
  volume={25},
  number={4},
  pages={320--328},
  year={1979},
  publisher={INFORMS}
}

@article{xia2002optimal,
  title={Optimal control of parallel queues with batch service},
  author={Xia, Cathy H and Michailidis, George and Bambos, Nicholas and Glynn, Peter W},
  journal={Probability in the Engineering and Informational Sciences},
  volume={16},
  number={3},
  pages={289--307},
  year={2002},
  publisher={Cambridge University Press}
}

@article{deb1978optimal,
  title={Optimal dispatching of a finite capacity shuttle},
  author={Deb, Rajat K},
  journal={Management Science},
  volume={24},
  number={13},
  pages={1362--1372},
  year={1978},
  publisher={INFORMS}
}

@article{browne1992gated,
  title={The gated infinite-server queue: Uniform service times},
  author={Browne, Sid and Coffman, Jr, EG and Gilbert, EN and Wright, Paul E},
  journal={SIAM Journal on Applied Mathematics},
  volume={52},
  number={6},
  pages={1751--1762},
  year={1992},
  publisher={SIAM}
}

@article{nagai2007stopping,
  title={Stopping problems of certain multiplicative functionals and optimal investment with transaction costs},
  author={Nagai, Hideo},
  journal={Applied Mathematics and Optimization},
  volume={55},
  number={3},
  pages={359--384},
  year={2007},
  publisher={Springer}
}

@article{bauerle2018risk,
  title={Risk-sensitive stopping problems for continuous-time Markov chains},
  author={B{\"a}uerle, Nicole and Popp, Anton},
  journal={Stochastics},
  volume={90},
  number={3},
  pages={411--431},
  year={2018},
  publisher={Taylor \& Francis}
}

@article{jelito2021risk,
  title={Risk sensitive optimal stopping},
  author={Jelito, Damian and Pitera, Marcin and Stettner, {\L}ukasz},
  journal={Stochastic Processes and their Applications},
  volume={136},
  pages={125--144},
  year={2021},
  publisher={Elsevier}
}

@article{palczewski2010finite,
  title={Finite horizon optimal stopping of time-discontinuous functionals with applications to impulse control with delay},
  author={Palczewski, Jan and Stettner, {\L}ukasz},
  journal={SIAM Journal on Control and Optimization},
  volume={48},
  number={8},
  pages={4874--4909},
  year={2010},
  publisher={SIAM}
}

@article{graves1982application,
  title={The application of queueing theory to continuous perishable inventory systems},
  author={Graves, Stephen C},
  journal={Management Science},
  volume={28},
  number={4},
  pages={400--406},
  year={1982},
  publisher={INFORMS}
}

@article{kaspi1983inventory,
  title={Inventory systems of perishable commodities},
  author={Kaspi, H and Perry, D},
  journal={Advances in Applied Probability},
  volume={15},
  number={3},
  pages={674--685},
  year={1983},
  publisher={Cambridge University Press}
}

\newpage
\appendix
\section*{Appendix}

Section~\ref{app:discounted-rs} presents two formulations of the discounted risk-sensitive problem, for which the exact decomposition no longer applies, together with numerical validation. Sections~\ref{app:RN}--\ref{app:RA} collect the supporting proofs for the risk-neutral, risk-seeking, and risk-averse analyses. Section~\ref{app:regime-comparison} proves Proposition~\ref{prop:regime-comparison}, and Section~\ref{app:renewal} provides the proofs supporting the renewal formulation of Section~\ref{sec:renewal}. Section~\ref{app:algorithms} states the exact computational algorithms induced by the threshold characterizations, while Section~\ref{app:sensitivity-theory} develops the comparative statics of the optimal risk-neutral threshold, including the reward--discount phase transition underlying Section~\ref{sec:numerical}.

\section{Discounted Risk-Sensitive Formulations}\label{app:discounted-rs}

This section develops the two discounted extensions of the entropic criterion announced in Remark~\ref{rem:why-delta-1}. Throughout, let $\delta < 1$. For an admissible stopping time $\tau$, define the discounted reward by
\[
\mathcal{R}_{\delta}(\tau) \;=\; \delta^{\tau}(\beta X_\tau - \alpha) - \sum_{t=0}^{\tau-1} \delta^t \bigl( c X_t + \delta \gamma A_t(X_t) \bigr),
\]
so that $J_{rn}(x,\tau) = \mathbb{E}_x[\mathcal{R}_{\delta}(\tau)]$ and $J_{r}(x,\tau) = \eta^{-1} \log \mathbb{E}_x[e^{\eta \mathcal{R}_{\delta}(\tau)}]$.

\begin{lemma}[Deterministic bounds under discounting]\label{lem:disc-bound}
For every $x \in \mathcal{S}$ and every admissible $\tau$, with $\bar t := \max_{t \ge 0} t \delta^t \le (e |\log \delta|)^{-1}$,
\begin{equation}
\begin{aligned}
-\alpha - (c + \delta\gamma)\Bigl( \frac{x}{1-\delta} + \frac{\delta}{(1-\delta)^2} \Bigr) \;\le\; \mathcal{R}_{\delta}(\tau)
\;\le\; \beta\,(x + \bar t) \qquad \mathbb{P}_x\text{-a.s.}
\end{aligned}
\end{equation}
Consequently, $\mathbb{E}_x[e^{\eta \mathcal{R}_{\delta}(\tau)}]$ is finite for every $\eta \in \mathbb{R}$, $x \in \mathcal{S}$, and every admissible $\tau$: under discounting, no exponential-moment restriction is required for either sign of $\eta$. The undiscounted case differs not because the value diverges, but because no discount modulus is available and the unregularized continuation branch lacks uniform control.
\end{lemma}

\begin{proof}[Proof of Lemma~\ref{lem:disc-bound}]
Since $X_{t+1} - X_t \le D_t \le 1$ almost surely, $X_t \le x + t$ and $A_t(X_t) \le X_t$ for all $t$. The lower bound follows by summing $\sum_{t \ge 0} \delta^t (x + t) = x/(1-\delta) + \delta/(1-\delta)^2$; the upper bound follows from $\delta^\tau(\beta X_\tau - \alpha) \le \beta\, \delta^\tau (x + \tau) \le \beta(x + \bar t)$.
\end{proof}

\subsection{Method A: Risk-Parameter Augmentation}

For $\lambda \neq 0$, define
\[
\mathcal{W}(x,\lambda) \;=\; \operatorname*{opt}_{\tau \in \mathfrak{T}}\, \mathbb{E}_x\!\bigl[e^{\lambda \mathcal{R}_{\delta}(\tau)}\bigr],
\qquad
\operatorname*{opt} = \begin{cases} \sup, & \lambda > 0, \\ \inf, & \lambda < 0. \end{cases}
\]

\begin{proposition}[Augmented Bellman recursion]\label{prop:cs-recursion}
For every $x \in \mathcal{S}$ and $\lambda \neq 0$,
\begin{equation}\label{eq:cs-bellman}
\begin{aligned}
\mathcal{W}(x,\lambda)
\;=\; \operatorname*{opt}\Bigl\{\,
 e^{\lambda(\beta x-\alpha)},  \mathbb{E}\!\bigl[ e^{-\lambda(cx + \delta\gamma A_x)}\, \mathcal{W}(\hat{X}_x,\, \lambda\delta) \bigr] \Bigr\},
\end{aligned}
\end{equation}
and the discounted entropic value is
\begin{equation*}
\begin{aligned}
J^r_\delta(x)
&:= \sup_{\tau\in\mathfrak T}\, \eta^{-1}\log \mathbb{E}_x\bigl[e^{\eta \mathcal{R}_{\delta}(\tau)}\bigr]= \eta^{-1}\log \mathcal{W}(x,\eta).
\end{aligned}
\end{equation*}
\end{proposition}

\begin{proof}[Proof of Proposition~\ref{prop:cs-recursion}]
Condition on the first period. If $\tau = 0$, the reward is $\beta x - \alpha$; if $\tau \ge 1$, write $\mathcal{R}_{\delta}(\tau) = -(cX_0 + \delta\gamma A_0) + \delta \mathcal{R}_{\delta}'$, where, conditionally on $(A_0, D_0)$, $\mathcal{R}_{\delta}'$ is distributed as $\mathcal{R}_{\delta}(\tau - 1)$ under the initial state $X_1 = \hat{X}_x$. Factoring the deterministic term out of the exponential and taking $\operatorname*{opt}$ over the continuation and then over the first-period decision yields \eqref{eq:cs-bellman}; the interchange is the optimality principle for exponential-utility control \citep{chung1987discounted}, whose moment assumptions hold by Lemma~\ref{lem:disc-bound}. For $\lambda < 0$, the branches combine by $\inf$ because $w \mapsto \lambda^{-1}\log w$ is decreasing; equivalently, the certainty equivalent is always maximized over $\tau$, since this decreasing map turns minimization of $\mathbb E[e^{\lambda\mathcal R_\delta}]$ into maximization of its certainty equivalent.
\end{proof}

The recursion \eqref{eq:cs-bellman} is the stopping-problem special case of the discounted exponential-utility dynamic program of \citet{chung1987discounted} (cf.\ \citealp{denardo1979optimal,bouakiz1992inventory}). The effective risk sensitivity at period $t$ is $\eta\delta^t$: the decision maker becomes asymptotically risk-neutral, and the optimal rule is stationary in the augmented state $(x, \eta\delta^t)$ but generally not a fixed threshold in $x$. The threshold structure is nonetheless retained: for $\eta > 0$ the optimal rule remains a control limit---now a moving boundary $x^*(\eta\delta^t)$---that converges to the risk-neutral threshold as the effective risk sensitivity vanishes (Proposition~\ref{prop:cs-boundary}).

\begin{proposition}[Boundary form of the optimal policy, $\eta > 0$]\label{prop:cs-boundary}
Let $\eta > 0$ and $V_\lambda(x) := \lambda^{-1} \log \mathcal{W}(x,\lambda)$. Then:
\begin{enumerate}
\item[(i)] $V_\lambda(x+1) - V_\lambda(x) \le \beta$ for all $x \ge 0$ and $\lambda \in (0, \eta]$;
\item[(ii)] the stopping region $\mathcal{S}^*(\lambda) := \{x : \mathcal{W}(x,\lambda) = e^{\lambda(\beta x-\alpha)}\}$ is nonempty and up-closed; writing $x^*(\lambda) := \min \mathcal{S}^*(\lambda)$, the policy that stops at the first $t$ with $X_t \ge x^*(\eta\delta^t)$ is optimal;
\item[(iii)] $x^*(\lambda) \to \theta^*_{rn}$ as $\lambda \downarrow 0$ whenever the optimal discounted risk-neutral threshold $\theta^*_{rn}$ of Theorem~\ref{thm:RN} is unique; in the tie case $\Delta v_{rn}(\theta^*_{rn}) = 0$, $x^*(\lambda) \in \{\theta^*_{rn}, \theta^*_{rn}+1\}$ for all sufficiently small $\lambda$ (both elements are optimal risk-neutral thresholds).
\end{enumerate}
\end{proposition}

\begin{proof}[Proof of Proposition~\ref{prop:cs-boundary}]
(i) Couple two copies started at $x+1$ and $x$ as in Lemma~\ref{lem:rs-increments}: the gap $\zeta_t \in \{0,1\}$ satisfies $\zeta_{t+1} = \zeta_t - B_t$, $B_t \sim \mathrm{Binomial}(\zeta_t, p)$, and for every stopping time $\tau$ applied to both copies,
\[
\mathcal{R}_{\delta}^{x+1} - \mathcal{R}_{\delta}^x
\;=\; \delta^\tau \beta \zeta_\tau - \sum_{t=0}^{\tau-1} \delta^t \bigl( c \zeta_t + \delta \gamma B_t \bigr)
\;\le\; \beta \quad \text{a.s.},
\]
because $\delta^\tau \beta \zeta_\tau \le \beta$ and the cost gap is nonnegative. For $\lambda > 0$ this gives $e^{\lambda \mathcal{R}_{\delta}^{x+1}} \le e^{\lambda\beta} e^{\lambda \mathcal{R}_{\delta}^x}$ pathwise; taking expectations, applying $\lambda^{-1}\log$, and then taking the supremum over $\tau$ yields (i).

(ii) Immediate stopping gives $V_\lambda(x+1) \ge \beta(x+1) - \alpha$; if $x \in \mathcal{S}^*(\lambda)$, then by (i), $V_\lambda(x+1) \le V_\lambda(x) + \beta = \beta(x+1) - \alpha$, hence equality and $x+1 \in \mathcal{S}^*(\lambda)$. For nonemptiness, continuing from $x$ at least one period yields at most $-cx + \delta\beta(x + 1 + \bar t)$ pathwise (the first-period holding cost plus the reward bound of Lemma~\ref{lem:disc-bound} from $x+1$), versus $\beta x - \alpha$ from stopping; the difference $x\bigl(c + \beta(1-\delta)\bigr) - \alpha - \delta\beta(1 + \bar t) \to \infty$, so stopping is strictly optimal at all large $x$. Optimality of the first-entrance rule follows from the one-step look-ahead property of the Snell envelope, as in Proposition~\ref{prop:rs-threshold-wlog}.

(iii) By Lemma~\ref{lem:disc-bound}, $|\mathcal{R}_{\delta}(\tau)| \le B(x,\delta)$ for a deterministic $B(x,\delta)$, uniformly over admissible $\tau$ and for almost every path, so the cumulant expansion
\[
\frac{1}{\lambda}\log\mathbb{E}_x\!\left[e^{\lambda \mathcal{R}_{\delta}(\tau)}\right]
= \mathbb{E}_x[\mathcal{R}_{\delta}(\tau)] + \frac{\lambda}{2}\,\operatorname{Var}_x\!\bigl(\mathcal{R}_{\delta}(\tau)\bigr)  + O(\lambda^2)
\]
holds uniformly over $\tau$, with $\operatorname{Var}_x(\mathcal{R}_{\delta}(\tau)) \le B(x,\delta)^2$; taking the supremum over $\tau$ yields $V_\lambda(x) = V_{rn}(x) + O(\lambda)$ for each fixed $x$, where $V_{rn}$ is the discounted risk-neutral value of Section~\ref{sec:rn}.

For the lower bound, Theorem~\ref{thm:RN}(iii) gives $V_{rn}(x) > S(x)$ at every $x < \theta^*_{rn}$; the expansion is uniform on the finite set $\{0, \dots, \theta^*_{rn} - 1\}$, so $V_\lambda(x) > S(x)$ there for small $\lambda$, and $x^*(\lambda) \ge \theta^*_{rn}$.

For the upper bound, let $V^{\mathrm{cont}}_{rn}(x) := -\mathcal{C} x + \delta\,\mathbb{E}[V_{rn}(\hat{X}_x)]$ and $\Lambda(x):=S(x)-V^{\mathrm{cont}}_{rn}(x)$; the dynamic programming principle gives $\mathbb{E}_x[\mathcal{R}_{\delta}(\tau)] \le V^{\mathrm{cont}}_{rn}(x)$ for every $\tau \ge 1$. By Lemma~\ref{lem:verification}, $\Lambda$ is nonnegative and nondecreasing for $x\ge\theta^*_{rn}$. In the unique-maximizer case, $\Delta v_{rn}(\theta^*_{rn})<0$, so the increment formula for $\Lambda$ gives $\Lambda(\theta^*_{rn})>0$, and monotonicity gives $\Lambda(x)>0$ for every $x\ge\theta^*_{rn}$. In the tie case, $\Lambda(\theta^*_{rn})=0$, but $h(\theta^*_{rn}+2)>0$ and the upward transition from $\theta^*_{rn}+1$ has positive probability; hence the expectation formula for $\Lambda$ gives $\Lambda(\theta^*_{rn}+1)>0$, and again $\Lambda(x)>0$ for every $x\ge\theta^*_{rn}+1$. With $\varepsilon_x := \Lambda(x) > 0$ on the relevant range, every $\tau \ge 1$ satisfies
\[
\frac{1}{\lambda}\log\mathbb{E}_x\!\left[e^{\lambda \mathcal{R}_{\delta}(\tau)}\right]
\;\le\; S(x) - \varepsilon_x + \frac{\lambda}{2}\,B(x,\delta)^2 + O(\lambda^2)
\;<\; S(x)
\]
for all sufficiently small $\lambda$, so stopping is strictly optimal at $x$. Hence $x^*(\lambda) \le \theta^*_{rn}$ in the unique-maximizer case and $x^*(\lambda) \le \theta^*_{rn} + 1$ in the tie case; together with the lower bound, this proves the claim.
\end{proof}

\begin{remark}[Computation]\label{rem:cs-computation}
Equation~\eqref{eq:cs-bellman} is solved by backward recursion on $\lambda_k = \eta\delta^k$, using the terminal approximation $\mathcal{W}(x,\lambda_K) = e^{\lambda_K V_{rn}(x)} + O(\lambda_K^2)$ from Proposition~\ref{prop:cs-boundary}(iii); each step requires one finite expectation over $(A_x, D)$ (validated in Section~\ref{app:disc-numerics}).
\end{remark}

\subsection{Method B: Stationary Bellman}

If instead the entropic measure is applied period by period, the Bellman equation is \eqref{eq:bellman-risk-primal} with $\delta$ inside the exponent:
\begin{equation}\label{eq:bellman-rs-disc}
V^{\delta}_{rs}(x) = \max\Bigl\{ \beta x - \alpha,\; Q^{\delta}_{rs}(x) \Bigr\}, \quad \text{where} 
\quad Q^{\delta}_{rs}(x) := -cx + \frac{1}{\eta} \log \mathbb{E}\!\Bigl[ e^{\eta\delta \bigl( V^{\delta}_{rs}(\hat{X}_x) - \gamma A_x \bigr)} \Bigr].
\end{equation}

\begin{proposition}\label{prop:routeB}
Let $\eta \neq 0$ and let $\mathcal{T}^{\delta}_{rs}$ denote the operator on the right-hand side of \eqref{eq:bellman-rs-disc}.
\begin{enumerate}
\item[(i)] On the truncated state space $\{0, \dots, X_{\max}\}$ with stopping imposed at $X_{\max}$, $\mathcal{T}^{\delta}_{rs}$ is a supremum-norm contraction with modulus $\delta$. On the untruncated space, for every $\rho\in(1,1/\delta)$, it is a contraction on $\mathcal B_\rho$ with modulus $\delta\rho$; hence \eqref{eq:bellman-rs-disc} has a unique fixed point in $\mathcal B_\rho$ and value iteration converges geometrically. Moreover, the fixed point satisfies $S(x)\le V^\delta_{rs}(x)\le\beta x+K_\delta$, where $K_\delta:=\delta\beta/(1-\delta)$.
\item[(ii)] At $\delta = 1$, \eqref{eq:bellman-rs-disc} coincides with \eqref{eq:bellman-risk-primal}: Method~B deforms continuously into the undiscounted model.
\item[(iii)] For $\eta > 0$, the fixed point satisfies $V^{\delta}_{rs}(x+1) - V^{\delta}_{rs}(x) \le \beta$, the stopping region is up-closed, and a stationary threshold policy is optimal.
\item[(iv)] For $\eta<0$, discounting also removes the need for the robust cap: part (i) gives a finite, linearly bounded fixed point directly for the recursive criterion, and the factor $\delta$ inside the exponent supplies the strict modulus. At $\delta=1$ the value remains finite, but the strict modulus is absent and the uncapped continuation branch lacks uniform lower control, motivating the regularized formulation of Section~\ref{sec:rs-averse}.
\end{enumerate}
\end{proposition}

\begin{proof}[Proof of Proposition~\ref{prop:routeB}]
(i) For any random variables $U,W$,
\begin{equation}\label{eq:entropic-lipschitz}
\left|
\frac1\eta\log\mathbb E[e^{\eta\delta U}]
-\frac1\eta\log\mathbb E[e^{\eta\delta W}]
\right|
\le\delta\|U-W\|_\infty .
\end{equation}
Indeed, this follows by comparing $U$ and $W$ after adding $\|U-W\|_\infty$. On the finite grid, where the stopping branch is fixed, this gives modulus $\delta$. On $\mathcal B_\rho$, the one-step
bound
\[
|u(\hat X_x)-w(\hat X_x)|
\le\|u-w\|_\rho\rho^{x+1}
\]
gives modulus $\delta\rho<1$. The finite one-step support also shows that the operator maps $\mathcal B_\rho$ into itself, so Banach's fixed-point theorem applies.

For the linear bound, start value iteration at $S$. If $S\le u$ and $u(x)\le\beta x+K_\delta$, then the outer maximum preserves the lower bound, while the certainty-equivalent upper bound gives
\[
Q^\delta_{rs}(x)
\le\delta[\beta(x+1)+K_\delta]-cx
\le\beta x+K_\delta,
\]
because $(1-\delta)K_\delta=\delta\beta$. These bounds pass to the fixed point. Part~(ii) is immediate.

(iii) Let $\eta>0$. Suppose $u(x+1)-u(x)\le\beta$. Using the coupling
$A_{x+1}=A_x+B$ from Lemma~\ref{lem:increments}, independence of $B$
and $(A_x,D)$ gives
\begin{equation}\label{eq:routeB-Q-increment}
\begin{aligned}
\Delta Q^\delta_{rs}(x)
&\le
-c+\frac1\eta\log\!\left[
(1-p)e^{\eta\delta\beta}+pe^{-\eta\delta\gamma}
\right]\\
&=
\delta\beta-c+\frac1\eta
\log\!\left(1-p+pe^{-\eta\delta(\beta+\gamma)}\right)\\
&=\delta\beta-c-\widetilde\gamma_{\eta\delta}
<\beta,
\end{aligned}
\end{equation}
where
\[
\widetilde\gamma_{\eta\delta} :=-\frac1\eta\log\!\left( 1-p+pe^{-\eta\delta(\beta+\gamma)} \right)>0.
\]
As in Lemma~\ref{lem:increments}, the outer maximum preserves
$\Delta u\le\beta$, and the bound passes to $V^\delta_{rs}$.

To show that stopping occurs, suppose continuation were optimal
everywhere and set $M(x):=\max_{y\le x}V^\delta_{rs}(y)$. For large
$x$,
\[
V^\delta_{rs}(x)\le-cx+M(x+1).
\]
Whenever $V^\delta_{rs}(x+1)>M(x)$, the increment bound therefore
gives $cx\le\beta$. Hence $M$ is eventually constant, contradicting
$V^\delta_{rs}(x)\ge\beta x-\alpha\to\infty$.

If stopping is optimal at $x$, then $Q^\delta_{rs}(x)\le S(x)$, and \eqref{eq:routeB-Q-increment} gives
\[
Q^\delta_{rs}(x+1)
<Q^\delta_{rs}(x)+\beta
\le S(x+1).
\]
Thus stopping is strictly optimal at every higher state. The stopping region is consequently $\{x\ge\theta\}$ for some finite $\theta$. 
Policy evaluation for this threshold rule has the same $\delta\rho$-contraction, so the fixed point is attained by the corresponding first-entrance policy. Part~(iv) follows directly from part~(i), whose argument used no likelihood-ratio cap.
\end{proof}

The interpretation changes accordingly: \eqref{eq:bellman-rs-disc} compounds one-period entropic risk evaluations (the recursive entropic criterion of \citealp{di1999risk,bauerle2014more}), and coincides with the entropic criterion of Section~\ref{sec:rs-criterion} if and only if $\delta = 1$.

\begin{remark}\label{rem:what-breaks}
For $\delta < 1$ the forward decomposition of Sections~\ref{sec:decoupling}--\ref{sec:rs-decoupling} is unavailable under either method: in Method~A the accumulation function would need a step-dependent exponent $\eta\delta^t$, and in Method~B the factor $\delta$ inside the exponent destroys the linearity behind $V_{rs}(x) = f_{rs}(x) + \mathrm{const.}$ Both are instead solved on a truncated grid (one backward pass over the $\lambda$-grid; value iteration of a $\delta$-contraction).
\end{remark}

\subsection{Numerical Validation}\label{app:disc-numerics}

All computations use the baseline calibration of Table~\ref{tab:baseline}. (a) At $\delta = 1$, $\eta = 0.1$, Method~B reproduces the values $V_{rs}(0) = 658.930146$ and $\theta^*_{rs} = 116$ of Theorem~\ref{thm:RS} within $2 \times 10^{-7}$. (b) Method~A agrees with an exact first-passage calculation (policy $\pi_1$: $-88.864963$ both ways) and with $10^6$ Monte Carlo replications ($155.801 \pm 0.083$ vs.\ $155.886$ at $\eta = 0.1$, $\delta = 0.99$). (c) At $\eta = 0.1$, $J^r_\delta(0) = 23.40, 155.89, 497.78, 637.32$ for $\delta = 0.95, 0.99, 0.999, 0.9999$, converging to the undiscounted limit $658.9301$. (d) At $\delta = 0.99$, both methods give thresholds decreasing in $\eta$ (Method~A: $x^* = 81, 65, 55, 49$ at $\eta = 0.5, 0.2, 0.1, 0.05$ and $37, 32, 27$ at $\eta = -0.05, -0.1, -0.2$; Method~B: identical for $\eta > 0$, one to two units smaller for $\eta < 0$), with Method~A values rising from $105.94$ at $\eta = -0.2$ to $211.61$ at $\eta = 0.5$. Both methods reduce as $\eta \to 0$ to the risk-neutral benchmark ($V(0) = 137.9232$, $\theta^*_{rn} = 43$), dominating it for $\eta > 0$ and dominated by it for $\eta < 0$, and the Method~A boundary $x^*(\eta\delta^t)$ relaxes toward $\theta^*_{rn}$ over time (e.g., $\eta = 0.2$: $x^* = 65, 58, 46$ at $t = 0, 50, 200$).

\section{Supporting Proofs Under the Risk-Neutral Criterion}\label{app:RN}

\subsection{Proofs for Section~\ref{sec:rn}}\label{app:RN-main}

We first record the tower-property computation underlying the cost reduction of Section~\ref{sec:rn-bellman}: since $\mathbb{E}[A_t(X_t) \mid \mathcal{F}_t] = pX_t$ and $\{t < \tau\} \in \mathcal{F}_t$,
\[
 \mathbb{E}_x\Bigl[\sum_{t<\tau} \delta^t \bigl(cX_t + \delta\gamma A_t(X_t)\bigr)\Bigr]= \mathbb{E}_x\Bigl[\sum_{t<\tau} \delta^t \bigl(cX_t + \delta\gamma\, \mathbb{E}[A_t(X_t) \mid \mathcal{F}_t]\bigr)\Bigr] = \mathcal{C}\, \mathbb{E}_x\Bigl[\sum_{t<\tau} \delta^t X_t\Bigr].
\]

\begin{lemma}[Pathwise cost and reward bounds]
\label{lem:pathwise-reward-envelope}
Define
\[
\bar U(x) :=\frac{(\beta+c/2)^2}{2c}-\alpha+\frac c2x(x-1),
\qquad x\in\mathbb Z_+ .
\]
For every initial state $x$ and every almost surely finite stopping time $\tau$,
set $M:=\max_{t\le\tau}X_t$. Then
\begin{equation}\label{eq:pathwise-level-sum}
\sum_{t=0}^{\tau-1}X_t
\ge \frac{M(M-1)-x(x-1)}2
\end{equation}
and
\begin{equation}\label{eq:pathwise-reward-envelope}
\mathcal R_1(\tau)
\le K+aM-\frac c2M^2
\le \bar U(x),
\end{equation}
where
\[
K:=-\alpha+\frac c2x(x-1),
\qquad
a:=\beta+\frac c2 .
\]
Moreover, if
$\mathbb E_x[\sum_{t<\tau}X_t]=\infty$,
then $\mathbb E_x[\mathcal R_1(\tau)^-]=\infty$.
\end{lemma}

\begin{proof}[Proof of Lemma~\ref{lem:pathwise-reward-envelope}]
By the skip-free property, every level $\ell\in\{x,\dots,M-1\}$ occurs among $X_0,\dots,X_{\tau-1}$, which gives \eqref{eq:pathwise-level-sum}. Since $X_\tau\le M$ and abandonment penalties are nonnegative,
\[
\mathcal R_1(\tau)
\le \beta M-\alpha-\frac c2
\bigl[M(M-1)-x(x-1)\bigr]
=K+aM-\frac c2M^2 .
\]
Maximizing the concave quadratic over $M\ge0$ gives \eqref{eq:pathwise-reward-envelope}.

To prove the final assertion, we establish the contrapositive. Define the total running cost
\[
C_\tau:=\sum_{t=0}^{\tau-1}
\bigl(cX_t+\gamma A_t(X_t)\bigr),
\qquad\text{so that}\qquad
\mathcal R_1(\tau)=\beta X_\tau-\alpha-C_\tau .
\]
Suppose that $\mathbb E_x[\mathcal R_1(\tau)^-]<\infty$. Because \eqref{eq:pathwise-reward-envelope} also bounds $\mathcal R_1(\tau)$ above by the finite deterministic constant $\bar U(x)$, the positive part $\mathcal R_1(\tau)^+$ is bounded, and hence $\mathcal R_1(\tau)$ is integrable. Moreover, the first inequality in
\eqref{eq:pathwise-reward-envelope} gives
\[
\mathcal R_1(\tau)^-
\ge \left(\frac c2M^2-aM-K\right)^+ .
\]
Integrability of the left-hand side therefore implies $\mathbb E_x[M^2]<\infty$; otherwise, the positive quadratic expression on the right-hand side would have infinite expectation. Consequently, $\mathbb E_x[M]<\infty$, and since $X_\tau\le M$, the terminal reward $\beta X_\tau-\alpha$ is integrable. It now follows from
\[
C_\tau=\beta X_\tau-\alpha-\mathcal R_1(\tau)
\]
that $\mathbb E_x[C_\tau]<\infty$. Finally, $C_\tau\ge c\sum_{t<\tau}X_t$, so
\[
\mathbb E_x\!\left[\sum_{t<\tau}X_t\right]
\le \frac1c\,\mathbb E_x[C_\tau]<\infty .
\]
Thus, if $\mathbb E_x[\sum_{t<\tau}X_t]=\infty$, then $\mathbb E_x[\mathcal R_1(\tau)^-]=\infty$.
\hfill $\Box$
\end{proof}

\begin{proof}[Proof of Proposition~\ref{prop:rn-bellman}]
(i) Since
\[
\mathbb E[\rho^{\hat X_x}] =\rho^x(1-p+p\rho^{-1})^x(1-q+q\rho) \le\rho^x(1-q+q\rho),
\]
the operator $\mathcal T_{rn}$ maps $\mathcal B_\rho$ into itself and satisfies
\[
|\mathcal T_{rn}u(x)-\mathcal T_{rn}w(x)|
\le\delta(1-q+q\rho)\rho^x\|u-w\|_\rho .
\]
The choice of $\rho$ makes the modulus strictly less than $1$. Stopping immediately gives $V_{rn}(x)\ge\beta x-\alpha$, while $X_t\le x+\sum_{s<t}D_s$ implies
\[
J_{rn}(x,\tau)
\le\beta\mathbb E_x\!\left[\sup_{t\ge0}\delta^tX_t\right]
\le\beta x+\frac{\beta q\delta}{(1-\delta)^2}.
\]
Thus $V_{rn}\in\mathcal B_\rho$. The dynamic programming principle, whose integrability condition follows from the same bound, gives \eqref{eq:bellman-equation}; uniqueness and convergence of value iteration follow from contraction.

(ii) The operator is order preserving and
$\mathcal T_{rn}u_0\ge u_0$, so $u_n\uparrow u_\infty$ pointwise. Monotone convergence through the one-step expectation shows that $u_\infty$ is a fixed point. Since $u_n$ is the $n$-stage value,
\[
u_n(x)=\sup_{\tau\in\mathfrak T} \mathbb E_x[\mathcal R_1(\tau\wedge n)]
\le V_{rn}(x),
\]
and hence $u_\infty\le V_{rn}$.

Conversely, fix $\tau$ with $J_{rn}(x,\tau)>-\infty$. Lemma~\ref{lem:pathwise-reward-envelope} implies $\mathbb E_x[\sum_{t<\tau}X_t]<\infty$: otherwise the expected negative part of the reward would be infinite. Fatou's lemma applied to $\beta X_{\tau\wedge n}-\alpha\ge-\alpha$, together with monotone convergence of the cost sum, therefore gives
\[
J_{rn}(x,\tau) \le\liminf_{n\to\infty}
\mathbb E_x[\mathcal R_1(\tau\wedge n)]
\le u_\infty(x).
\]
Taking the supremum over $\tau$ yields $V_{rn}=u_\infty$. Finally, any fixed point $u\ge u_0$ satisfies $u=\mathcal T_{rn}u\ge\mathcal T_{rn}u_n=u_{n+1}$ for every $n$, so $u\ge u_\infty=V_{rn}$, proving minimality.
\end{proof}

\begin{proof}[Proof of Lemma~\ref{lem:increments}]
Set $u_0(x) = \beta x - \alpha$ and $u_{n+1} = \mathcal{T}_{rn} u_n$ for the Bellman operator~\eqref{eq:bellman-operator}. By Proposition~\ref{prop:rn-bellman}, $u_n \to V_{rn}$ pointwise for $\delta < 1$ and $u_n \uparrow V_{rn}$ for $\delta = 1$, so it suffices to verify that (i) $u(x) \ge \beta x - \alpha$ and (ii) $u(x+1) - u(x) \le \beta$ are preserved by $\mathcal{T}_{rn}$; both hold for $u_0$ and pass to the pointwise limit. Property (i) is preserved by construction. For (ii), write the continuation branch as $(Q u)(x) = -\mathcal{C}\, x + \delta\, \mathbb{E}[u(\hat{X}_x)]$ and couple $A_{x+1} = A_x + B$, where $B \sim \operatorname{Bernoulli}(p)$ is independent of $(A_x, D)$, so that $\hat{X}_{x+1} = \hat{X}_x + 1 - B$. Then
\[
\begin{aligned}
(Q u)(x+1) - (Q u)(x)
& = -\mathcal{C} + \delta\, \mathbb{E}\bigl[u(\hat{X}_x + 1 - B) - u(\hat{X}_x)\bigr]\\
& = -\mathcal{C} + \delta (1-p)\, \mathbb{E}\bigl[u(\hat{X}_x+1) - u(\hat{X}_x)\bigr]\\
& \leq -\mathcal{C} + \delta (1-p) \beta \;<\; \beta,
\end{aligned}
\]
where the second equality holds because the integrand vanishes on $\{B = 1\}$, and the inequality uses $u(\hat{X}_x+1) - u(\hat{X}_x) \le \beta$. Taking the maximum with the stop branch preserves the bound: the only nontrivial combination is stopping at $x$ and continuing at $x+1$, in which case $(Qu)(x) \le \beta x - \alpha$ and hence
\[
\begin{aligned}
(\mathcal{T}_{rn}u)(x+1) - (\mathcal{T}_{rn}u)(x)
& = (Q u)(x+1) - (\beta x - \alpha)\\
& \le (Qu)(x) + \beta - (\beta x - \alpha)
\le \beta.
\end{aligned}
\qedhere
\]
\end{proof}

\begin{proof}[Proof of Proposition~\ref{prop:threshold-wlog}]
Suppose $V_{rn}(x) = \beta x - \alpha$. Since stopping immediately is admissible, $V_{rn}(x+1) \ge \beta(x+1) - \alpha$. Combining this with Lemma~\ref{lem:increments} gives $V_{rn}(x+1) \le V_{rn}(x) + \beta = \beta(x+1) - \alpha \le V_{rn}(x+1)$, so $x+1 \in \mathcal{S}^*$. Hence, whenever $\mathcal{S}^*$ is nonempty, it has the form $\{\theta^*, \theta^*+1, \dots\}$ and the first-entrance rule into $\mathcal{S}^*$ is precisely $\pi_{\theta^*}$. Nonemptiness of $\mathcal{S}^*$ and optimality of the first-entrance rule are established in the proof of Theorem~\ref{thm:RN}(iii) in Section~\ref{sec:RN}.
\end{proof}

\begin{proof}[Proof of Proposition~\ref{prop:value-decomposition-and-objective}]
Write $V_\theta:=V_{rn,\pi_\theta}$ and $v_\theta:=V_\theta(0)$, and recall that $\hat{X}_x=x-A_x+D$. The case $\theta\le x$ is immediate, since $\pi_\theta$ stops immediately.

(i) By Remark~\ref{rem:f-rearrange}, $f$ and $g$ satisfy the fixed-point equations
\eqref{eq:f-fixed}--\eqref{eq:g-fixed}. Define
\[
W_\theta(x):=g(x)v_\theta+f(x).
\]
Multiplying \eqref{eq:g-fixed} by $v_\theta$ and adding \eqref{eq:f-fixed} gives
\[
W_\theta(x)
=-\mathcal{C}x+\delta\,\mathbb{E}[W_\theta(\hat{X}_x)].
\]
For $x<\theta$, the value function under $\pi_\theta$ satisfies the same continuation equation:
\[
V_\theta(x)
=-\mathcal{C}x+\delta\,\mathbb{E}[V_\theta(\hat{X}_x)].
\]
In either equation, $\hat{X}_x=x+1$ with probability $(1-p)^xq>0$, while all other realizations of $\hat{X}_x$ lie in $\{0,\dots,x\}$. Hence the value at $x+1$ is uniquely determined by the values at $0,\dots,x$. Since
\[
W_\theta(0)=g(0)v_\theta+f(0)=v_\theta=V_\theta(0),
\]
induction gives $W_\theta(x)=V_\theta(x)$ for every $x\le\theta$, proving (i).

(ii) At the stopping state $\theta$,
\[
g(\theta)v_\theta+f(\theta)
=W_\theta(\theta)
=V_\theta(\theta)
=\beta\theta-\alpha.
\]
Since $g(\theta)\ge1>0$ by Lemma~\ref{lem:monotone-increasing}, solving for $v_\theta=v_{rn}(\theta)$ yields \eqref{eq:reduced-objective}.

(iii) Let $\delta=1$. The fixed-point equation for $g$ becomes
\[
g(x)=\mathbb{E}[g(\hat{X}_x)].
\]
Suppose $g(y)=1$ for every $y\le x$. Then
\[
g(x)
=(1-p)^x\bigl[(1-q)+qg(x+1)\bigr]
+\mathbb{P}(A_x\ge1).
\]
Because $\mathbb{P}(A_x\ge1)=1-(1-p)^x$, this reduces to
\[
(1-p)^xq\,[g(x+1)-1]=0,
\]
so $g(x+1)=1$. Since $g(0)=1$, induction gives $g\equiv1$, and \eqref{eq:reduced-objective} reduces to \eqref{eq:undiscounted-reduced-objective}.
\end{proof}

\subsection{Proofs for Section~\ref{sec:RN}}

\begin{proof}[Proof of Lemma~\ref{lem:monotone-increasing}]
For $u:\mathbb{Z}_{+}\to\mathbb{R}$, write
\[
\Delta u(x):=u(x+1)-u(x),
\qquad
e_k(u):=\Delta^k u(0).
\]
By the binomial-difference identity \eqref{eq:binomial-difference},
\[
\begin{aligned}
&\Delta^k\,\mathbb{E}_{Y_x}\!\left[q u(Y_x+1)+(1-q)u(Y_x)\right] = (1-p)^k\,\mathbb{E}_{Y_x}\!\Bigl[
q\Delta^k u(Y_x+1)+(1-q)\Delta^k u(Y_x)
\Bigr].
\end{aligned}
\]
Since $Y_0=0$ almost surely and
$\Delta^k u(1)=e_k(u)+e_{k+1}(u)$, evaluating at $x=0$ gives
\begin{equation}\label{eq:boundary-difference-recursion}
\begin{aligned}
&\Delta^k\,\mathbb{E}_{Y_x}\!\left[q u(Y_x+1)+(1-q)u(Y_x)\right]\bigg|_{x=0} = (1-p)^k\bigl(e_k(u)+q e_{k+1}(u)\bigr).
\end{aligned}
\end{equation}

Consider first $g$. Applying \eqref{eq:boundary-difference-recursion}
to the fixed-point equation~\eqref{eq:g-fixed} yields, for
every $k\ge0$,
\[
e_k(g)=\delta(1-p)^k\bigl(e_k(g)+q e_{k+1}(g)\bigr),
\quad\text{so}\quad
e_{k+1}(g)=r_k e_k(g),
\quad
r_k:=\frac{1-\delta(1-p)^k}{\delta q(1-p)^k}, \quad
e_0(g)=g(0)=1.
\]
Because $r_k\ge0$ for every $k\ge0$, all boundary differences of $g$
are nonnegative.

For $f$, the fixed-point equation~\eqref{eq:f-fixed} at $x=0$ gives
$0=f(0)=\delta\bigl(qf(1)+(1-q)f(0)\bigr)$,
so $f(1)=0$ and hence $e_0(f)=e_1(f)=0$. Taking one forward
difference of \eqref{eq:f-fixed} at $x=0$ gives
\[
e_1(f)
=
-\mathcal{C}
+\delta(1-p)\bigl(e_1(f)+q e_2(f)\bigr),
\qquad\text{so}\qquad
e_2(f)=\frac{\mathcal{C}}{\delta(1-p)q}>0.
\]
For $k\ge2$, the affine term $-\mathcal{C}x$ has vanishing $k$-th
difference, so
$e_k(f)=\delta(1-p)^k\bigl(e_k(f)+q e_{k+1}(f)\bigr)$,
hence $e_{k+1}(f)=r_k e_k(f)$ for $k\ge2$, where $r_k>0$. It follows
that $e_k(f)\ge0$ for every $k$.

Finally, Newton's forward-difference formula gives, for
$u\in\{f,g\}$,
\[
\Delta u(x)=\sum_{k=0}^{x}\binom{x}{k}e_{k+1}(u)\ge0,
\qquad x\in\mathbb{Z}_{+}.
\]
Thus both $f$ and $g$ are nondecreasing on $\mathbb{Z}_{+}$.
\end{proof}

\begin{proof}[Proof of Lemma~\ref{lem:strict-f-convex}]
Following the notation of the proof of Lemma~\ref{lem:monotone-increasing} (Section~\ref{app:RN}), the boundary differences of $f$ satisfy $e_0(f) = e_1(f) = 0$, $e_2(f) = \mathcal{C}/(\delta(1-p)q) > 0$, and $e_k(f) \ge 0$ for all $k$. Newton's forward-difference formula gives
\[
\Delta^2 f(x) = \sum_{k=0}^{x} \binom{x}{k} e_{k+2}(f) \;\ge\; 0, \qquad x \in \mathbb{Z}_{+},
\]
so $f$ is convex;
moreover, the summand indexed by $k = 0$ equals $\binom{x}{0}\, e_2(f) = e_2(f) > 0$ and all remaining summands are nonnegative, so the inequality is strict and $f$ is strictly convex.
\end{proof}

\begin{proof}[Proof of Lemma~\ref{lem:strict-g-log-convex}]
The proof of Lemma~\ref{lem:monotone-increasing} gives $e_0(g) = 1$ and $e_{k+1}(g) = r_k\, e_k(g)$ with $r_k := \bigl(1 - \delta (1-p)^k\bigr)/\bigl(\delta q (1-p)^k\bigr) \ge 0$, so $e_k(g) \ge 0$, and Newton's forward-difference formula gives $\Delta^2 g(x) = \sum_{k=0}^x \binom{x}{k} e_{k+2}(g) \ge 0$, proving convexity. For $\delta = 1$ there is nothing more to prove, since $g \equiv 1$ by Proposition~\ref{prop:value-decomposition-and-objective}(iii). Assume $\delta < 1$ and write $e_k := \Delta^k g(0)$; then
\[
r_k := \frac{1 - \delta (1-p)^k}{\delta q (1-p)^k} \;\ge\; \frac{1-\delta}{\delta q (1-p)^k} \;>\; 0,
\qquad
r_k - r_{k-1} = \frac{1}{\delta q (1-p)^k} - \frac{1}{\delta q (1-p)^{k-1}} = \frac{p}{\delta q (1-p)^k} > 0,
\]
so $e_k > 0$ for all $k$ and $\{e_k\}$ is strictly log-convex. Symmetrizing Newton's expansion $\Delta^j g(x) = \sum_{k=0}^x \binom{x}{k} e_{k+j}$ gives
\[
\begin{aligned}
&g(x)\, \Delta^2 g(x) - \big(\Delta g(x)\big)^2 = \frac{1}{2} \sum_{i=0}^x \sum_{j=0}^x \binom{x}{i}\binom{x}{j}\,
\big[ e_i e_{j+2} + e_j e_{i+2} - 2 e_{i+1} e_{j+1} \big].
\end{aligned}
\]
For $i = j$, the bracket is $2\big(e_i e_{i+2} - e_{i+1}^2\big) > 0$ by strict log-convexity of $\{e_k\}$. For $j > i$, factor out $e_{i+1} e_{j+1} > 0$: since $r_{j+1}/r_i > r_j/r_{i+1}$ (because $r_{j+1} > r_j$ and $r_i < r_{i+1}$), and since $z + z^{-1} \ge 2$ for $z := r_j/r_{i+1} > 0$,
\[
\begin{aligned}
e_i e_{j+2} + e_j e_{i+2} - 2 e_{i+1} e_{j+1}
&= e_{i+1} e_{j+1} \left( \frac{r_{j+1}}{r_i} + \frac{r_{i+1}}{r_j} - 2 \right) \\
&> e_{i+1} e_{j+1} \left( \frac{r_j}{r_{i+1}} + \frac{r_{i+1}}{r_j} - 2 \right) \ge 0.
\end{aligned}
\]
Hence $g(x)\, \Delta^2 g(x) > (\Delta g(x))^2$, which rearranges to $g(x)\, g(x+2) > g(x+1)^2$.
\end{proof}

\begin{proof}[Proof of Lemma~\ref{lem:f-express-g}]
Seek a linear particular solution $f_p(x) = ax + b$ of \eqref{eq:f-fixed}. Since $\mathbb{E}[Y_x] = (1-p)x$, matching coefficients in $f_p(x) = -\mathcal{C}\, x + \delta\big(a(1-p)x + aq + b\big)$ gives $a = -\mathcal{C}/\big(1-\delta(1-p)\big) = -\bar c$ and $b = \delta a q/(1-\delta) = -K$. The difference $h := f - f_p$ satisfies the homogeneous equation~\eqref{eq:g-fixed}, in which $h(x+1)$ appears with coefficient $\delta (1-p)^x q > 0$ and all remaining terms involve states $\le x$, so induction gives $h = h(0)\, g$ (recall $g(0) = 1$). Since $h(0) = f(0) - f_p(0) = K$, one has $h = K g$; recombining gives the claimed representation, and comparison with \eqref{eq:c-infty} shows $K = C_\infty$.
\end{proof}

\begin{proof}[Proof of Lemma~\ref{lem:verification}]
Set
\[
U(x):=g(x)v^*+f(x),
\qquad
h(x):=U(x)-S(x).
\]
The matching condition $S(x)=g(x)v_{rn}(x)+f(x)$ gives
\begin{equation}\label{eq:h-representation}
h(x)=g(x)\bigl(v^*-v_{rn}(x)\bigr)\ge0.
\end{equation}
Because $\theta^*$ is the smallest maximizer of $v_{rn}$, $h(\theta^*)=0$ and $h(x)>0$ for $x<\theta^*$. Moreover, $\Delta^2h=v^*\Delta^2g+\Delta^2f>0$, so $h$ is strictly convex; since $\Delta h(\theta^*-1)=-h(\theta^*-1)<0$, $h$ is strictly decreasing on $\{0,\dots,\theta^*\}$.

Multiplying \eqref{eq:g-fixed} by $v^*$ and adding \eqref{eq:f-fixed} gives
\begin{equation}\label{eq:U-fixed}
U(x)=-\mathcal{C}x+\delta\,\mathbb{E}[U(\hat{X}_x)].
\end{equation}
Also,
\[
\widehat V=S+h\,\mathbf{1}_{\{\cdot\le\theta^*\}},
\]
so $\widehat V\ge S$, with equality exactly when $x\ge\theta^*$.

For $x<\theta^*$, one has $\hat{X}_x\le x+1\le\theta^*$, hence $\widehat V(\hat{X}_x)=U(\hat{X}_x)$. Therefore, by \eqref{eq:U-fixed},
\[
-\mathcal{C}x+\delta\,\mathbb{E}[\widehat V(\hat{X}_x)]
=U(x)=\widehat V(x)>S(x),
\]
so continuation is optimal at $x$.

For $x\ge\theta^*$, define the stopping advantage
\[
\Lambda(x)
:=S(x)+\mathcal{C}x
-\delta\,\mathbb{E}[\widehat V(\hat{X}_x)].
\]
Stopping is optimal at $x$ exactly when $\Lambda(x)\ge0$. Using \eqref{eq:U-fixed},
\[
\Lambda(x)
=\delta\,\mathbb{E}\!\left[
h(\hat{X}_x)\mathbf{1}_{\{\hat{X}_x>\theta^*\}}
\right]-h(x).
\]
At $x=\theta^*$, only the event $\hat{X}_{\theta^*}=\theta^*+1$ contributes, so
\begin{equation}\label{eq:Lambda-at-theta}
\begin{aligned}
\Lambda(\theta^*)
&=\delta q(1-p)^{\theta^*}h(\theta^*+1)=-\delta q(1-p)^{\theta^*}g(\theta^*+1)\Delta v_{rn}(\theta^*)
\ge0.
\end{aligned}
\end{equation}

For monotonicity, write instead
\[
\Lambda(x)=\mu(x)
-\delta\,\mathbb{E}\!\left[
h(\hat{X}_x)\mathbf{1}_{\{\hat{X}_x\le\theta^*\}}
\right],
\qquad
\mu(x)
:=S(x)+\mathcal{C}x-\delta\,\mathbb{E}[S(\hat{X}_x)].
\]
Since $\mathbb{E}[\hat{X}_x]=(1-p)x+q$, $\mu(x)=\bigl[\beta(1-\delta(1-p))+\mathcal{C}\bigr]x-\delta\beta q-(1-\delta)\alpha$ is strictly increasing. Moreover, coupling $A_{x+1}=A_x+B$, with $B\sim\operatorname{Bernoulli}(p)$, gives $\hat{X}_{x+1}=\hat{X}_x+1-B\ge\hat{X}_x$, and because $h$ is decreasing below $\theta^*$, the expected lower-tail term is nonincreasing in $x$. Hence $\Lambda$ is nondecreasing on $x\ge\theta^*$, and \eqref{eq:Lambda-at-theta} implies $\Lambda(x)\ge0$.

Thus continuation is optimal below $\theta^*$, stopping is optimal from $\theta^*$ upward, and $\widehat V$ solves \eqref{eq:bellman-equation}.
\end{proof}

\section{Supporting Proofs Under the Risk-Seeking Criterion}\label{app:RS}

\subsection{Proofs for Section~\ref{sec:risk-sensitive}}
\begin{proof}[Proof of Proposition~\ref{prop:rs-dpp}]
(i) Lemma~\ref{lem:pathwise-reward-envelope} gives
$\mathcal R_1(\tau)\le\bar U(x)$ almost surely for every almost surely finite stopping time $\tau$. Applying $\eta^{-1}\log\mathbb E_x[\,\cdot\,]$ and taking the supremum over $\tau$ gives $V_{rs}(x)\le\bar U(x)$. Immediate stopping gives the lower bound $V_{rs}(x)\ge\beta x-\alpha$.

For the second claim, let $M_t:=\max_{s\le t}X_s$. Since $\bar U$ is nondecreasing on $\mathbb Z_+$ and abandonment penalties are nonnegative,
\[
\begin{aligned}
V_{rs}(X_t)-\sum_{s=0}^{t-1}(cX_s+\gamma A_s)
&\le \bar U(M_t)-c\sum_{s=0}^{t-1}X_s\\
&\le \bar U(M_t)
-\frac c2\bigl[M_t(M_t-1)-x(x-1)\bigr]\\
&=\bar U(x),
\end{aligned}
\]
where the second inequality uses \eqref{eq:pathwise-level-sum} and the last equality follows from the definition of $\bar U$.

(ii) If $\tau\ge1$, undiscounted multiplicative separability gives
\[
e^{\eta\mathcal R_1(\tau)}
=e^{-\eta cx}\,e^{-\eta\gamma A_0}\,
 e^{\eta\mathcal R'_1(\tau-1)},
\]
where $\mathcal R'_1$ is the reward functional of the time-shifted chain from $X_1=\hat X_x$. Part~(i) implies that all exponential expectations appearing below are finite. Conditioning on $(A_0,D_0)$ shows that the optimal value over $\tau\ge1$ does not exceed the continuation branch in \eqref{eq:bellman-risk-primal}. Conversely, concatenating the first period with $\varepsilon$-optimal stopping rules from each state in the finite one-step support gives the reverse inequality up to $\varepsilon$. Letting $\varepsilon\downarrow0$ and combining continuation with immediate stopping yields \eqref{eq:bellman-risk-primal}.

(iii) The operator $\mathcal T_{rs}$ is order preserving, and $u_0\le\mathcal T_{rs}u_0$. Since $u_0\le V_{rs}$, induction gives $u_n\uparrow u_\infty$ pointwise for some $u_\infty\le V_{rs}\le\bar U$. Monotone convergence through the one-step exponential expectation shows that $u_\infty=\mathcal T_{rs}u_\infty$.

Because $u_n$ is the $n$-stage value, $J_r(x,\tau\wedge n)\le u_n(x)$ for every fixed $\tau\in\mathfrak T$. Fatou's lemma applied to $e^{\eta\mathcal R_1(\tau\wedge n)}$ therefore gives
\[
J_r(x,\tau)
\le\liminf_{n\to\infty}J_r(x,\tau\wedge n)
\le u_\infty(x).
\]
Taking the supremum over $\tau$ yields $V_{rs}\le u_\infty$, and hence $u_\infty=V_{rs}$.

Finally, if $\widehat V$ solves \eqref{eq:bellman-risk-primal} and $\widehat V\ge u_0$, then induction gives
\[
\widehat V=\mathcal T_{rs}\widehat V
\ge\mathcal T_{rs}u_n=u_{n+1}
\qquad\text{for every }n.
\]
Letting $n\to\infty$ gives $\widehat V\ge u_\infty=V_{rs}$, proving minimality.
\end{proof}

\begin{proof}[Proof of Proposition~\ref{prop:risk-V0}]
Write $C_{x \to y}$ for the always-continue cost from $x$ until the first passage to $y > x$ (Remark~\ref{rem:rs-ce}). Since the chain is skip-free to the right (Lemma~\ref{lem:skip-free}), $T_x \le T_\theta$ almost surely for $x \le \theta$, and the strong Markov property at $T_x$ gives the independent decomposition $C_{0 \to \theta} = C_{0 \to x} + C_{x \to \theta}$, so
\begin{equation}\label{eq:rs-mass}
\mathbb{E}_0\!\left[ e^{-\eta C_{0 \to \theta}} \right] = \mathbb{E}_0\!\left[ e^{-\eta C_{0 \to x}} \right] \, \mathbb{E}_x\!\left[ e^{-\eta C_{x \to \theta}} \right].
\end{equation}
Under $\pi_\theta$ from $x \le \theta$, $\tau_\theta = T_\theta$ and $X_{T_\theta} = \theta$ almost surely, so $\mathcal{R}_1(\tau_\theta) = \beta\theta - \alpha - C_{x \to \theta}$ and policy evaluation gives
\[
e^{\eta V_{rs,\pi_\theta}(x)} = e^{\eta (\beta \theta - \alpha)}\, \mathbb{E}_x\!\left[ e^{-\eta C_{x \to \theta}} \right], \qquad x \le \theta.
\]
Dividing by the corresponding value at $x = 0$ and using \eqref{eq:rs-mass} with $e^{\eta V_{rs,\pi_\theta}(x)} = H(x)\, e^{\eta V_{rs,\pi_\theta}(0)}$ from Proposition~\ref{prop:rs-forward}(iii) gives $H(x) = 1 / \mathbb{E}_0[ e^{-\eta C_{0 \to x}} ]$, that is,
\begin{equation}\label{eq:fr-ce}
f_{rs}(x) = -\frac{1}{\eta} \log \mathbb{E}_0\!\left[ e^{-\eta C_{0 \to x}} \right], \qquad x \ge 0,
\end{equation}
independent of $\theta$. Since $e^{\eta f_{rs}(x)} = H(x)$, taking $\eta^{-1}\log$ in Proposition~\ref{prop:rs-forward}(iii) gives the decomposition \eqref{eq:risk-V-decomp}, part (i); evaluating it at $x = \theta$, where $V_{rs,\pi_\theta}(\theta) = \beta \theta - \alpha$, gives \eqref{eq:risk-V0-closed-form}.
\end{proof}

\begin{proof}[Proof of Lemma~\ref{lem:rs-increments}]
Couple copies $X^{x+1}$ and $X^x$ using the same arrival sequence $(D_t)$ and, at states $u \ge \ell$, abandonment counts $A^\ell \sim \operatorname{Binomial}(\ell, p)$ and $A^u = A^\ell + B$ with $B \sim \operatorname{Binomial}(u - \ell, p)$ independent of $A^\ell$. The gap $\zeta_t := X^{x+1}_t - X^x_t$ satisfies $\zeta_0 = 1$ and $\zeta_{t+1} = \zeta_t - B_t$, $B_t \sim \operatorname{Binomial}(\zeta_t, p)$, so $\zeta_t \in \{0, 1\}$. For $y \in \{x,x+1\}$ and any stopping time $\tau$ on the common filtration, let $U^y(\tau) := \beta X^y_{\tau} - \alpha - \sum_{t=0}^{\tau-1} \bigl(c X^y_t + \gamma A^y_t\bigr)$. Applying the same $\tau$ to both copies gives
\begin{equation*}
U^{x+1}(\tau) - U^x(\tau)
= \beta \zeta_\tau - \sum_{t=0}^{\tau-1} \bigl(c \zeta_t + \gamma B_t\bigr)
\le \beta
\end{equation*}
almost surely, because $\beta \zeta_\tau\le\beta$ and the running-cost gap is nonnegative. Since $\eta>0$, $e^{\eta U^{x+1}(\tau)} \le e^{\eta\beta}e^{\eta U^x(\tau)}$ pathwise; taking expectations (all quantities are finite by Proposition~\ref{prop:rs-dpp}(i)), applying $\eta^{-1}\log$, and taking the supremum over $\tau$ yields
\[
V_{rs}(x+1)=\sup_\tau J_r(x+1,\tau)
\le \sup_\tau J_r(x,\tau)+\beta=V_{rs}(x)+\beta.
\qedhere
\]
\end{proof}

\begin{proof}[Proof of Proposition~\ref{prop:rs-threshold-wlog}]
Let \[Q_{r}(x) = -c x + \eta^{-1} \log \mathbb{E}\!\left[ e^{\eta (V_{rs}(\hat{X}_x) - \gamma A_x)} \right]\] denote the continuation value.
Suppose, toward a contradiction, that continuation were strictly optimal everywhere, and set $M(x) = \max_{y \le x} V_{rs}(y)$, finite by Proposition~\ref{prop:rs-dpp}(i). For $\eta > 0$ the entropic certainty equivalent does not exceed the essential supremum, so at any continuation state $x$,
\[
\begin{aligned}
V_{rs}(x) = Q_{r}(x)
&\le -c x + \max_{y \le x+1} V_{rs}(y)= -c x + M(x+1),
\end{aligned}
\]
where the equality holds by Proposition~\ref{prop:rs-dpp}(ii).
If $V_{rs}(x+1) > M(x)$, then $M(x+1) = V_{rs}(x+1)$, and the display with Lemma~\ref{lem:rs-increments} gives $c x \le V_{rs}(x+1) - V_{rs}(x) \le \beta$; hence no new maximum is attained beyond $\beta / c$, so $M(x) = M < \infty$ for large $x$, contradicting $V_{rs}(x) \ge \beta x - \alpha \to \infty$. Thus stopping is optimal at some state and $\theta^*_{rs}$ exists.

Suppose that stopping is optimal at $x$, i.e., $\beta x - \alpha \ge Q_{r}(x)$. Under the coupling of Lemma~\ref{lem:rs-increments}, $A_{x+1} = A_x + B$ and $V_{rs}(\hat{X}_{x+1}) \le V_{rs}(\hat{X}_x) + \beta$ pathwise; since $B$ is independent of $(A_x, D)$,
\[
Q_{r}(x+1) \le Q_{r}(x) + \beta - c - \gamma_\eta, \quad\text{where} \quad
\gamma_\eta := -\frac{1}{\eta} \log \bigl( 1 - p + p e^{-\eta \gamma} \bigr) > 0.
\]
Therefore $\beta (x+1) - \alpha \ge Q_{r}(x) + \beta \ge Q_{r}(x+1) + c + \gamma_\eta > Q_{r}(x+1)$: stopping is strictly optimal at $x+1$, and inductively at every $y > x$. The stopping region is thus $\{x \ge \theta^*_{rs}\}$.

It remains to verify that stopping at the first entrance $\tau^{*} := \inf\{t \ge 0 : X_t \ge \theta^*_{rs}\}$ attains the value; $\tau^{*}$ is a.s.\ finite under every $\mathbb{P}_x$ (as recorded after Lemma~\ref{lem:skip-free}). Consider the multiplicative continuation process
\[
Z_t := \exp\!\Bigl(\eta\Bigl(V_{rs}(X_t) - \sum_{s=0}^{t-1}\bigl(cX_s + \gamma A_s(X_s)\bigr)\Bigr)\Bigr).
\]
On $\{t < \tau^{*}\}$ the Bellman equation~\eqref{eq:bellman-risk-primal} holds with equality, so $\mathbb{E}[Z_{t+1} \mid \mathcal{F}_t] = Z_t$ and $\{Z_{t \wedge \tau^{*}}\}$ is a martingale. By Proposition~\ref{prop:rs-dpp}(i), $Z_{t \wedge \tau^{*}} \le e^{\eta \bar U(x)}$ a.s., so optional sampling at $\tau^{*}$ is legitimate; since $V_{rs}(X_{\tau^{*}}) = \beta X_{\tau^{*}} - \alpha$,
\[
e^{\eta V_{rs}(x)} = \mathbb{E}_x\!\left[Z_{\tau^{*}}\right]
= \mathbb{E}_x\!\left[e^{\eta \mathcal{R}_1(\tau^{*})}\right]
= e^{\eta J_r(x, \tau^{*})}.
\]
Thus $V_{rs}(x) = J_r(x, \tau^{*})$ for every $x$: first-entrance stopping is optimal, and by the skip-free property $\tau^{*}$ is the stopping time of $\pi_{\theta^*_{rs}}$. Finally, $V_{rs}(0) = v_{rs}(\theta^*_{rs}) \le \max_{\theta \ge 0} v_{rs}(\theta) \le V_{rs}(0)$, the last inequality because each $v_{rs}(\theta) = J_r(0, \tau_\theta)$ is the value of an admissible policy; equality holds throughout.
\end{proof}

\subsection{Proofs for Section~\ref{sec:RS}}
\begin{proof}[Proof of Proposition~\ref{prop:rs-forward}]
(i) We induct on $x$. At $x = 0$, $A_0 = 0$ a.s., so $F_{rs}(0) = 0$ and the argument of the logarithm equals $e^{\eta f_{rs}(0)}\, \omega_{rs}(0) = q = H(0)\bigl(e^{0} - 1 + q\bigr) > 0$; hence $f_{rs}(1) = 0$ and $H(1) = 1 = H(0)$. Assume $H(0) \le \dots \le H(x)$ for some $x \ge 1$. On $\{A_x \ge 1\}$, $\hat{X}_x \le x$, so the induction hypothesis and $e^{-\eta\gamma A_x} \le 1$ give
\[
F_{rs}(x) \le H(x)\, \mathbb{P}(A_x \ge 1) = H(x)\bigl(1 - (1-p)^x\bigr).
\]
The logarithm's argument in \eqref{eq:fr-recurrence} therefore satisfies
\[
\begin{aligned}
e^{\eta f_{rs}(x)}\,\omega_{rs}(x) - F_{rs}(x)
& \ge H(x)\Bigl(\omega_{rs}(x) - \bigl(1 - (1-p)^x\bigr)\Bigr) \\
&= H(x)\bigl(e^{\eta c x} - 1 + (1-p)^x q\bigr) > 0,
\end{aligned}
\]
since $e^{\eta c x} \ge 1$ and $(1-p)^x q > 0$. Exponentiating \eqref{eq:fr-recurrence} and applying the same bound gives
\[
\begin{aligned}
H(x+1)
&= \frac{e^{\eta f_{rs}(x)}\,\omega_{rs}(x) - F_{rs}(x)}{(1-p)^x q}\ge H(x)\, \frac{e^{\eta c x} - 1 + (1-p)^x q}{(1-p)^x q}\ge H(x),
\end{aligned}
\]
completing the induction. Nonnegativity and monotonicity of $f_{rs} = \eta^{-1} \log H$ follow, and $\Delta f_{rs}(0) = f_{rs}(1) - f_{rs}(0) = 0$.

(ii) Multiplying the exponentiated \eqref{eq:fr-recurrence} by $(1-p)^x q$ and adding $F_{rs}(x) + (1-p)^x(1-q)H(x)$ to both sides gives
\[
\begin{aligned}
e^{\eta c x} H(x)
&= (1-p)^x \bigl[ q H(x+1) + (1-q) H(x) \bigr] + F_{rs}(x) = \mathbb{E}_{A_x, D}\!\left[ e^{-\eta \gamma A_x}\, H(\hat{X}_x) \right],
\end{aligned}
\]
where the second equality splits the expectation on $\{A_x = 0\}$ and $\{A_x \ge 1\}$ and uses \eqref{eq:Fr-operator}. Dividing by $e^{\eta c x}$ yields the recursion; the coefficient of $H(x+1)$ is read off the $\{A_x = 0, D = 1\}$ atom.

(iii) Write $h_n := e^{\eta V_{rs,\pi_\theta}(n)}$; the claim is $h_n = H(n)\, h_0$ for $n \le \theta$, proved by induction. The base case $n = 0$ holds since $H(0) = 1$. Assume the identity at $0, \dots, n$ for some $n < \theta$. Since $n$ lies in the continuation region of $\pi_\theta$, one-step policy evaluation (Proposition~\ref{prop:rs-dpp}(ii) applied to $\tau_\theta \ge 1$) gives
\[
h_n = e^{-\eta c n}\, \mathbb{E}_{A_n, D}\!\left[ e^{-\eta \gamma A_n}\, h_{\hat{X}_n} \right].
\]
Splitting the expectation on $\{A_n = 0\}$ and $\{A_n \ge 1\}$ and applying the induction hypothesis on $\{A_n \ge 1\}$, where $\hat{X}_n \le n$,
\[
\begin{aligned}
\mathbb{E}\!\left[ e^{-\eta \gamma A_n}\, h_{\hat{X}_n} \right]
& = (1-p)^n \bigl[ q\, h_{n+1}  + (1-q)\, H(n)\, h_0 \bigr] + F_{rs}(n)\, h_0.
\end{aligned}
\]
Substituting $h_n = H(n)\, h_0$ and solving for $h_{n+1}$ gives $h_{n+1} = H(n+1)\, h_0$ by the exponentiated \eqref{eq:fr-recurrence}, completing the induction, which remains valid up to $n + 1 = \theta$ as an identity for the boundary value.
\end{proof}

\begin{proof}[Proof of Proposition~\ref{prop:rs-ce}]
Applying \eqref{eq:fr-ce} at $\theta$ and $\theta + 1$ with the strong-Markov factorization $C_{0 \to \theta + 1} = C_{0 \to \theta} + C_{\theta \to \theta + 1}$ gives \eqref{eq:dfr-ce}. From state $0$, no cost accrues before the first arrival, so $C_{0 \to 1} = 0$ a.s.\ and $f_{rs}(1) = 0$; since $C_{0 \to x} \ge 0$, $f_{rs}(x) \ge 0$. Jensen's inequality gives $-\eta^{-1} \log \mathbb{E}[ e^{-\eta C} ] \le \mathbb{E}[ C ]$ for nonnegative, a.s.-finite $C$; applied to $C_{0 \to x}$ and $C_{\theta \to \theta + 1}$, this yields $f_{rs}(x) \le f(x)$ and $\Delta f_{rs}(\theta) \le \Delta f(\theta)$, since at $\delta = 1$ the risk-neutral quantities are $f(x) = \mathbb{E}[ C_{0 \to x} ]$ and $\Delta f(\theta) = \mathbb{E}[ C_{\theta \to \theta + 1} ]$ (Remark~\ref{rem:fg-meaning}). Part (iii) then follows from part (ii) and the first-order characterizations of Theorems~\ref{thm:RS}(ii) and~\ref{thm:RN}(iii).
\end{proof}

\begin{proof}[Proof of Proposition~\ref{prop:fr-strict-convex}]
By \eqref{eq:dfr-ce}, $\Delta f_{rs}(\theta)$ is the certainty equivalent of the excursion cost $C_{\theta \to \theta + 1}$; we compare consecutive excursions pathwise. Couple always-continue chains $L$ from $\theta$ and $U$ from $\theta + 1$ as in Lemma~\ref{lem:rs-increments}: the gap $\zeta_t = U_t - L_t \in \{0,1\}$ satisfies $\zeta_{t+1} = \zeta_t - B_t$, $B_t \sim \mathrm{Binomial}(\zeta_t, p)$. With $T^L = \inf\{t : L_t = \theta + 1\}$ and $T^U = \inf\{t : U_t = \theta + 2\}$, the skip-free property gives $L_t \le \theta$ for $t < T^L$, hence $U_t \le \theta + 1 < \theta + 2$ and $T^U \ge T^L$. The excursion costs therefore satisfy, pathwise,
\[
\begin{aligned}
C^U - C^L
& = \sum_{t=0}^{T^L - 1} \bigl( c \zeta_t + \gamma B_t \bigr) + \sum_{t = T^L}^{T^U - 1} \bigl( c U_t + \gamma A^U_t \bigr) \ge c \zeta_0 = c,
\end{aligned}
\]
because every summand is nonnegative and the first-period gap equals $c \zeta_0 = c$. Since $z \mapsto e^{-\eta z}$ is strictly decreasing, $\mathbb{E}[ e^{-\eta C^U} ] \le e^{-\eta c} \, \mathbb{E}[ e^{-\eta C^L} ]$, and applying $-\eta^{-1} \log$ gives $\Delta f_{rs}(\theta + 1) \ge \Delta f_{rs}(\theta) + c$, i.e., $\Delta^2 f_{rs}(\theta) \ge c$. The remaining claims follow from $\Delta f_{rs}(0) = 0$ by telescoping.
\end{proof}

\section{Supporting Proofs Under the Risk-Averse Criterion}\label{app:RA}
\subsection{Proofs for Section~\ref{sec:rs-averse}}
\begin{proof}[Proof of Lemma~\ref{lem:ra-continuation-diverges}]
Fix $u\in\mathscr L_K$. Since $\eta<0$ and
$u(y)\le\beta y+K$, with $y=\hat X_x=x-A_x+D$,
\[
\begin{aligned}
\eta\{u(\hat X_x)-\gamma A_x\}
&\ge
\eta\beta(x-A_x+D)+\eta K-\eta\gamma A_x\\
&=
\eta K+\eta\beta x
-\eta(\beta+\gamma)A_x+\eta\beta D .
\end{aligned}
\]
Taking expectations and using the independence of
$A_x\sim\mathrm{Binomial}(x,p)$ and
$D\sim\mathrm{Bernoulli}(q)$ gives
\[
\begin{aligned}
\mathbb E\!\left[
e^{\eta\{u(\hat X_x)-\gamma A_x\}}
\right]
&\ge
e^{\eta K}e^{\eta\beta x}\,
\mathbb E\!\left[e^{-\eta(\beta+\gamma)A_x}\right]\,
\mathbb E\!\left[e^{\eta\beta D}\right]\\
&=
e^{\eta K}
\left[
e^{\eta\beta}
\bigl(1-p+p e^{-\eta(\beta+\gamma)}\bigr)
\right]^x
\bigl(1-q+qe^{\eta\beta}\bigr)\\
&=
e^{\eta K}\bigl(1-q+qe^{\eta\beta}\bigr)\,
\Lambda_\eta^x .
\end{aligned}
\]
Because $\Lambda_\eta>1$, the right-hand side tends to $+\infty$.
Since $\eta<0$, applying $\eta^{-1}\log$ reverses the inequality and
yields
\[
\begin{aligned}
\frac1\eta\log
\mathbb E\!\left[
e^{\eta\{u(\hat X_x)-\gamma A_x\}}
\right]
&\le
\frac1\eta
\left[
x\log\Lambda_\eta
+\log(1-q+qe^{\eta\beta})
+\eta K
\right].
\end{aligned}
\]
Adding $-cx$ gives the displayed bound. The right-hand side depends
on $u$ only through the fixed constant $K$, so the estimate is
uniform over $\mathscr L_K$. Finally,
\[
\log\Lambda_\eta
=\eta\beta
+\log\bigl(1-p+p e^{-\eta(\beta+\gamma)}\bigr)
=-\eta\gamma+\log p+o(1)
\qquad(\eta\to-\infty),
\]
so $\log\Lambda_\eta\to+\infty$; hence $\Lambda_\eta>1$ for all
sufficiently large $|\eta|$.
\end{proof}

\begin{proof}[Proof of Proposition~\ref{prop:kappa-limit}]
Fix $\eta < 0$. For (i), if $\kappa < \kappa'$, then $\mathcal{Q}_\kappa(x) \subseteq \mathcal{Q}_{\kappa'}(x)$, so $\mathcal{T}_{ra}^{\kappa'} \le \mathcal{T}_{ra}^{\kappa}$ pointwise; the ordering propagates through value iteration from $V_0 = S$, and Proposition~\ref{prop:ra-wellposed}(i) gives $V_{ra}^{\kappa'} \le V_{ra}^{\kappa}$. Since $S \le V_{ra}^{\kappa'} \le V_{ra}^{\kappa}$, every stopping state of the $\kappa$-problem is also a stopping state of the $\kappa'$-problem, so $\theta^*_{ra}(\kappa') \le \theta^*_{ra}(\kappa)$.

For (ii), recall that $\mathbb{P}_x \in \mathcal{Q}_\kappa(x)$ for every $\kappa \ge 1$; evaluating the adversarial infimum at the nominal measure gives $\mathcal{T}_{ra}^{\kappa} \le \mathcal{T}_{rn}$ pointwise (cf.\ \eqref{eq:ra-le-rn}), hence $V_{ra}^{\kappa} \le V_{rn}$. Conversely, for any $\mathbb{Q}_x \in \mathcal{Q}_\kappa(x)$, writing $m = d\mathbb{Q}_x/d\mathbb{P}_x \in [0, \kappa]$ and $\varphi(a,d) = V_{ra}^{\kappa}(x - a + d) - \gamma a$,
\[
\begin{aligned}
\mathbb{E}_{\mathbb{Q}_x}[\varphi] - \tfrac{1}{\eta} D_{\mathrm{KL}}(\mathbb{Q}_x \,\|\, \mathbb{P}_x)
& \ge \mathbb{E}_{\mathbb{Q}_x}[\varphi]
= \mathbb{E}_{\mathbb{P}_x}[\varphi] + \mathbb{E}_{\mathbb{P}_x}\bigl[ \varphi\, (m - 1) \bigr] \\
&\ge \mathbb{E}_{\mathbb{P}_x}[\varphi] - 2(\kappa - 1)\, \|\varphi\|_\infty,
\end{aligned}
\]
since $-1/\eta > 0$ and $D_{\mathrm{KL}} \ge 0$, and since $\mathbb{E}_{\mathbb{P}_x}|m-1| = 2\,\mathbb{E}_{\mathbb{P}_x}[(m-1)^+] \le 2(\kappa - 1)$. By (i), $V_{ra}^{\kappa} \uparrow V^* \le V_{rn}$ as $\kappa \downarrow 1$. Passing to the limit in $V_{ra}^{\kappa}(x) = (\mathcal{T}_{ra}^{\kappa} V_{ra}^{\kappa})(x)$ at fixed $x$---the convergence is monotone on the finite one-step support $\{0, \dots, x+1\}$, and the discrepancy displayed above vanishes since $\|\varphi\|_\infty \le \max_{y \le x+1} |V_{rn}(y)| + \alpha + \gamma x$ uniformly in $\kappa$---shows that $V^*$ is a fixed point of $\mathcal{T}_{rn}$ dominating $S$, so $V^* = V_{rn}$ by minimality (Proposition~\ref{prop:rn-bellman}(ii)). For the threshold claim, $V_{rn}(x) > S(x)$ at every $x < \theta^*_{rn}$, and pointwise convergence on the finite set $\{0, \dots, \theta^*_{rn} - 1\}$ gives $V_{ra}^{\kappa}(x) > S(x)$ there for $\kappa$ close to $1$, hence $\theta^*_{ra}(\kappa) \ge \theta^*_{rn}$. Conversely, $V_{ra}^{\kappa} \le V_{rn}$ and $V_{rn} = S$ on $\{\theta^*_{rn}, \theta^*_{rn} + 1, \dots\}$ (Theorem~\ref{thm:RN}(iii)) give $V_{ra}^{\kappa} \le S$ there, hence equality, so $\theta^*_{ra}(\kappa) \le \theta^*_{rn}$.
\end{proof}

\subsection{Proofs for Section~\ref{sec:RA}}
\begin{proof}[Proof of Lemma~\ref{lem:ra-drift}]
Let $m = d\mathbb{Q}_x/d\mathbb{P}_x$ denote the joint likelihood ratio; since $0 \le m \le \kappa$, one has $|m-1| \le K_\kappa := \max\{1,\kappa-1\}$. Centering $A_x$ at its nominal mean $px$ and applying Cauchy--Schwarz,
\[
\begin{aligned}
\mathbb{E}_{\mathbb{Q}_x}[A_x]
 & = px + \mathbb{E}_{\mathbb{P}_x}\!\left[(m-1)(A_x-px)\right] \\
&\ge px - K_\kappa\,\mathbb{E}_{\mathbb{P}_x}|A_x-px|
\ge px - K_\kappa\sqrt{p(1-p)}\,\sqrt{x},
\end{aligned}
\]
which is the first claim with $C_\kappa = K_\kappa\sqrt{p(1-p)}$. Since $D \le 1$ almost surely,
$\mathbb{E}_{\mathbb{Q}_x}[w(\hat{X}_x)] \le (1-p)x + C_\kappa\sqrt{x} + 1 + L$.
For every $\xi \in (1-p,1)$, the difference $(\xi-(1-p))x - C_\kappa\sqrt{x} - 1 - (1-\xi)L$ is positive for all sufficiently large $x$; choose $M$ so that it is positive for $x > M$, and set $b = M + 1 + L$. Then \eqref{eq:ra-drift} holds for $x > M$ by the preceding display, and for $x \le M$ by $\mathbb{E}_{\mathbb{Q}_x}[w(\hat{X}_x)] \le x + 1 + L \le b$. The $m$-step bound follows by induction: $\mathbb{E}[w(X_m)] \le \xi^m w(x) + b \sum_{j < m} \xi^j \le \xi^m w(x) + b/(1-\xi)$.
\end{proof}

\begin{proof}[Proof of Lemma~\ref{lem:ra-concavity}]
In terms of the likelihood ratios $y_i = (d\mathbb{Q}/d\mathbb{P})(i)$, the feasible set $\{ y \in \mathbb{R}^{\mathcal{X}} : 0 \le y_i \le \kappa,\ \sum_i \mathbb{P}(i)\, y_i = 1 \}$ is compact and convex, and the objective $\sum_i \mathbb{P}(i) \bigl( y_i Z_i - \eta^{-1} y_i \log y_i \bigr)$ is continuous on it (with the convention $0 \log 0 = 0$), so the infimum is attained. Monotonicity and concavity in $Z$ follow from taking pointwise infima of the affine, monotone maps $Z \mapsto \mathbb{E}_{\mathbb{Q}}[Z] - \eta^{-1} D_{\mathrm{KL}}(\mathbb{Q} \,\|\, \mathbb{P})$.
\end{proof}

\begin{proof}[Proof of Lemma~\ref{lem:ra-fullsupport}]
Suppose $m^*(a_0, d_0) = 0$ at some atom with $\mathbb{P}_x(a_0, d_0) =: \pi_0 > 0$, where $m^* = d\mathbb{Q}^*/d\mathbb{P}_x$. For small $\epsilon > 0$, set $\mathbb{Q}_\epsilon = (1-\epsilon) \mathbb{Q}^* + \epsilon \mathbb{P}_x$; its likelihood ratio $m_\epsilon = (1-\epsilon) m^* + \epsilon \in (0, \kappa]$, so $\mathbb{Q}_\epsilon \in \mathcal{Q}_\kappa(x)$. In $\mathcal{G}(\mathbb{Q}_\epsilon) - \mathcal{G}(\mathbb{Q}^*)$, the expectation term contributes $O(\epsilon)$, while in $D_{\mathrm{KL}}(\mathbb{Q}_\epsilon \| \mathbb{P}_x) = \mathbb{E}_{\mathbb{P}_x}[m_\epsilon \log m_\epsilon]$ the atom $(a_0, d_0)$ contributes $\pi_0 \, \epsilon \log \epsilon$, while every other atom contributes $O(\epsilon)$ (by differentiability of $m \log m$ on $(0, \kappa]$, or a further nonpositive term of order $\epsilon \log \epsilon$ if $m^* = 0$ there as well). Hence the objective difference is at most $\pi_0 \, \epsilon \log \epsilon + O(\epsilon) < 0$ for small $\epsilon > 0$, contradicting the optimality of $\mathbb{Q}^*$. Thus $m^* > 0$ at every atom, consistent with \eqref{eq:ra-gibbs}.
\end{proof}

\begin{proof}[Proof of Lemma~\ref{lem:ra-increments}]
We induct on $n$ to show that $V_n = \mathcal{T}_{ra}^n S$ satisfies $\Delta V_n(x) \le \beta$ and that the continuation branch built from $V_n$ satisfies the strict decrement bound; both bounds then pass to the limit $V = \lim_n V_n$. The base case $\Delta V_0 = \Delta S = \beta$ is immediate.

For the inductive step, write $\varphi^x(a,d) = V_n(x - a + d) - \gamma a$ and $\mathcal{R}_x(\varphi) = \inf_{\mathbb{Q} \in \mathcal{Q}_\kappa(x)} \{ \mathbb{E}_{\mathbb{Q}}[\varphi] - \eta^{-1} D_{\mathrm{KL}}(\mathbb{Q} \| \mathbb{P}_x) \}$, so the continuation branch is $Q_{ra,n+1}(x) = -cx + \mathcal{R}_x(\varphi^x)$. Assume $\Delta V_n \le \beta$ and couple $A_{x+1} = A_x + B$ with $B \sim \mathrm{Bernoulli}(p)$ independent of $(A_x, D)$. If $B = 0$,
\begin{equation}\label{eq:inc-B0}
\begin{aligned}
\varphi^{x+1}(a, d)
= V_n(x+1-a+d) - \gamma a
&= \varphi^x(a,d) + \Delta V_n(x - a + d)\\
&\le \varphi^x(a,d) + \beta,
\end{aligned}
\end{equation}
If $B = 1$ instead,
\begin{equation}\label{eq:inc-B1}
\begin{aligned}
\varphi^{x+1}(a+1, d)
&= V_n(x - a + d) - \gamma(a+1) = \varphi^x(a,d) - \gamma.
\end{aligned}
\end{equation}
Let $\mathbb{Q}^*_x \in \mathcal{Q}_\kappa(x)$ attain $\mathcal{R}_x(\varphi^x)$ (Lemma~\ref{lem:ra-concavity}), and let $\mathbb{Q}^\dagger$ be the law of $(A_x + B, D)$ when $(A_x, D) \sim \mathbb{Q}^*_x$ and $B \sim \mathrm{Bernoulli}(p)$ independently. The likelihood ratio of $\mathbb{Q}^\dagger$ relative to $\mathbb{P}_{x+1}$ at each atom is a weighted average of the likelihood ratios under $\mathbb{Q}^*_x$, hence lies in $[0, \kappa]$, so $\mathbb{Q}^\dagger \in \mathcal{Q}_\kappa(x+1)$; by the data-processing inequality applied to $(a, b, d) \mapsto (a + b, d)$,
\[
\begin{aligned}
D_{\mathrm{KL}}(\mathbb{Q}^\dagger \| \mathbb{P}_{x+1})
\le D_{\mathrm{KL}}\bigl(\mathbb{Q}^*_x \otimes \mathrm{Bernoulli}(p) \,\big\|\, \mathbb{P}_x \otimes \mathrm{Bernoulli}(p)\bigr) = D_{\mathrm{KL}}(\mathbb{Q}^*_x \| \mathbb{P}_x).
\end{aligned}
\]
Evaluating $\mathcal{R}_{x+1}(\varphi^{x+1})$ at $\mathbb{Q}^\dagger$ and using \eqref{eq:inc-B0}--\eqref{eq:inc-B1},
\[
\begin{aligned}
\mathcal{R}_{x+1}(\varphi^{x+1})
&\le \mathbb{E}_{\mathbb{Q}^*_x}\!\left[ \varphi^x(A_x, D) \right] + \beta\,\mathbb{P}(B=0) - \gamma\,\mathbb{P}(B=1)- \frac{1}{\eta} D_{\mathrm{KL}}(\mathbb{Q}^*_x \| \mathbb{P}_x) \\
&= \mathcal{R}_x(\varphi^x) + \beta (1-p) - \gamma p.
\end{aligned}
\]
Consequently,
\begin{equation}\label{eq:inc-Q}
\begin{aligned}
\Delta Q_{ra,n+1}(x)
&= -c + \mathcal{R}_{x+1}(\varphi^{x+1}) - \mathcal{R}_x(\varphi^x) \le \beta - c - \gamma p,
\end{aligned}
\end{equation}
which is the strict decrement bound. For the value function: if $V_{n+1}(x) = Q_{ra,n+1}(x) > S(x)$, then
\[
\begin{aligned}
\Delta V_{n+1}(x)
& \le \max\bigl\{ S(x+1),\, Q_{ra,n+1}(x+1) \bigr\} - Q_{ra,n+1}(x) \\
&\le \max\bigl\{ S(x+1) - Q_{ra,n+1}(x),\, \Delta Q_{ra,n+1}(x) \bigr\} < \beta,
\end{aligned}
\]
using $Q_{ra,n+1}(x) > S(x)$, so $S(x+1) - Q_{ra,n+1}(x) < \beta$; and if $V_{n+1}(x) = S(x) \ge Q_{ra,n+1}(x)$, then \eqref{eq:inc-Q} gives $Q_{ra,n+1}(x+1) \le S(x) + \beta - c - \gamma p < S(x+1)$, so $V_{n+1}(x+1) = S(x+1)$ and $\Delta V_{n+1}(x) = \beta$. This completes the induction, and the bounds pass to $V = \lim_n V_n$ and its continuation branch by the pointwise convergence of Proposition~\ref{prop:ra-wellposed}(i).
\end{proof}

\begin{proof}[Proof of Corollary~\ref{cor:ra-threshold-wlog}]
(i) Suppose the stopping region were empty, so $V(x) = Q_{ra}(x)$ at every $x$. Iterating the strict decrement bound from Lemma~\ref{lem:ra-increments} gives $V(x) = Q_{ra}(x) \le Q_{ra}(0) + (\beta - c - \gamma p) x$, while $V(x) \ge S(x) = \beta x - \alpha$. Hence $(c + \gamma p) x \le Q_{ra}(0) + \alpha$ for every $x$, which is impossible because $c + \gamma p > 0$.

(ii) Suppose $V(x) = S(x)$. By Lemma~\ref{lem:ra-increments}, $V(x+1) \le V(x) + \beta = S(x+1)$, while $V(x+1) \ge S(x+1)$ always holds; hence $V(x+1) = S(x+1)$, and induction gives the claim.

(iii) For $x < \theta^*$, continuation is optimal by definition of $\theta^*$, and strictness follows from $V(x) > S(x)$. For $x > \theta^*$, the strict decrement bound applied at $x - 1 \ge \theta^*$ gives
\[
\begin{aligned}
Q_{ra}(x) \le Q_{ra}(x-1) + \beta - c - \gamma p
&\le S(x-1) + \beta - c - \gamma p\\
&= S(x) - (c + \gamma p)
< S(x),
\end{aligned}
\]
using $Q_{ra}(x-1) \le V(x-1) = S(x-1)$.
\end{proof}

\begin{proof}[Proof of Lemma~\ref{lem:ra-uniqueness}]
By Proposition~\ref{prop:ra-wellposed}(ii), $\tilde{D} := \tilde V - V \ge 0$, and by hypothesis with $V \ge S$, $0 \le \tilde{D}(x) \le \beta x + K + \alpha \le C w(x)$ for some $C < \infty$.

For each $x$, let $\widehat{\mathbb{Q}}_x \in \mathcal{Q}_\kappa(x)$ attain the infimum defining the continuation branch of $\mathcal{T}_{ra} V$ at $x$ (Lemma~\ref{lem:ra-concavity}); it is unique by strict convexity of the KL term, and by Lemma~\ref{lem:ra-fullsupport} it assigns positive mass to every $y \le x+1$. We claim
\begin{equation}\label{eq:D-subharmonic}
\tilde{D}(x) \le \mathbb{E}_{\widehat{\mathbb{Q}}_x}\!\left[ \tilde{D}(\hat{X}_x) \right] \qquad \text{for every } x.
\end{equation}
Indeed, if $\tilde V(x) = S(x)$, then $\tilde{D}(x) = S(x) - V(x) \le 0$, hence $\tilde{D}(x) = 0$ and~\eqref{eq:D-subharmonic} is trivial. If $\tilde V(x) > S(x)$, then $\tilde V(x)$ equals its continuation branch; evaluating its infimum at $\widehat{\mathbb{Q}}_x$ gives, with $f_u$ as in \eqref{eq:ra-fu},
\[
\tilde V(x) \le -cx + f_{\tilde V}(\widehat{\mathbb{Q}}_x),
\qquad
V(x) \ge -cx + f_{V}(\widehat{\mathbb{Q}}_x),
\]
the second display being the value of $V$'s own continuation branch at its minimizer. Subtracting, the cost, abandonment, and KL terms cancel and~\eqref{eq:D-subharmonic} follows.

By Lemma~\ref{lem:ra-fullsupport} the chain driven by $\widehat{\mathbb{Q}}$ is irreducible on $\mathbb{Z}_{+}$, and by Lemma~\ref{lem:ra-drift} it satisfies the Foster--Lyapunov drift \eqref{eq:ra-drift}; hence it visits every state infinitely often $\widehat{\mathbb{P}}_x$-a.s., with $\mathbb{E}_x[w(X_n)] \le \xi^n w(x) + b/(1-\xi)$.

Iterating \eqref{eq:D-subharmonic}, $M_n := \tilde{D}(X_n)$ is a nonnegative submartingale with $\sup_n \mathbb{E}[M_n] \le C \sup_n \mathbb{E}_x[w(X_n)] < \infty$; hence $M_n$ converges a.s., and since the chain visits every state infinitely often, $\tilde{D} \equiv d \ge 0$ is constant.

Suppose $d > 0$. The continuation branch shifts by $d$ under a constant shift, so $(\mathcal{T}_{ra}(V+d))(x) = \max\{ S(x), Q_{ra}(x) + d \}$, while $(V + d)(x) = \max\{ S(x), Q_{ra}(x) \} + d$. By Corollary~\ref{cor:ra-threshold-wlog} stopping is strictly optimal for $V$ at any $\bar x > \theta^*$, i.e., $V(\bar x) = S(\bar x) > Q_{ra}(\bar x)$. At $\bar x$,
\[
\begin{aligned}
(\mathcal{T}_{ra}(V+d))(\bar x)
&= \max\{ S(\bar x), Q_{ra}(\bar x) + d \} < S(\bar x) + d = (V+d)(\bar x),
\end{aligned}
\]
contradicting $\tilde V = \mathcal{T}_{ra} \tilde V$. Hence $d = 0$.
\end{proof}

\begin{proof}[Proof of Proposition~\ref{prop:ra-verification}]
Let $\mathbb Q^*_x\in\mathcal Q_\kappa(x)$ attain the infimum
defining $Q_{ra}(x)$; by Lemma~\ref{lem:ra-fullsupport}, it has full
support.

\noindent\textit{(i)}
The Bellman equation gives
$V(x)\ge-cx+f_V(\mathbb Q^*_x)$, with equality whenever
$V(x)>S(x)$; equivalently,
\[
\begin{aligned}
\mathbb E[V(X_{t+1})\mid\mathcal F_t]
\le{}&
V(X_t)+cX_t+\gamma\mathbb E[A_t\mid\mathcal F_t]\\
&-\mathcal K_\eta(
\mathbb Q^*_{X_t}\,\|\,\mathbb P_{X_t}
).
\end{aligned}
\]
Thus, under $\mathbb Q^*$,
\[
V(X_{t\wedge\tau})
-\sum_{s<t\wedge\tau}
\left(
cX_s+\gamma A_s
-\mathcal K_\eta(
\mathbb Q^*_{X_s}\,\|\,\mathbb P_{X_s}
)
\right)
\]
is a supermartingale. Since $V\ge S$, evaluating the process at $N$
and using the tower property gives
\[
V(x)\ge\mathbb E^{\mathbb Q^*}_x[R_N(\tau)]
\qquad\text{for every }N.
\]
Taking $\liminf$ yields the stronger strategy-specific bound
\begin{equation}\label{eq:ra-qstar-upper}
V(x)\ge J(x,\tau,\mathbb Q^*).
\end{equation}
Taking the infimum over adversarial strategies then gives
$V(x)\ge\underline J(x,\tau)$.

\noindent\textit{(ii)}
For the $n$-stage game, the reachable state space is finite and all
one-step payoffs are bounded, so backward induction gives
\begin{equation}\label{eq:ra-finite-horizon}
V_n(x)=(\mathcal T_{ra}^nS)(x)
=\sup_{\tau\le n}\inf_{\mathbb Q}
\mathbb E^{\mathbb Q}_x[R_n(\tau)]
=\sup_{\tau\le n}\underline J(x,\tau),
\end{equation}
with attainment by Lemma~\ref{lem:ra-concavity}. Letting
$n\to\infty$ in \eqref{eq:ra-finite-horizon} and using
Proposition~\ref{prop:ra-wellposed}(i) gives
$V(x)\le\sup_{\tau\in\mathfrak T}\underline J(x,\tau)$; the reverse
inequality is part~(i).

\noindent\textit{(iii)}
Write $\theta^*=\theta^*_{ra}$ and $\tau^*=\tau_{\theta^*}$. Since
$\mathbb Q^*$ has full support on the finite continuation region, the
hitting time $\tau^*$ has a geometric tail. At continuation states
the inequality in part~(i) holds with equality, so the corresponding
process is a martingale up to $\tau^*$. Hence
\[
\begin{aligned}
V(x)
&=\mathbb E^{\mathbb Q^*}_x[R_N(\tau^*)]
+\mathbb E^{\mathbb Q^*}_x\!\left[
(V-S)(X_N)\mathbf 1_{\{\tau^*>N\}}
\right]\\
&\longrightarrow J(x,\tau^*,\mathbb Q^*),
\end{aligned}
\]
because $V-S$ is bounded on the continuation region. Thus
\[
\underline J(x,\pi_{\theta^*})
\le J(x,\tau^*,\mathbb Q^*)=V(x).
\]

For the reverse inequality, define the committed-policy best-response
values
\[
\widehat V_N(x)
:=\inf_{\mathbb Q}\mathbb E^{\mathbb Q}_x[R_N(\tau^*)],
\qquad x\in\{0,\dots,\theta^*\}.
\]
Backward induction gives $\widehat V_N=\Psi^N S$, where
\[
(\Psi u)(x):=
\begin{cases}
-cx+\inf_{\mathbb Q_x\in\mathcal Q_\kappa(x)}f_u(\mathbb Q_x),
&x<\theta^*,\\[2mm]
S(\theta^*),&x=\theta^*.
\end{cases}
\]
Since continuation is optimal for $V$ below $\theta^*$,
$\Psi V=V$ on $\{0,\dots,\theta^*\}$. Hence
\[
\widehat V_N=\Psi^N S\le\Psi^N V=V.
\]
Let $W:=\liminf_N\widehat V_N$. Since
$\|\Psi u-\Psi u'\|_\infty\le\kappa\|u-u'\|_\infty$ by the
likelihood-ratio cap, $\Psi$ is continuous on this finite state
space. It is also order preserving, so passing to the $\liminf$ in
$\widehat V_{N+1}=\Psi\widehat V_N$ gives $\Psi W\le W$.

Set $\Xi:=V-W\ge0$; then $\Xi(\theta^*)=0$. For $x<\theta^*$,
choosing $\widehat{\mathbb Q}_x$ to attain the infimum in
$(\Psi W)(x)$ and subtracting
$W(x)\ge-cx+f_W(\widehat{\mathbb Q}_x)$ from
$V(x)\le-cx+f_V(\widehat{\mathbb Q}_x)$ gives
\[
\Xi(x)\le
\mathbb E_{\widehat{\mathbb Q}_x}
[\Xi(\hat X_x)].
\]
Let $\bar x$ be the largest maximizer of $\Xi$ on
$\{0,\dots,\theta^*\}$. Then
\[
\Xi(\bar x)
\le\mathbb E_{\widehat{\mathbb Q}_{\bar x}}[\Xi]
\le\Xi(\bar x),
\]
so $\Xi$ is constant on the support of
$\widehat{\mathbb Q}_{\bar x}$. Full support forces
$\Xi(\bar x+1)=\Xi(\bar x)$ if $\bar x<\theta^*$, contradicting
maximality. Hence $\bar x=\theta^*$, $\Xi\equiv0$, $W=V$, and
$\widehat V_N\to V$.

Finally, for every strategy $\mathbb Q$,
\[
J(x,\tau^*,\mathbb Q)
\ge\liminf_{N\to\infty}\widehat V_N(x)=V(x),
\]
so $\underline J(x,\pi_{\theta^*})\ge V(x)$.

\noindent\textit{Tie case and uniqueness.}
If $Q_{ra}(\theta^*)=S(\theta^*)$, continuation at $\theta^*$ is also
optimal, and the same verification argument applies to
$\pi_{\theta^*+1}$. If $\theta<\theta^*$, then
\[
\underline J(\theta,\pi_\theta)=S(\theta)<V(\theta),
\]
so $\pi_\theta$ is not optimal. If $\theta>\theta^*+1$, choose $x$
with $\theta^*<x<\theta$ and set
$U(y):=J(y,\tau_\theta,\mathbb Q^*)$. Finite-horizon convergence on
$\{0,\dots,\theta\}$ gives
\[
U(y)=-cy+f_U(\mathbb Q^*_y),
\qquad y<\theta.
\]
Moreover, the strategy-specific bound \eqref{eq:ra-qstar-upper},
applied to $\tau_\theta$, gives $U\le V$. Therefore
\[
\underline J(x,\pi_\theta)
\le U(x)
\le-cx+f_V(\mathbb Q^*_x)
=Q_{ra}(x)
<S(x)=V(x),
\]
and no threshold beyond $\theta^*+1$ is optimal.
\end{proof}

\section{Proof of Proposition~\ref{prop:regime-comparison}}\label{app:regime-comparison}

\begin{proof}[Proof of Proposition~\ref{prop:regime-comparison}]
(i) Proposition~\ref{prop:rs-ce}(iii) gives
$\theta^*_{rs,\eta}\ge\theta^*_{rn}$. For the risk-averse problem, $\mathcal T_{ra,\eta'}\le\mathcal T_{rn}$ by \eqref{eq:ra-le-rn}. 
Value iteration from $S$ therefore gives $S\le V_{ra,\eta'}\le V_{rn}$. Consequently, every risk-neutral stopping state is also risk-averse stopping state, and taking minima gives $\theta^*_{ra,\eta'}\le\theta^*_{rn}$.

\medskip
\noindent\textit{(ii) Risk-seeking limit.}
From any $y<\theta$, the probability of reaching $\theta$ within $\theta$ steps is at least $\bigl(q(1-p)^{\theta-1}\bigr)^\theta$. The usual block argument therefore gives geometric hitting-time tails and finite local exponential moments for every excursion cost. Hence the cumulant expansions used below are valid.

By \eqref{eq:fr-ce} and \eqref{eq:dfr-ce}, H\"older's inequality
shows that $f_{rs,\eta}(\theta)$ and $\Delta f_{rs,\eta}(\theta)$ are
nonincreasing in $\eta>0$, while expansion at $\eta=0$ gives
\[
f_{rs,\eta}(\theta)\uparrow f(\theta),
\qquad
\Delta f_{rs,\eta}(\theta)\uparrow\Delta f(\theta)
\qquad(\eta\downarrow0).
\]
Thus $v_{rs,\eta}(\theta)\downarrow v_{rn}(\theta)$ pointwise. Both
maximizations are effectively over
$\bar\Theta:=\{0,\dots,\lceil\beta/c\rceil\}$: the risk-seeking bound
is in Theorem~\ref{thm:RS}(ii), and the risk-neutral bound follows
from $\Delta f(\theta)\ge c\theta$. Convergence of maxima over this
finite set gives
\begin{equation}\label{eq:rs-origin-limit}
V_{rs,\eta}(0)\downarrow V_{rn}(0).
\end{equation}

Next, $\Delta f(0)=0<\beta$, so $\theta^*_{rn}\ge1$. For every
$\theta\ge1$, the excursion cost is nondegenerate and strict Jensen
gives $\Delta f_{rs,\eta}(\theta)<\Delta f(\theta)$. In the unique
risk-neutral maximizer case, $\Delta f(\theta^*_{rn})>\beta$, so
convergence yields $\theta^*_{rs,\eta}=\theta^*_{rn}$ eventually. In
the tie case, strict Jensen at $\theta^*_{rn}$ and convergence at
$\theta^*_{rn}+1$ instead give
$\theta^*_{rs,\eta}=\theta^*_{rn}+1$ eventually.

Operator-level Lyapunov monotonicity and value iteration from $S$
show that $V_{rs,\eta}(x)$ is nondecreasing in $\eta$. The additive
decomposition and the eventual threshold stabilization now give
pointwise convergence: if $x<\theta^*_{rn}$, then eventually
\[
V_{rs,\eta}(x)
=f_{rs,\eta}(x)+V_{rs,\eta}(0)
\downarrow f(x)+V_{rn}(0)=V_{rn}(x).
\]
At $x=\theta^*_{rn}$, either the same display applies or both values
equal $S(x)$. If $x>\theta^*_{rn}$, then eventually
$x\ge\theta^*_{rs,\eta}$, so
$V_{rs,\eta}(x)=S(x)=V_{rn}(x)$.

\medskip
\noindent\textit{(iii) Risk-averse limit.}
Set $k_\eta:=-1/\eta$. For bounded $Z$, define
\[
\mathcal R_{\eta,x}(Z)
:=\inf_{\mathbb Q_x\in\mathcal Q_\kappa(x)}
\left\{
\mathbb E_{\mathbb Q_x}[Z]
+k_\eta D_{\mathrm{KL}}(\mathbb Q_x\|\mathbb P_x)
\right\}.
\]
Comparing the capped and uncapped infima and evaluating at the
nominal measure gives
\[
-k_\eta\log\mathbb E_{\mathbb P_x}[e^{-Z/k_\eta}]
\le\mathcal R_{\eta,x}(Z)
\le\mathbb E_{\mathbb P_x}[Z].
\]
Both sides converge to $\mathbb E_{\mathbb P_x}[Z]$ as
$\eta\uparrow0$, and monotonicity follows from the increasing KL
multiplier. Hence
$\mathcal R_{\eta,x}(Z)\uparrow\mathbb E_{\mathbb P_x}[Z]$.

The corresponding value functions increase to some
$V^*\le V_{rn}$. On each finite one-step support, the uniform bound
\[
|\mathcal R_{\eta,x}(Z)-\mathcal R_{\eta,x}(Z')|
\le\kappa\|Z-Z'\|_\infty
\]
allows passage to the limit in the Bellman equation. Thus $V^*$ is a
fixed point of $\mathcal T_{rn}$ dominating $S$, and minimality gives
$V^*=V_{rn}$.

Finally, part~(i) gives
$\theta^*_{ra,\eta}\le\theta^*_{rn}$. Since
\[
\min_{x<\theta^*_{rn}}\{V_{rn}(x)-S(x)\}>0,
\]
pointwise convergence on this finite set implies
$V_{ra,\eta}(x)>S(x)$ for every $x<\theta^*_{rn}$ once $\eta$ is
sufficiently close to $0$. Therefore
$\theta^*_{ra,\eta}\ge\theta^*_{rn}$ eventually, and equality holds.
\end{proof}

\section{Proofs for Section~\ref{sec:renewal}}\label{app:renewal}

\begin{proof}[Proof of Lemma~\ref{lem:t-recursion}]
Conditioning on $\{A_x=0\}$ and solving for $t(x+1)$ gives the equivalence of \eqref{eq:t-recurrence} and \eqref{eq:t-fixed}, as in Remark~\ref{rem:f-rearrange}. Fix $\theta\ge1$ and set $m_\theta(x):=\mathbb{E}_x[T_\theta]$ and $u(x):=t(\theta)-t(x)$ for $x\le\theta$. Both solve the same Poisson equation $u(x)=1+\mathbb{E}[u(\hat X_x)]$ on $\{0,\dots,\theta-1\}$ with boundary value $0$ at $\theta$ (for $u$, by \eqref{eq:t-fixed}). Since $\theta$ is reached with positive probability from every $x<\theta$, this finite Poisson equation has a unique solution; hence $m_\theta=u$, and evaluating at $0$ gives $\mathbb{E}_0[T_\theta]=t(\theta)-t(0)=t(\theta)$.
\end{proof}

\begin{proof}[Proof of Proposition~\ref{prop:renewal-objectives}]
(i) Under $\pi_\theta$, the cycles are i.i.d. with expected length $t(\theta)$ by Lemma~\ref{lem:t-recursion}. At $\delta=1$, $g\equiv1$ and Remark~\ref{rem:fg-meaning} gives $\mathbb{E}_0[R_\theta]=N_1(\theta)$. The renewal--reward theorem therefore yields \eqref{eq:avg-objective}.

(ii) The one-cycle discounted value is $v_{rn}(\theta;\delta)=N_\delta(\theta)/g(\theta;\delta)$ by Proposition~\ref{prop:value-decomposition-and-objective}. Independence of the cycle lengths and rewards gives a geometric series with ratio $\mathbb{E}_0[\delta^{T_\theta}]=1/g(\theta;\delta)$, from which \eqref{eq:disc-renewal-objective} follows.

(iii) By Remark~\ref{rem:fg-meaning},
\[
f(\theta;\delta) = \frac{C_\theta(\delta)}{\mathbb{E}_0[\delta^{T_\theta}]}, \quad
C_\theta(\delta) :=
\mathbb{E}_0\!\left[ \sum_{k<T_\theta}\delta^k\bigl(cX_k+\delta\gamma A_k(X_k)\bigr)
\right].
\]
Monotone convergence gives $C_\theta(\delta)\uparrow C_\theta(1)$ and $\mathbb{E}_0[\delta^{T_\theta}]\uparrow1$ as $\delta\uparrow1$; hence $f(\theta;\delta)\to f(\theta;1)$ (no monotonicity of the ratio itself is needed). Similarly,
\[
\frac{g(\theta;\delta)-1}{1-\delta} = \frac{1-\mathbb{E}_0[\delta^{T_\theta}]}{(1-\delta)\,\mathbb{E}_0[\delta^{T_\theta}]} = \frac{\mathbb{E}_0\!\left[\sum_{k<T_\theta}\delta^k\right]}
{\mathbb{E}_0[\delta^{T_\theta}]} \longrightarrow
\mathbb{E}_0[T_\theta] = t(\theta).
\]
Combining these two limits gives \eqref{eq:abelian-consistency}.
\end{proof}

\begin{lemma}\label{lem:renewal-denominators}
\begin{enumerate}
\item[(i)] The time-accumulation function is strictly increasing and strictly convex on $\mathbb{Z}_+$: $\Delta t(x)\ge 1/q$ and $\Delta^2 t(x)\ge p/\bigl(q^2(1-p)\bigr)>0$ for every $x\ge0$.
\item[(ii)] For $\delta\in(0,1)$, the discount-accumulation function is strictly increasing and strictly convex on $\mathbb{Z}_+$: $\Delta g(x;\delta)\ge(1-\delta)/(\delta q)>0$ and $\Delta^2 g(x;\delta)>0$ for every $x\ge0$; in particular, $\Delta g(x;\delta)\to\infty$ as
$x\to\infty$.
\end{enumerate}
\end{lemma}

\begin{proof}[Proof of Lemma~\ref{lem:renewal-denominators}]
(i) Apply the boundary-difference calculus of Lemma~\ref{lem:monotone-increasing} (Section~\ref{app:RN}) to \eqref{eq:t-fixed}, writing $e_k(t):=\Delta^k t(0)$. Since $t(0)=0$, $e_0(t)=0$; evaluating \eqref{eq:t-fixed} at $x=0$ gives $e_1(t)=t(1)=1/q>0$. For $k\ge1$, the constant term has vanishing $k$-th forward difference, and \eqref{eq:boundary-difference-recursion} yields
\[
e_k(t)=(1-p)^k\bigl(e_k(t)+q\,e_{k+1}(t)\bigr),
\qquad\text{hence}\qquad
e_{k+1}(t)=\frac{1-(1-p)^k}{q(1-p)^k}\,e_k(t)>0.
\]
By induction, $e_k(t)>0$ for every $k\ge1$, with $e_2(t)=p/\bigl(q^2(1-p)\bigr)$. Newton's forward-difference formula then gives, for every $x\ge0$,
\[
\Delta t(x)=\sum_{k=0}^{x}\binom{x}{k}e_{k+1}(t)\ge e_1(t)=\frac1q>0,
\qquad
\Delta^2 t(x)=\sum_{k=0}^{x}\binom{x}{k}e_{k+2}(t)\ge e_2(t)>0.
\]

(ii) For $\delta\in(0,1)$, the proof of Lemma~\ref{lem:strict-g-log-convex} (Section~\ref{app:RN}) shows that $e_0(g)=1$ and $e_{k+1}(g)=r_k\,e_k(g)$ with $r_k=\bigl(1-\delta(1-p)^k\bigr)/\bigl(\delta q(1-p)^k\bigr)>0$, so $e_k(g)>0$ for every $k\ge0$. Newton's formula gives $\Delta g(x;\delta)\ge e_1(g)=r_0=(1-\delta)/(\delta q)>0$ and $\Delta^2 g(x;\delta)\ge e_2(g)>0$. Finally, $g(\theta;\delta)=1/\mathbb{E}_0[\delta^{T_\theta}]\ge\delta^{-\theta}$ because $T_\theta\ge\theta$; the nondecreasing increments of a convex sequence dominating a geometrically growing one are unbounded, so $\Delta g(x;\delta)\to\infty$.
\end{proof}

\begin{lemma}[Drift identity]\label{lem:drift-identity}
At $\delta=1$, the cost- and time-accumulation functions are affinely related:
\begin{equation}\label{eq:drift-identity}
f(\theta;1)=\frac{\mathcal{C}_1}{p}\,\bigl(q\,t(\theta)-\theta\bigr), \qquad \theta\ge0.
\end{equation}
Equivalently, in view of Remark~\ref{rem:fg-meaning} and Lemma~\ref{lem:t-recursion}, $\mathbb{E}_0\!\left[\sum_{t<T_\theta}X_t\right] =\bigl(q\,\mathbb{E}_0[T_\theta]-\theta\bigr)/p$.
\end{lemma}

\begin{proof}[Proof of Lemma~\ref{lem:drift-identity}]
Set $\tilde f(x):=\frac{\mathcal{C}_1}{p}\bigl(q\,t(x)-x\bigr)$. Using \eqref{eq:t-fixed} and $\mathbb{E}[\hat X_x]=(1-p)x+q$,
\[
-\mathcal{C}_1 x+\mathbb{E}\!\left[\tilde f(\hat X_x)\right] =-\mathcal{C}_1 x+\frac{\mathcal{C}_1}{p}\Bigl(q\bigl(t(x)+1\bigr)-(1-p)x-q\Bigr)
=\frac{\mathcal{C}_1}{p}\bigl(q\,t(x)-x\bigr)
=\tilde f(x),
\]
so $\tilde f$ satisfies \eqref{eq:f-fixed} at $\delta=1$, with $\tilde f(0)=0=f(0)$. The difference $h:=f-\tilde f$ solves $h(x)=\mathbb{E}[h(\hat X_x)]$ with $h(0)=0$; since $h(x+1)$ appears with coefficient $(1-p)^xq>0$ and all remaining terms involve states $\le x$, recursion gives $h\equiv0$, proving \eqref{eq:drift-identity}. The probabilistic statement is optional sampling of the martingale $X_n-\sum_{t<n}(q-pX_t)$ at $T_\theta$ (via $T_\theta\wedge n$, then monotone and dominated convergence using the geometric tail of $T_\theta$): $\mathbb{E}_0[X_{T_\theta}]=\theta=q\,\mathbb{E}_0[T_\theta]-p\,\mathbb{E}_0[\sum_{t<T_\theta}X_t]$.
\end{proof}

\begin{lemma}[Marginal-rate functions]\label{lem:marginal-rate}
Define, for $\theta\ge1$, the marginal rates
\begin{equation}\label{eq:xi-def}
\xi_1(\theta):=\frac{\Delta N_1(\theta)}{\Delta t(\theta)},
\qquad
\xi_\delta(\theta):=\frac{\Delta N_\delta(\theta)}{\Delta g(\theta;\delta)}
\quad(\delta\in(0,1)).
\end{equation}
Both denominators are strictly positive by Lemma~\ref{lem:renewal-denominators}, and
\begin{equation}\label{eq:xi-closed}
\xi_1(\theta)=\frac{\beta+\mathcal{C}_1/p}{\Delta t(\theta)}-\frac{\mathcal{C}_1 q}{p},
\qquad
\xi_\delta(\theta)=\frac{\beta+\bar c}{\Delta g(\theta;\delta)}-C_\infty .
\end{equation}
In particular, $\xi_1$ and $\xi_\delta$ are strictly decreasing on $\{\theta\ge1\}$, with $\xi_1(\theta)\to-\mathcal{C}_1 q/p<0$ and $\xi_\delta(\theta)\to-C_\infty<0$ as $\theta\to\infty$.
\end{lemma}

\begin{proof}[Proof of Lemma~\ref{lem:marginal-rate}]
Both denominators are strictly positive by Lemma~\ref{lem:renewal-denominators}. For the average criterion, Lemma~\ref{lem:drift-identity} gives
\[
\Delta N_1(\theta)=\beta-\Delta f(\theta;1)
=\beta+\frac{\mathcal{C}_1}{p}-\frac{\mathcal{C}_1 q}{p}\,\Delta t(\theta),
\]
and dividing by $\Delta t(\theta)>0$ yields the first expression in \eqref{eq:xi-closed}; since $\Delta t$ is strictly increasing with $\Delta t(\theta)\to\infty$ (Lemma~\ref{lem:renewal-denominators}(i)), $\xi_1$ is strictly decreasing with limit $-\mathcal{C}_1q/p$. For the discounted criterion, Lemma~\ref{lem:f-express-g} gives
\[
\Delta N_\delta(\theta)=\beta-\Delta f(\theta;\delta)
=\beta+\bar c-K\,\Delta g(\theta;\delta),
\qquad K=C_\infty,
\]
and dividing by $\Delta g(\theta;\delta)>0$ yields the second expression; since $\Delta g(\cdot;\delta)$ is strictly increasing with $\Delta g(\theta;\delta)\to\infty$ (Lemma~\ref{lem:renewal-denominators}(ii)), $\xi_\delta$ is strictly decreasing with limit $-C_\infty$.
\end{proof}

\begin{proof}[Proof of Proposition~\ref{prop:renewal-threshold-wlog}]
For $\delta<1$, the discounted Bellman equation is
\begin{equation}\label{eq:renewal-discounted-bellman}
V_\delta(x)=
\begin{cases}
\delta\mathbb E[V_\delta(D)],&x=0,\\[1mm]
\max\{\beta x-\alpha+V_\delta(0),\,(Q_\delta V_\delta)(x)\},&x\ge1,
\end{cases}
\end{equation}
where $(Q_\delta u)(x):=-\mathcal Cx+\delta\mathbb E[u(\hat X_x)]$.
Let $\mathcal U_\delta$ denote the right-hand side. Choose
\[
1<\rho<1+\frac{1-\delta}{\delta q},
\qquad
1\le a<\frac{1-\delta(1-q)}{\delta q\rho},
\]
and equip functions with the norm
\[
\|u\|_{\rho,a}
:=\max\Bigl\{a|u(0)|,\sup_{x\ge1}\rho^{-x}|u(x)|\Bigr\}.
\]
Direct substitution shows that the reset branch, the continuation
branch, and the transition at $x=0$ have respective moduli at most
\[
\frac1{a\rho},\qquad
\delta(1-q+q\rho),\qquad
\delta(1-q+aq\rho),
\]
all strictly less than $1$. Thus $\mathcal U_\delta$ is a contraction,
and its unique fixed point $V_\delta$ is the optimal discounted value.

Starting value iteration from zero, the coupling in
Lemma~\ref{lem:increments} gives
\[
\Delta(Q_\delta u_n)(x)
\le-\mathcal C+\delta(1-p)\beta<\beta,
\qquad x\ge1,
\]
whenever $\Delta u_n\le\beta$. Since the clearing branch has increment
$\beta$, the maximum preserves this bound, which therefore passes to
$V_\delta$. If clearing is optimal at $x$, then
\[
(Q_\delta V_\delta)(x+1)
<(Q_\delta V_\delta)(x)+\beta
\le\beta(x+1)-\alpha+V_\delta(0),
\]
so the clearing region is up-closed. Profitability gives a threshold
policy with positive value, whereas always continuing has negative
value. Hence the clearing region is nonempty and equals
$\{x\ge\theta\}$ for some $\theta\ge1$; the corresponding
first-entrance policy is optimal.

For the average-reward problem, define
\[
\rho_\delta:=(1-\delta)V_\delta(0),
\qquad
h_\delta(x):=V_\delta(x)-V_\delta(0).
\]
The always-continue policy bounds $\rho_\delta$ below, while
discounting processing rewards back to arrival periods gives
$\rho_\delta\le\beta q\delta\le\beta q$. Moreover,
\[
\beta x-\alpha\le h_\delta(x)
\le h_\delta(1)+\beta(x-1),\qquad x\ge1,
\]
and \eqref{eq:renewal-discounted-bellman} at $x=0$ gives
$h_\delta(1)=\rho_\delta/(\delta q)$. Thus $h_\delta(x)$ is bounded
uniformly in $\delta$ for each fixed $x$. Along a suitable sequence
$\delta_n\uparrow1$, suppose $\rho_{\delta_n}\to\rho$ and
$h_{\delta_n}(x)\to h(x)$ for every $x$. Then $h(0)=0$, $h$ has linear
growth, and $\Delta h(x)\le\beta$ for $x\ge1$.

Subtracting $V_\delta(0)$ from
\eqref{eq:renewal-discounted-bellman} gives, for $x\ge1$,
\[
h_\delta(x)
=\max\{\beta x-\alpha,\,(Q_\delta h_\delta)(x)-\rho_\delta\}.
\]
Passing to the limit yields
\begin{equation}\label{eq:renewal-average-bias}
h(x)+\rho=
\begin{cases}
(Q_1h)(0),&x=0,\\[1mm]
\max\{\beta x-\alpha+\rho,\,(Q_1h)(x)\},&x\ge1,
\end{cases}
\end{equation}
or, equivalently,
$h(x)=\max\{\beta x-\alpha,(Q_1h)(x)-\rho\}$ for $x\ge1$.
In particular, $\rho=(Q_1h)(0)=qh(1)$.

For every stationary clearing policy $\pi$, equation
\eqref{eq:renewal-average-bias} implies
\[
\rho\ge r^\pi(x)+(P^\pi h)(x)-h(x),
\]
with equality for maximizing actions. Here $P^\pi$ includes the reset
and forced transition from state $0$ after clearing; on a clearing
state, $(P^\pi h)(x)=(Q_1h)(0)=\rho$. Summing this inequality and
dividing by $n$ proves that $\rho$ dominates every stationary average
reward, with equality for a maximizing policy. The boundary term
vanishes because $h$ has linear growth and the Foster--Lyapunov drift
is uniform over stationary clearing policies.

Finally, $\Delta(Q_1h)(x)<\beta$ for $x\ge1$, so the maximizing
clearing region is up-closed by the same increment argument.
Profitability at $\delta=1$ and Abelian consistency give
\[
\rho\ge
\lim_{\delta\uparrow1}
(1-\delta)v^{\mathrm{ren}}_\delta(\theta_0)
=\bar v(\theta_0)>0
\]
for some $\theta_0\ge1$. Clearing therefore cannot be absent from the
maximizing region, since permanent continuation has average reward
$-\mathcal C_1q/p<0$. The maximizing clearing region consequently
equals $\{x\ge\theta\}$ for some $\theta\ge1$, and $\pi_\theta$ is
optimal.
\end{proof}

\begin{proof}[Proof of Theorem~\ref{thm:renewal-unimodal}]
Write $v(\theta):=N(\theta)/D(\theta)$, $\theta\ge1$, where $N:=N_1$ and $D:=t$ for the average-reward criterion, $N:=N_\delta$ and $D:=g(\cdot;\delta)-1$ for the discounted criterion, and $\xi:=\xi_1$ or $\xi:=\xi_\delta$ accordingly. Three properties drive the argument: $D(\theta)>0$ for $\theta\ge1$ (since $t(\theta)=\mathbb{E}_0[T_\theta]\ge\theta$ and $\mathbb{E}_0[\delta^{T_\theta}]<1$); $\Delta D$ is positive and strictly increasing (Lemma~\ref{lem:renewal-denominators}); and $\Delta N(\theta)=\beta-\Delta f(\theta;\delta)$ is strictly decreasing (Lemma~\ref{lem:strict-f-convex}), so $\xi=\Delta N/\Delta D$ is strictly decreasing by Lemma~\ref{lem:marginal-rate}. Rearranging $\Delta v(\theta)=v(\theta+1)-v(\theta)$ gives the sign identity
\begin{equation}\label{eq:renewal-sign-xi}
\operatorname{sgn}\Delta v(\theta)
=\operatorname{sgn}\Bigl(\Delta N(\theta)-v(\theta)\,\Delta D(\theta)\Bigr)
=\operatorname{sgn}\bigl(\xi(\theta)-v(\theta)\bigr),
\qquad \theta\ge1,
\end{equation}
and $v(\theta+1)$ is the mediant of $v(\theta)$ and $\xi(\theta)$:
\begin{equation}\label{eq:renewal-mediant}
v(\theta+1)=\frac{N(\theta)+\Delta N(\theta)}{D(\theta)+\Delta D(\theta)}
\quad\text{satisfies}\quad
\min\{v(\theta),\xi(\theta)\}\le v(\theta+1)\le\max\{v(\theta),\xi(\theta)\},
\end{equation}
strictly so unless $v(\theta)$ and $\xi(\theta)$ coincide.

\medskip
\noindent\textit{Single crossing.} Suppose $\Delta v(\theta)\le0$ for some $\theta\ge1$. By \eqref{eq:renewal-sign-xi}, $v(\theta)\ge\xi(\theta)$; the mediant property \eqref{eq:renewal-mediant} then gives $v(\theta+1)\ge\xi(\theta)>\xi(\theta+1)$, and \eqref{eq:renewal-sign-xi} at $\theta+1$ yields $\Delta v(\theta+1)<0$. Iterating, every subsequent increment is strictly negative: $\Delta v$ changes sign at most once, from positive to negative, and a zero increment can occur only at the crossing index itself.

\medskip
\noindent\textit{A crossing occurs.} By Assumption~\ref{ass:profitability}, $v(\theta_0)>0$ for some $\theta_0\ge1$. If every increment were positive, then $v(\theta)\ge v(\theta_0)>0$ for all $\theta\ge\theta_0$, and \eqref{eq:renewal-sign-xi} would force $\xi(\theta)>v(\theta_0)>0$ for every $\theta$, contradicting the negative limit of $\xi$: $\xi_1(\theta)\to-\mathcal C_1q/p$ in the average-reward case and $\xi_\delta(\theta)\to-C_\infty$ in the discounted case (Lemma~\ref{lem:marginal-rate}). Hence the set in \eqref{eq:renewal-argmax} is nonempty.

\medskip
With $\theta^u:=\min\{\theta\ge1:\Delta v(\theta)\le0\}$, the increments are strictly positive below $\theta^u$, nonpositive at $\theta^u$, and strictly negative above it, proving parts (i) and (ii). For part (iii), $\theta^u$ is the smallest maximizer of $v$; the only other possible maximizer is $\theta^u+1$, in the tie case $\Delta v(\theta^u)=0$. Proposition~\ref{prop:renewal-threshold-wlog} upgrades the resulting threshold policies to optimality among all stationary clearing policies.
\end{proof}

\begin{proof}[Proof of Proposition~\ref{prop:renewal-ordering}]
Because the chain is skip-free to the right, reaching $\theta+1$ requires first reaching $\theta$ and at least one further step, so $T_{\theta+1}\ge T_\theta+1$ almost surely; hence $t$ is strictly increasing and, for $\delta<1$, $g(\cdot;\delta)=1/\mathbb{E}_0[\delta^{T_\theta}]$ is strictly increasing.

For the average criterion, let $\theta^0:=\theta^*_{rn}(1)$ and fix $\theta>\theta^0$. By Assumption~\ref{ass:profitability} at $\delta=1$, $N_1(\theta^0)=\max_{\theta'}N_1(\theta')>0$, so $N_1(\theta)\le N_1(\theta^0)$, while $t(\theta)>t(\theta^0)$. If $N_1(\theta)\le0$, then $N_1(\theta)/t(\theta)\le0<N_1(\theta^0)/t(\theta^0)$ immediately; if $N_1(\theta)>0$, then
\[
\frac{N_1(\theta)}{t(\theta)}
\le
\frac{N_1(\theta^0)}{t(\theta)} < \frac{N_1(\theta^0)}{t(\theta^0)},
\]
using $N_1(\theta^0)>0$ and $t(\theta)>t(\theta^0)$ (the inequality remains strict in the tie case $N_1(\theta)=N_1(\theta^0)$). Hence no threshold above $\theta^0$ maximizes $\bar v$.

For the discounted criterion, let $\theta^\delta:=\theta^*_{rn}(\delta)$ and fix $\theta>\theta^\delta$. By Assumption~\ref{ass:profitability} at $\delta$, $v_{rn}(\theta^\delta;\delta)>0$, and since $\theta^\delta$ is the smallest maximizer of $v_{rn}(\cdot;\delta)=N_\delta/g_\delta$,
\[
\frac{N_\delta(\theta)}{g(\theta;\delta)}
\le
\frac{N_\delta(\theta^\delta)}{g(\theta^\delta;\delta)},
\]
with equality in the tie case. The map $r\mapsto r/(r-1)$ is strictly decreasing on $(1,\infty)$, and $g(\theta;\delta)> g(\theta^\delta;\delta)>1$, so when $N_\delta(\theta)/g(\theta;\delta)\ge0$,
\[
\frac{N_\delta(\theta)}{g(\theta;\delta)-1} = \frac{N_\delta(\theta)}{g(\theta;\delta)}
  \frac{g(\theta;\delta)}{g(\theta;\delta)-1} < \frac{N_\delta(\theta^\delta)}{g(\theta^\delta;\delta)}
  \frac{g(\theta^\delta;\delta)}{g(\theta^\delta;\delta)-1} = \frac{N_\delta(\theta^\delta)}{g(\theta^\delta;\delta)-1},
\]
and when $N_\delta(\theta)/g(\theta;\delta)<0$ the same inequality is immediate because the right-hand side is positive. Hence no threshold above $\theta^\delta$ maximizes the repeated-clearing objective.
\end{proof}

\section{Computational Algorithms}\label{app:algorithms}

\subsection{Forward-Recursive Threshold Algorithms}
\label{app:forward-algorithms}

The risk-neutral, risk-seeking, and repeated-clearing problems admit the same forward-ascent implementation. For each regime $j$, let $\mathcal{A}_j$ denote its accumulation sequences and $v_j$ its reduced objective, as specified in Table~\ref{tab:forward-ascent}. Starting from the indicated initial index $x_j$, the recurrences are extended one state at a time. Theorems~\ref{thm:RN}, \ref{thm:RS}, and \ref{thm:renewal-unimodal} imply that the first nonpositive increment of $v_j$ identifies the smallest optimal threshold; in the equality case,
the adjacent threshold $x+1$ is also optimal.

\begin{algorithm}[htbp]
\caption{Generic Forward Ascent for an Optimal Threshold}
\label{alg:forward-ascent}
\begin{algorithmic}[1]
\Require A regime $j$ and its initialization, forward recurrences, and
reduced objective in Table~\ref{tab:forward-ascent}.
\Ensure Smallest optimal threshold $\theta_j^*$ and optimal value
$v_j(\theta_j^*)$.
\State Initialize $x \gets x_j$, the accumulation sequences
$\mathcal{A}_j$, and $v_j(x_j)$.
\Loop
    \State Compute the transition-weighted sums using the stored values
    of $\mathcal{A}_j$ at states $0,\ldots,x$.
    \State Extend $\mathcal{A}_j$ to state $x+1$ using the indicated
    forward recurrences.
    \State Evaluate $v_j(x+1)$.
    \If{$v_j(x+1) \le v_j(x)$}
        \Return $\theta_j^* \gets x$ and $v_j(x)$.
    \EndIf
    \State $x \gets x+1$.
\EndLoop
\end{algorithmic}
\end{algorithm}

\begin{table}[htbp]
\centering
\caption{Regime-specific inputs to Algorithm~\ref{alg:forward-ascent}.}
\label{tab:forward-ascent}
\small
\begin{tabular}{@{}p{.17\textwidth}p{.29\textwidth}p{.18\textwidth}p{.26\textwidth}@{}}
\toprule
Regime & Initialization & Forward update & Reduced objective \\
\midrule
Risk-neutral
& $x_j=0$; $f(0)=0$, $g(0)=1$
& \eqref{eq:f-recurrence}--\eqref{eq:g-recurrence}
& $\displaystyle v_{rn}(\theta)=
\frac{\beta\theta-\alpha-f(\theta)}{g(\theta)}$
\\[6pt]
Risk-seeking
& $x_j=1$; $f_{rs}(0)=f_{rs}(1)=0$
& \eqref{eq:fr-recurrence}
& $\displaystyle v_{rs}(\theta)=
\beta\theta-\alpha-f_{rs}(\theta)$
\\[6pt]
Repeated clearing, average reward
& $x_j=1$; $(f,t)(0)=(0,0)$; compute $f(1),t(1)$
& \eqref{eq:f-recurrence} at $\delta=1$, and \eqref{eq:t-recurrence}
& $\displaystyle \bar v(\theta)=
\frac{\beta\theta-\alpha-f(\theta)}{t(\theta)}$
\\[6pt]
Repeated clearing, discounted
& $x_j=1$; $(f,g)(0)=(0,1)$; compute $f(1),g(1)$
& \eqref{eq:f-recurrence}--\eqref{eq:g-recurrence}
& $\displaystyle v^{\mathrm{ren}}_{\delta}(\theta)=
\frac{\beta\theta-\alpha-f(\theta)}{g(\theta)-1}$
\\
\bottomrule
\end{tabular}
\end{table}

Each transition-weighted sum at state $x$ requires $\mathcal{O}(x)$
operations. Consequently, every case in Table~\ref{tab:forward-ascent}
uses $\mathcal{O}((\theta_j^*)^2)$ time and $\mathcal{O}(\theta_j^*)$
space. Finite termination follows from the corresponding optimality
theorem; additionally, $\theta_{rs}^*\le\lceil\beta/c\rceil$.
At $\delta=1$, set $g\equiv1$ directly rather than evaluating its
numerically sensitive recursion.

\subsection{Computational Algorithm for the Risk-Averse Problem}
For the regularized risk-averse problem ($\eta < 0$, $\delta = 1$), no forward recursion is available; we solve \eqref{eq:bellman-robust} by value iteration from $V_0 = S$, which increases pointwise to $V_{ra}^{\kappa}$ (Theorem~\ref{thm:RA}(i)). Proposition~\ref{prop:regime-comparison} provides an exact truncation: imposing stopping at $X_{\max} = \theta^*_{rn}$ leaves the solution unchanged on $\{0, \dots, X_{\max}\}$.

At each state $x$, the continuation value requires the worst-case distortion $\mathbb{Q}_x^* \in \mathcal{Q}_\kappa(x)$ of the joint law $\mathbb{P}_x$ over $(a, d) \in \{0, \dots, x\} \times \{0, 1\}$. By Theorem~\ref{thm:RA}(iii), the likelihood ratio is the truncated tilt $m^*(a,d) = \min\bigl(\kappa, \lambda_x e^{\eta \varphi(a,d)}\bigr)$ with $\varphi(a,d) = V_{ra}^{\kappa}(x-a+d) - \gamma a$, so the variational problem reduces to the scalar equation
\[
\sum_{a=0}^{x}\sum_{d=0}^{1} \mathbb{P}_x(a,d)\, \min\bigl(\kappa,\, \lambda_x e^{\eta \varphi(a,d)}\bigr) = 1,
\]
whose left-hand side is continuous in $\lambda_x$, increases from $0$ to $\kappa > 1$, and is strictly increasing at the root; we locate $\log \lambda_x$ by bisection over an adaptively expanded bracket.

\begin{algorithm}[htbp]
\caption{Value Iteration for the Regularized Risk-Averse Bellman Equation}
\label{alg:robust-vi}
\begin{algorithmic}[1]
\Require Parameters $(\alpha,\beta,c,\gamma,p,q)$; risk-sensitivity $\eta < 0$; cap $\kappa > 1$; truncation $X_{\max} = \theta^*_{rn}$; tolerance $\varepsilon > 0$.
\Ensure Fixed point $V_{ra}^{\kappa}$ on $\{0, \dots, X_{\max}\}$; smallest optimal threshold $\theta^*_{ra}$.
\State \textbf{Initialization:} $V(x) \gets \beta x - \alpha$ for $x \in \{0, \dots, X_{\max}\}$; set the imposed boundary $V(X_{\max}+1) \gets \beta(X_{\max}+1) - \alpha$.
\Repeat
    \For{$x = 0,1,\dots,X_{\max}$}
        \State $V_{\text{stop}} \gets \beta x - \alpha$
        \State \textbf{Nominal joint law and continuation payoffs:}
        \For{$a = 0,1,\dots,x$ and $d \in \{0,1\}$}
            \State $\mathbb{P}_x(a,d) \gets \binom{x}{a} p^a (1-p)^{x-a} \cdot q^{d}(1-q)^{1-d}$
            \State $\varphi(a,d) \gets V(x-a+d) - \gamma a$ \Comment{$x-a+d \le X_{\max}+1$ uses the imposed boundary}
        \EndFor
        \State \textbf{Worst-case distortion (bisection on $\log \lambda_x$):} find $\lambda_x > 0$ such that
        \Statex \hspace{2.5em} $\sum_{a=0}^{x}\sum_{d=0}^{1} \mathbb{P}_x(a,d)\, m^*(a,d) = 1$,
        \Statex \hspace{2.5em} where $m^*(a,d) \gets \min\bigl(\kappa,\, \lambda_x e^{\eta \varphi(a,d)}\bigr)$
        \State $\mathbb{Q}_x^*(a,d) \gets \mathbb{P}_x(a,d)\, m^*(a,d)$ for all $(a,d)$
        \State $D_{\mathrm{KL}} \gets \sum_{a,d} \mathbb{Q}_x^*(a,d) \log m^*(a,d)$
        \State $V_{\text{cont}} \gets -c x + \sum_{a,d} \mathbb{Q}_x^*(a,d)\, \varphi(a,d) - \frac{1}{\eta} D_{\mathrm{KL}}$
        \State $V_{\text{new}}(x) \gets \max\{V_{\text{stop}},\; V_{\text{cont}}\}$
    \EndFor
    \State $\Delta \gets \max_x |V_{\text{new}}(x) - V(x)|$
    \State $V \gets V_{\text{new}}$
\Until{$\Delta < \varepsilon$}
\State $\theta^*_{ra} \gets \min\{x : V(x) \le \beta x - \alpha + \sqrt{\varepsilon}\}$
\State \Return $V$, $\theta^*_{ra}$
\end{algorithmic}
\end{algorithm}

\noindent{\it Computational complexity.} Each sweep costs $\mathcal{O}(X_{\max}^2)$ operations ($2(x+1)$ atoms per state; bisection contributes a constant factor) and $\mathcal{O}(X_{\max})$ memory. At the baseline calibration, convergence to $\varepsilon = 10^{-8}$ requires on the order of $10^3$ sweeps, reducible tenfold by warm-starting from a nearby parameter value.

\section{Theoretical Sensitivity Results}\label{app:sensitivity-theory}

This section develops the comparative statics of the optimal risk-neutral threshold $\theta^*_{rn}$ of Theorem~\ref{thm:RN}: Lemma~\ref{lem:sensitivity-costs} covers the cost parameters $(c, \gamma)$ and $\alpha$; the remaining results concern the interaction of $\beta$ with $\delta$, governed by the sign of a single phase function. Throughout, $f$ and $g$ are the accumulation functions of Section~\ref{sec:decoupling}; we write $g(\theta;\delta)$ when the dependence on $\delta$ matters.

\begin{lemma}\label{lem:sensitivity-costs}
Let $\theta^*_{rn}$ be the optimal risk-neutral threshold characterized in Theorem~\ref{thm:RN}, and fix all other primitives.
\begin{enumerate}
    \item[(i)] $\theta^*_{rn}$ is nonincreasing in each of the cost parameters $c$ and $\gamma$.
    \item[(ii)] For $\delta < 1$, $\theta^*_{rn}$ is nondecreasing in $\alpha$. For $\delta = 1$, $v_{rn}(\theta) = \beta\theta - \alpha - f(\theta)$, so $\alpha$ shifts the objective by a constant and $\theta^*_{rn}$ is independent of $\alpha$.
\end{enumerate}
\end{lemma}

\begin{proof}[Proof of Lemma~\ref{lem:sensitivity-costs}]
(i) With $\mathbb E_{\mathbb P_x}[A_x]=px$, write the risk-neutral
Bellman operator with cost parameters $(c,\gamma)$ as
\[
\mathcal T_{c,\gamma}V(x)
=\max\!\left\{
S(x),\,
-cx-\delta\gamma px
+\delta\mathbb E_{\mathbb P_x}[V(\hat X_x)]
\right\}.
\]
Let $(c,\gamma)\le(c',\gamma')$ componentwise. Since $S$ does not
depend on $(c,\gamma)$ and the continuation branch is componentwise
nonincreasing in $(c,\gamma)$, induction on the value iterates gives
\[
\mathcal T_{c',\gamma'}^nS
\le\mathcal T_{c,\gamma}^nS
\qquad\text{for every }n\ge0.
\]
Passing to the limit by Proposition~\ref{prop:rn-bellman} yields
\[
V_{rn}(\cdot\,;c',\gamma')
\le V_{rn}(\cdot\,;c,\gamma)
\qquad\text{componentwise.}
\]
By Theorem~\ref{thm:RN}(iii), the optimal stopping region is
\[
\mathcal S^*(c,\gamma)
=\{\theta^*_{rn}(c,\gamma),
\theta^*_{rn}(c,\gamma)+1,\dots\}.
\]
Since $V_{rn}\ge S$ everywhere, the displayed value inequality implies
\[
\mathcal S^*(c,\gamma)
\subseteq\mathcal S^*(c',\gamma'),
\]
and hence
$\theta^*_{rn}(c',\gamma')\le\theta^*_{rn}(c,\gamma)$.

\medskip
(ii) For $\delta<1$,
Proposition~\ref{prop:value-decomposition-and-objective} gives
\[
v_{rn}(\theta;\alpha)
=\frac{\beta\theta-\alpha-f(\theta)}{g(\theta)},
\]
where $f$ and $g$ do not depend on $\alpha$. Let $\alpha_1<\alpha_2$
and set $\theta_1:=\theta^*_{rn}(\alpha_1)$, the smallest maximizer of
$v_{rn}(\cdot\,;\alpha_1)$. For every $\theta<\theta_1$,
\[
\begin{aligned}
&v_{rn}(\theta;\alpha_2)-v_{rn}(\theta_1;\alpha_2)\\
&\quad=
\bigl[
v_{rn}(\theta;\alpha_1)-v_{rn}(\theta_1;\alpha_1)
\bigr]
+(\alpha_2-\alpha_1)
\left(
\frac1{g(\theta_1)}-\frac1{g(\theta)}
\right)<0.
\end{aligned}
\]
Indeed, the first bracket is negative because $\theta_1$ is the
smallest maximizer and $v_{rn}(\cdot\,;\alpha_1)$ is strictly
increasing on $\{0,\dots,\theta_1\}$ by Theorem~\ref{thm:RN}(i). The
second term is negative because $g$ is strictly increasing for
$\delta<1$: the proof of Lemma~\ref{lem:monotone-increasing} gives
$e_1(g)>0$ and $e_k(g)\ge0$ for every $k$, and Newton's
forward-difference formula therefore gives $\Delta g(x)>0$ for every
$x$. Hence no $\theta<\theta_1$ maximizes
$v_{rn}(\cdot\,;\alpha_2)$, so
$\theta^*_{rn}(\alpha_2)\ge\theta_1$.

For $\delta=1$, the claim is immediate from
\[
v_{rn}(\theta)=\beta\theta-\alpha-f(\theta),
\]
by Proposition~\ref{prop:value-decomposition-and-objective}(iii).
\end{proof}

\begin{remark}[The probability parameters $p$ and $q$]\label{rem:sensitivity-pq}
Varying $p$ or $q$ changes the transition kernel itself rather than a one-step cost, so the operator comparison above does not apply. The directions reported in Section~\ref{sec:num-sensitivity}---$\theta^*_{rn}$ nonincreasing in $p$ and nondecreasing in $q$ on the computed grid (Figure~\ref{fig:sensitivity})---are therefore numerical observations; an analytic proof remains open.
\end{remark}

The next lemma records the only discount-monotonicity statement needed below.

\begin{lemma}[Discount monotonicity of the accumulation ratio]\label{lem:g_delta_decrease}
For each $x \ge 0$, the ratio $\delta \mapsto g(x+1;\delta)/g(x;\delta)$ is strictly decreasing on $(0,1]$.
\end{lemma}
\begin{proof}[Proof of Lemma~\ref{lem:g_delta_decrease}]
By Remark~\ref{rem:fg-meaning}, $1/g(x;\delta) = \mathbb{E}_0[\delta^{T_x}]$ with $T_x := \min\{t \ge 0 : X_t = x\}$. Since the chain is skip-free to the right (Lemma~\ref{lem:skip-free}), reaching $x+1$ requires first reaching $x$, so the strong Markov property gives
\[
\frac{g(x+1;\delta)}{g(x;\delta)}
= \frac{\mathbb{E}_0[\delta^{T_x}]}{\mathbb{E}_0[\delta^{T_{x+1}}]}
= \frac{1}{\mathbb{E}_x\bigl[\delta^{T_{x \to x+1}}\bigr]},
\qquad
T_{x \to x+1} := \min\{t \ge 1 : X_t = x+1\} \ \text{under } \mathbb{P}_x
\]
(the excursion analogue of $T_x$). Since $T_{x \to x+1} \ge 1$ almost surely and $\mathbb{P}_x(T_{x \to x+1} = 1) = (1-p)^x q > 0$, the map $\mathbb{E}_x[\delta^{T_{x \to x+1}}] = \sum_{k \ge 1} \delta^k\, \mathbb{P}_x(T_{x \to x+1} = k)$ is a power series in $\delta$ with nonnegative coefficients and a positive leading coefficient; it is strictly increasing on $(0,1]$, and its reciprocal is strictly decreasing.
\end{proof}

We turn to the interaction of $\beta$ and $\delta$. The $\beta$-dependent part of $v_{rn}(\theta) = \bigl(\beta\theta - \alpha - f(\theta)\bigr)/g(\theta;\delta)$ is $\beta\,\theta/g(\theta;\delta)$, so the comparative statics of $\theta^*_{rn}$ in $\beta$ are governed by the increment
\[
\Phi(\theta;\delta) \;:=\; \Delta\!\left[ \frac{\theta}{g(\theta;\delta)} \right]
= \frac{\theta+1}{g(\theta+1;\delta)} - \frac{\theta}{g(\theta;\delta)},
\]
which we call the \emph{phase function}; it depends neither on $\beta$ nor on the cost parameters.

\begin{proposition}[Sign structure of the phase function]\label{prop:beta-monotonicity}
For each $\theta \ge 1$, the following hold.
\begin{enumerate}
    \item[(i)] (Factorization) With $T_{\theta \to \theta+1}$ as above,
    \[
    \Phi(\theta;\delta)
    = \mathbb{E}_0\bigl[\delta^{T_\theta}\bigr]\, \Bigl( (\theta+1)\,\mathbb{E}_\theta\bigl[\delta^{T_{\theta \to \theta+1}}\bigr] - \theta \Bigr);
    \]
    the first factor is positive, so $\Phi(\theta;\delta)$ has the sign of the bracketed term. Moreover, $\Phi(\theta;\cdot)$ is continuous on $(0,1]$, with $\Phi(\theta;1) = 1$.
    \item[(ii)] (Low discount factors) $\Phi(\theta;\delta) < 0$ whenever $\delta < \theta/(\theta+1)$; in particular, $\Phi(\theta;\delta) < 0$ for all $\theta \ge 1$ whenever $\delta < 1/2$.
    \item[(iii)] (Phase boundary) $\Phi(\theta;\cdot)$ has exactly one zero $\delta_\theta$ in $(0,1)$, with $\Phi(\theta;\delta) < 0$ for $\delta < \delta_\theta$ and $\Phi(\theta;\delta) > 0$ for $\delta > \delta_\theta$. Moreover, $\delta_\theta \ge \theta/(\theta+1)$; in particular, $\delta_\theta \to 1$ as $\theta \to \infty$, and for every fixed $\delta < 1$ all $\theta > \delta/(1-\delta)$ have $\Phi(\theta;\delta) < 0$.
\end{enumerate}
\end{proposition}

\begin{proof}[Proof of Proposition~\ref{prop:beta-monotonicity}]
(i) By Remark~\ref{rem:fg-meaning} and the strong Markov property for the skip-free chain (Lemma~\ref{lem:skip-free}), $\mathbb{E}_0[\delta^{T_{\theta+1}}] = \mathbb{E}_0[\delta^{T_\theta}]\,\mathbb{E}_\theta[\delta^{T_{\theta \to \theta+1}}]$; substituting $1/g(x;\delta) = \mathbb{E}_0[\delta^{T_x}]$ into the definition of $\Phi$ gives the desired identity, and $\mathbb{E}_0[\delta^{T_\theta}] > 0$. Each factor is a power series in $\delta$ with nonnegative coefficients, hence continuous on $(0,1]$. As $\delta \uparrow 1$, $\mathbb{E}_\theta[\delta^{T_{\theta \to \theta+1}}] \uparrow \mathbb{P}_\theta(T_{\theta \to \theta+1} < \infty) = 1$ (the chain is positive recurrent), so $\Phi(\theta;1) = (\theta+1) - \theta = 1$.

\medskip
(ii) Since $T_{\theta \to \theta+1} \ge 1$ almost surely, $\mathbb{E}_\theta[\delta^{T_{\theta \to \theta+1}}] \le \delta$, so the bracket in part~(i) is at most $(\theta+1)\delta - \theta < 0$ whenever $\delta < \theta/(\theta+1)$; the uniform claim follows because $\theta/(\theta+1) \ge 1/2$ for all $\theta \ge 1$.

\medskip
(iii) By part~(i), $\Phi(\theta;\delta)$ has the sign of the bracket $b_\theta(\delta) := (\theta+1)\,\mathbb{E}_\theta\bigl[\delta^{T_{\theta \to \theta+1}}\bigr] - \theta$, which is strictly increasing on $(0,1]$ by Lemma~\ref{lem:g_delta_decrease}, with $b_\theta(\delta) < 0$ for $\delta < \theta/(\theta+1)$ by part~(ii) and $b_\theta(1) = 1 > 0$ by part~(i). Hence $\Phi(\theta;\cdot)$ has exactly one zero $\delta_\theta \in [\theta/(\theta+1), 1)$, and the sign claims follow from strict monotonicity; the final claims follow because $\delta < \theta/(\theta+1)$ is equivalent to $\theta > \delta/(1-\delta)$.
\end{proof}

\begin{corollary}[Comparative statics in the per-item reward]\label{cor:beta-statics}
Fix $(c, \gamma, \alpha)$ and $\delta$, and write $\theta^*_{rn}(\beta)$ for the optimal threshold of Theorem~\ref{thm:RN} as a function of $\beta$.
\begin{enumerate}
    \item[(i)] If $\Phi(\theta;\delta) \le 0$ for all $\theta \ge \theta^*_{rn}(\beta_1)$---in particular, if $\delta < 1/2$---then $\theta^*_{rn}(\beta_2) \le \theta^*_{rn}(\beta_1)$ for all $\beta_2 > \beta_1$.
    \item[(ii)] If $\Phi(\theta;\delta) \ge 0$ for all $\theta \le \theta^*_{rn}(\beta_1)$, then $\theta^*_{rn}(\beta_2) \ge \theta^*_{rn}(\beta_1)$ for all $\beta_2 > \beta_1$.
\end{enumerate}
\end{corollary}

\begin{proof}[Proof of Corollary~\ref{cor:beta-statics}]
The reduced objective is affine in $\beta$: $v_{rn}(\theta;\beta) = \bigl(\beta\theta - \alpha - f(\theta)\bigr)/g(\theta;\delta)$ (with $g \equiv 1$ when $\delta = 1$), with $f$ and $g$ independent of $\beta$ (Proposition~\ref{prop:value-decomposition-and-objective}). Hence, for $\beta_2 > \beta_1$ and $\theta_1 := \theta^*_{rn}(\beta_1)$,
\[
v_{rn}(\theta;\beta_2) - v_{rn}(\theta_1;\beta_2) = \bigl[ v_{rn}(\theta;\beta_1) - v_{rn}(\theta_1;\beta_1) \bigr] + (\beta_2 - \beta_1)\Bigl( \frac{\theta}{g(\theta;\delta)} - \frac{\theta_1}{g(\theta_1;\delta)} \Bigr).
\]

(i) For any $\theta > \theta_1$, the first bracket is nonpositive because $\theta_1$ maximizes $v_{rn}(\cdot\,;\beta_1)$, and the second term is nonpositive because $\theta \mapsto \theta/g(\theta;\delta)$ is nonincreasing on $\{\theta_1, \theta_1+1, \dots\}$ by hypothesis. Hence the smallest maximizer of $v_{rn}(\cdot\,;\beta_2)$ is at most $\theta_1$.

(ii) For any $\theta < \theta_1$, the first bracket is strictly negative because $\theta_1$ is the smallest maximizer and $v_{rn}(\cdot\,;\beta_1)$ is strictly increasing on $\{0, \dots, \theta_1\}$ (Theorem~\ref{thm:RN}(i)), and the second term is nonpositive because $\theta \mapsto \theta/g(\theta;\delta)$ is nondecreasing on $\{0, \dots, \theta_1\}$ by hypothesis. Hence no $\theta < \theta_1$ maximizes $v_{rn}(\cdot\,;\beta_2)$, so $\theta^*_{rn}(\beta_2) \ge \theta_1$.
\end{proof}

\begin{remark}
At $\delta = 1$, $\Phi(\theta;1) = 1$ (Proposition~\ref{prop:beta-monotonicity}(i)), so Corollary~\ref{cor:beta-statics}(ii) applies without restriction: $\theta^*_{rn}$ is nondecreasing in $\beta$ when $\delta = 1$.
\end{remark}

\subsection{Reward--Discount Phase Transition}\label{sec:num-phase}
By Corollary~\ref{cor:beta-statics}, the sign of $\Phi(\theta;\delta)$ determines how $\theta^*_{rn}$ responds to $\beta$: nonincreasing where $\Phi < 0$, nondecreasing where $\Phi > 0$. Since $\Phi$ depends on neither $\beta$ nor the cost parameters, the $(\theta, \delta)$ plane is partitioned by its zero set; by Proposition~\ref{prop:beta-monotonicity}(iii) the zero $\delta_\theta$ is unique, lies in $[\theta/(\theta+1), 1)$, and converges to $1$. In all our experiments, $\delta_\theta$ is increasing in $\theta$, so $\Phi(\theta;\delta) > 0$ exactly for $\theta < \bar\theta(\delta) := \min\{\theta : \delta_\theta \ge \delta\}$, and the conditions of Corollary~\ref{cor:beta-statics} reduce to comparing $\theta^*_{rn}$ with $\bar\theta(\delta)$.

\begin{figure}[htbp]
\centering
\includegraphics[width=\textwidth]{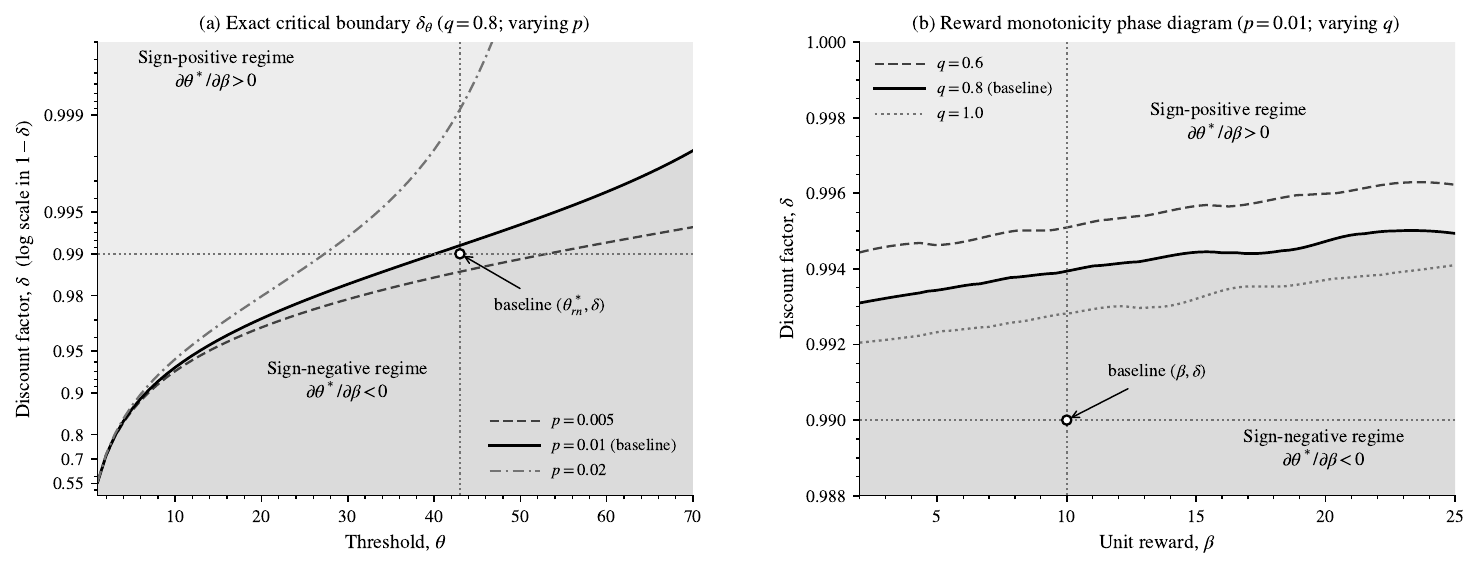}
\caption{The reward--discount phase boundary. Panel~(a): the critical boundary $\delta_\theta$ of Proposition~\ref{prop:beta-monotonicity}(iii) as a function of $\theta$ (log scale in $1-\delta$) for $q = 0.8$ and $p \in \{0.005, 0.01, 0.02\}$; $\Phi$ is positive above the boundary ($\theta^*_{rn}$ nondecreasing in $\beta$) and negative below ($\theta^*_{rn}$ nonincreasing in $\beta$). Panel~(b): the boundary in the $(\beta, \delta)$ plane for $p = 0.01$ and $q \in \{0.6, 0.8, 1.0\}$, with the sign of $\Phi$ evaluated at $\theta^*_{rn}(\beta)$. In panels (a) and (b), the markers show the baseline values $(\theta^*_{rn}, \delta) = (43, 0.99)$ and $(\beta, \delta) = (10, 0.99)$, respectively, just inside the sign-negative region.}
\label{fig:phase_diagram}
\end{figure}

Figure~\ref{fig:phase_diagram}(a) displays the boundary $\delta_\theta$ (log scale in $1-\delta$) for $q = 0.8$ and $p \in \{0.005, 0.01, 0.02\}$. It rises toward $1$ as $\theta$ grows, in line with $\delta_\theta \ge \theta/(\theta+1)$ ($\delta_1 \approx 0.56$, versus $1/2$; $\delta_{43} = 0.9865, 0.9913, 0.9991$ for the three values of $p$, versus $43/44 \approx 0.977$), and moves up as $p$ increases. Economically, $\beta\,\Phi(\theta;\delta)$ is the marginal discounted-reward content of one additional item: in the sign-negative region, waiting for a larger batch discounts the per-item reward more heavily than the batch size grows, so a higher $\beta$ weakly lowers the threshold.

Figure~\ref{fig:phase_diagram}(b) redraws the boundary in the $(\beta, \delta)$ plane for $p = 0.01$ and $q \in \{0.6, 0.8, 1.0\}$ ($q = 1$ being the deterministic-arrival limit), tracing where $\theta^*_{rn}(\beta)$---recomputed at each grid point---coincides with $\bar\theta(\delta)$. The boundary decreases slightly as $q$ rises, with visibly smaller sensitivity to $q$ than to $p$ in panel~(a).

\end{document}